\documentclass[11pt, reqno]{article}
\pdfoutput=1
\usepackage{amssymb,amsmath,amsthm,amsfonts,latexsym}
\usepackage{epsfig,esint}
\usepackage{mathrsfs}
\usepackage{mathtools,textcomp}
\usepackage{graphicx}
\usepackage{epic,eepic,wrapfig,color, ifthen} 
\usepackage{url,cite,geometry}
\usepackage[colorlinks=true, citecolor=blue, filecolor=cyan, linkcolor=red, urlcolor=magenta]{hyperref}
\usepackage{pmboxdraw}
\usepackage{verbatim}
\usepackage[normalem]{ulem}
\usepackage{enumerate}
\usepackage[bottom]{footmisc}

\def\nn{\nonumber}

\def\Dini{{\mathscr D_0}}
\def\DDini{{\mathscr D_{00}}}
  
 \def\E {{\mathbb E}}

\def\P {{\mathbb P}}  
\def\bS {{\mathbb S}}

\def\sA {{\mathcal A}}
\def\sL{{\mathcal L}}

\def\1{{\rm{1}}}
\def\nn{\nonumber}

\def\qed{{\hfill $\Box$ \bigskip}}
\def\eps{\varepsilon}
\def\wt{\widetilde}
\def\wh{\widehat}

\def\ind {{\bf 1}}

\def\sA {{\cal A}} \def\sB {{\cal B}} \def\sC {{\cal C}}

\def\sE {{\cal E}} 
\def\sF {{\cal F}}
  
 \def\sK {{\cal K}}
\def\sL {{\cal L}}

\def\R {{\mathbb R}} 
 
\def\E {{\mathbb E}}  \def \P{{\mathbb P}}

\def \diag{{\textrm{\rm diag}}}
\def \diam{{\textrm{\rm diam}}}

\numberwithin{equation}{section}
\def\qed{{\hfill $\Box$ \bigskip}}

\def\eps{\varepsilon}
\def\wh{\widehat}
\def\wt{\widetilde}
\def\pf{\noindent{\bf Proof. }}

\def\vp{{\varphi}}

\theoremstyle{plain}
\newtheorem{thm}{Theorem}[section]
\newtheorem{lemma}[thm]{Lemma}
\newtheorem{cor}[thm]{Corollary}
\newtheorem{remark}[thm]{Remark}
\newtheorem{prop}[thm]{Proposition}
\newtheorem{defn}[thm]{Definition}
\newtheorem{definition}[thm]{Definition}

\newtheorem{example}[thm]{Example}

\theoremstyle{definition}
\newtheorem*{eg*}{Example}
\newtheorem*{egs*}{Examples}
\newtheorem*{def*}{Definition}
\theoremstyle{remark}
\begin{document}
	\title{Approximate 	factorizations
	for  non-symmetric jump processes}

\author{
	{\bf Soobin Cho}
	\quad {\bf Renming Song\thanks{Research supported in part by a grant from the Simons Foundation (\#960480, Renming Song)}
	}
}

\date{}

	\maketitle
	
	\begin{abstract}
		In this paper, we first extend the approximate factorization for purely discontinuous Markov process established in \cite{CKSV20} by getting rid of some of the conditions imposed in \cite{CKSV20}. Then we apply the approximate factorization to obtain 
		sharp two-sided 
		heat kernel estimates for  three classes of processes: stable-like processes with critical killings in $C^{1, {\rm Dini}}$
		open sets; killed stable-like processes in the setting of  \cite{KW24} in $C^{1, \eps
			}$ open sets; and non-symmetric stable processes in what we call $C^{1,2{\text - \rm Dini}}$ open sets.  In particular, we obtain explicit sharp two-sided heat kernel estimates of  
		killed $\alpha$-stable processes in $C^{1, {\rm Dini}}$ open sets for all $\alpha\in (0, 2)$ and of censored $\alpha$-stable processes in $C^{1, {\rm Dini}}$ open sets for all $\alpha\in (1, 2)$.
		
		\medskip
		
		\medskip
		
		
	\end{abstract}
	\allowdisplaybreaks

	\tableofcontents
	\section{Introduction}\label{s:intro}

	An isotropic $\alpha$-stable process, $\alpha\in (0, 2)$, on $\R^d$ is an $\R^d$-valued L\'evy process $(Y_t, \P_x)$ such that
	$$
	\E_x\left[ e^{i\xi\cdot (Y_t-Y_0)}\right]=e^{-t|\xi|^\alpha}, \quad \text{ for all } t>0,\, x, \xi\in \R^d.
	$$
	The infinitesimal generator of $(Y_t, \P_x)$ 	is
	the fractional Laplacian $\Delta^{\alpha/2}:=-(-\Delta)^{\alpha/2}$, which is a prototype of non-local operators. The fractional Laplacian can be written in the form 
	$$
	\Delta^{\alpha/2}u(x)=\lim_{\varepsilon\downarrow 0}\int_{\{y\in \R^d: |y-x|>\varepsilon\}}(u(y)-u(x))\frac{A}{|x-y|^{d+\alpha}}dy,
	$$
	where
	$$
	A:=A(d, \alpha):=\frac{\alpha\Gamma((d+\alpha)/2)}{2^{1-\alpha}\pi^{d/2}\Gamma(1-\alpha/2)}.
	$$
	It is a classical result that the $\alpha$-stable process $(Y_t, \P_x)$ admits a transition density $q(t, x, y)$ with respect to the Lebesgue measure and that $q(t, x, y)$ satisfies the following two-sided estimates:
	there exists $C>1$ such that 
	$$
	C^{-1}\bigg( t^{-d/\alpha} \wedge \frac{t}{|x-y|^{d+\alpha}}\bigg)\le q(t, x, y)\le C\bigg( t^{-d/\alpha} \wedge \frac{t}{|x-y|^{d+\alpha}}\bigg) \quad \text{ for all } t>0,\, x, y\in \R^d.
	$$
	Analytically, $q(t, x, y)$ is the heat kernel of the fractional Laplacian $\Delta^{\alpha/2}$.
	For any open set $D\subset \R^d$, let $\tau_D:=\inf\{t>0: Y_t\notin D\}$. Define a  process $Y^D$ by $Y^D_t=Y_t$ if $t<\tau_D$ and $Y^D_t=\partial$ if 	$t\ge \tau_D$, where $\partial$ is a cemetery point added to $D$. 
	$Y^D$ is called the 	killed process of $Y$ in $D$. 
	The generator of $Y^D$ is the following 
	regional fractional Laplacian --
	killed
	fractional Laplacian:
	$$
	\Delta^{\alpha/2}|_{k, D}u(x)
	:=\lim_{\varepsilon\downarrow 0}\int_{\{y\in D: |y-x|>\varepsilon\}}(u(y)-u(x))\frac{A}{|x-y|^{d+\alpha}}dy
	+\kappa_Du(x), \quad \;\;u\in C^2_c(D), \, x\in D,
	$$
	where
	$$
	\kappa_D(x)=\int_{D^c}\frac{A}{|x-y|^{d+\alpha}}dy.
	$$
	The process $Y^D$
	also admits a transition density with respect to the Lebesgue measure and this transition density is also the heat kernel of the killed
	fractional Laplacian $\Delta^{\alpha/2}|_{k, D}$. 
	
	In \cite{BBC03}, Bogdan, Burdzy and Chen introduced the so-called censored $\alpha$-stable process in $D$ by suppressing all the jumps to $D^c$.  The generator of the censored $\alpha$-stable process in $D$ is the following regional fractional Laplacian --
	censored
	fractional Laplacian:
	$$
	\Delta^{\alpha/2}|_{c, D}u(x)
	:=\lim_{\varepsilon\downarrow 0}\int_{\{y\in D: |y-x|>\varepsilon\}}(u(y)-u(x))\frac{A}{|x-y|^{d+\alpha}}dy, 
	\quad \;\;u\in C^2_c(D), \, x\in D.
	$$
	The censored $\alpha$-stable process admits a transition density with respect to the Lebesgue measure and this transition density is also the heat kernel of the censored
	fractional Laplacian $\Delta^{\alpha/2}|_{c, D}$.

	The killed fractional Laplacian and the censored
	fractional Laplacian are two prototypes of regional non-local operators and they are important in analysis and partial differential equations. Unlike the fractional Laplacian (in the whole space), sharp two-sided heat kernel estimates for the killed fractional Laplacian and the censored
	fractional Laplacian are much more difficult to establish. Explicit sharp two-sided heat kernel estimates for the killed fractional Laplacian in $C^{1, 1}$ open set were established by Chen, Kim and Song \cite{CKS10a}. An approximate factorization for the heat kernel of  the killed fractional Laplacian in $\kappa$-fat open sets were obtained 	by Bogdan, Grzywny and Ryznar in \cite{BGR10}. 
	Explicit sharp two-sided heat kernel estimates for the censored fractional Laplacian in $C^{1, 1}$ open set were established by Chen, Kim and Song \cite{CKS10}. The fact that $C^{1, 1}$ open sets satisfy the interior and exterior ball conditions played a crucial role in obtaining the explicit heat kernel estimates in \cite{CKS10a, CKS10}. 
	
	Subsequently, efforts have been made to extend the explicit sharp two-sided heat kernel estimates to open sets with less regularity and to general operators of the form
	$$
	Lu(x)=\lim_{\varepsilon\downarrow 0}\int_{\{y\in \R^d: |y-x|>\varepsilon\}}(u(y)-u(x))\frac{K(x, y)}{|x-y|^{d+\alpha}}dy
	$$
	with $K(x, y)$ a symmetric function 
	bounded between 
	two positive constants and satisfying certain conditions.
	The related topic on the regularity up to the boundary of the solution to
	\begin{align}\label{e:Poissoneq}
		\begin{cases}
			Lu=f & \text{ in } D,\\
			u=0, & \text{ in } \R^d\setminus D
		\end{cases}
	\end{align}
	was also intensively studied.
	Explicit sharp two-sided heat kernel estimates for $\Delta^{\alpha/2}|_{k, D}$ were extended to 	$C^{1, \varepsilon}$ open sets with $\varepsilon>\alpha/2$ 
	in \cite{KK, KK2}. Explicit sharp two-sided heat kernel estimates for $\Delta^{\alpha/2}|_{k, D}$ and for $\Delta^{\alpha/2}|_{c, D}$ were extended to 	$C^{1, \varepsilon}$ open sets with $\varepsilon> (\alpha-1)_+$ in \cite{SWW}. 
	In the recent paper \cite{KW24}, Kim and Weidner established explicit sharp two-sided estimates for the Green function of killed stable-like operators in divergence form in bounded $C^{1, \varepsilon}$ open sets for any $\varepsilon>0$. In \cite{KW24}, only a very weak H\"{o}lder continuous condition is assumed for 	the coefficient $K$, see \eqref{e:low-regular-2} below. 
	It was mentioned in \cite{KW24} that it was an intriguing problem to establish Dirichlet heat kernel estimates in the setting of  \cite{KW24}. On the regularity up to the boundary of solutions to \eqref{e:Poissoneq},  Ros-Oton and Serra \cite{ROS} proved that, in the case when $L=\Delta^{\alpha/2}$ and $D$ is  a bounded $C^{1, 1}$, if $f$ is bounded, then the solution $u$ to \eqref{e:Poissoneq} is $C^{\alpha/2}(\overline{D})$.  
	This result was extended to more the case when $K$ is translation invariant in  
	\cite{f-Gru, Gru, ROS16a, ROS16b, ROS17, AROS20}. 
	In \cite{ROW}, Ros-Oton and Weidner proved that, when $L$ is translation invariant and $D$ is 	$C^{1, \varepsilon}$ for any $\varepsilon>0$, 
	the solution $u$ to \eqref{e:Poissoneq} is $C^{\alpha/2}(\overline{D})$. In \cite{KW24}, Kim and Weidner extended this result to case when
	$K$ is not necessarily translation invariant.  Boundary regularity of the solution to \eqref{e:Poissoneq} in the case when $L$ is the generator of a non-symmetric stable process was studied in \cite{DRSV22, Ju}.
	
	In this paper, we first extend the approximate factorization for purely discontinuous Markov process established in \cite{CKSV20} by getting rid of Assumptions {\bf A} and {\bf U} imposed in \cite{CKSV20}, and 
	show that, under some natural conditions, the approximate factorization of the heat kernel is equivalent to the approximate factorization of the Green function. Then we apply the approximate factorization to obtain  	sharp two-sided heat kernel estimates for  
	three classes of processes: stable-like process with critical killings in $C^{1, {\rm Dini}}$
	open sets; killed stable-like processes in the setting of  \cite{KW24}; and non-symmetric stable processes in $C^{1, 2{\text -\rm Dini}}$
	open sets. See Definition \ref{d:Dini} below for the definitions of $C^{1, {\rm Dini}}$ and $C^{1, 2{\text -\rm Dini}}$ open sets. In particular, we obtain explicit sharp two-sided heat kernel estimates of 
	$\Delta^{\alpha/2}|_{k, D}$ and for $\Delta^{\alpha/2}|_{c, D}$ for  $C^{1, {\rm Dini}}$ open sets. 
	In \cite{CS25}, we will show that $C^{1, {\rm Dini}}$ regularity assumption on $D$ for sharp two-sided heat kernel estimates of 
	$\Delta^{\alpha/2}|_{k, D}$ is optimal.
	
	This paper is organized as follows. In Section \ref{s:setp-prel}, we introduce the general setup for the approximate factorization and present some preliminaries. In Section \ref{s:surv}, we establish  survival estimates for general Markov processes in our class.
	The main results on approximate factorization are stated and proved in Section \ref{s:fact}.  In Section \ref{s:appl}, we apply the approximate factorizations to obtain sharp two-sided heat kernel estimates for three classes of processes. Subsection \ref{ss:5.1} deals with stable-like processes with critical killings in $C^{1, {\rm Dini}}$ open sets.  In this subsection, the coefficient is assumed to satisfy a H\"{o}lder continuity condition near the diagonal (see \eqref{e:censored-2}). 
	Subsection \ref{ss:5.2} deals with non-symmetric stable processes in $C^{1, 2{\text -\rm Dini}}$ open sets. A key step in both Subsections   \ref{ss:5.1} and  \ref{ss:5.2} is the construction of appropriate barrier functions. Subsection \ref{ss:5.3} deals with explicit Dirichlet heat kernel estimates in the setting of \cite{KW24}. 
	The Appendix collects some auxiliary results on Dini and double Dini functions.
	
	\noindent \textit{Notation:} 
Throughout this paper,   lower case letters  $c$, $c_i$, $i=0,1,2,...$, denote constants in the proofs, their values remain fixed in each proof and their labelling starts anew in each proof unless otherwise specified to represent particular values. Upper case letters $C$, $C_i$, $i=0,1,2,...$, are used for constants in the statements of results
and their labelling starts anew in each result.
 We write $\R^d_+:=\left\{y=(\wt y, y_d)\in \R^d: y_d>0\right\}$ for the upper half-space.
We use the notations $a\wedge b:=\min\{a,b\}$ and $a\vee b:=\max\{a,b\}$ for $a,b\in \R$.  The notation $f(x) \asymp g(x)$ means that there exist constants $c_2 \ge c_1>0$ such that $c_1g(x)\leq f (x)\leq c_2 g(x)$ for a specified range of the variable $x$. 
	
\section{General
	setup and preliminaries}\label{s:setp-prel} 

\subsection{General setup}

	Let	$(M,d)$ be a locally compact separable metric space  and $\mu$ be a  Radon measure on $M$ with full support.	Let $B(x,r)=\{y\in M: d(x,y)<r\}$ be the open ball of radius $r$ centered at $x\in M$ and   $V(x,r):= \mu(B(x,r))$.  We assume that $(M,d,\mu)$ satisfies the volume doubling and reverse volume doubling conditions: there exist constants 
	$\beta_2\ge \beta_1>0$ and 
	$A\ge 1$ such that
	\begin{equation}\label{e:VD}
	\frac{V(x,R)}{V(x,r)}\le  
	A\left(\frac{R}{r}\right)^{\beta_2}
		\quad \text{for all}\  x\in M \text{ and }\ 0<r\le R
	\end{equation}
	and
		\begin{equation}\label{e:RVD}
	\frac{V(x,R)}{V(x,r)} \ge 		
	A^{-1}\left(\frac{R}{r}\right)^{\beta_1}
		\quad \text{for all}\  x\in M \text{ and }\ 0<r\le R< \diam(M),
	\end{equation}
where $\diam(M)$ is the diameter of $M$.	The reverse volume doubling condition  \eqref{e:RVD} implies the existence of a constant $k_0>1$ such that $V(x,r) \ge 2 V(x,r/k_0)$ for all $x\in M$ and $r\in (0,\diam(M))$. In particular, it holds that 
	\begin{align}\label{e:uniformlyperfect}
		B(x,r) \setminus B(x, r/k_0) \neq \emptyset \quad \text{for all $x\in D$ and $r\in (0, \diam(M))$.}
	\end{align}

	Let $Y=(Y_t, \P_x, \zeta^Y)$ be a Hunt process on $M$. Here $\zeta^Y$ is the lifetime of $Y$, so that $Y_t=\partial$ for $t\ge \zeta^Y$, where $\partial$ is a cemetery point added to $M$. We extend any function $f$ on $M$ to $M\cup \{\partial\}$ by setting $f(\partial)=0$. The transition semigroup 
	$(Q_t)_{t\ge 0}$ of $Y$ is defined by
	\begin{align*}
		Q_t f(x) = \E_x[f(Y_t)], \quad t\ge 0, \; x \in M,
	\end{align*}
whenever the expectation makes sense.	We assume that 	the semigroup $(Q_t)_{t\ge 0}$   admits a strictly positive and jointly continuous transition density $q(t,x,y)$ with respect to $\mu$.

	Let $\phi:(0,\infty) \to (0,\infty)$ be a strictly increasing continuous function satisfying the following weak scaling condition: there exist constants 	$\alpha_2 \ge \alpha_1>0$ and 
	$A\ge 1$
	such that
	\begin{align}\label{e:scale-phi}
		A^{-1}\left(\frac{R}{r}\right)^{\alpha_1} \le \frac{\phi(R)}{\phi(r)} \le A
		\left(\frac{R}{r}\right)^{\alpha_2} 
		\quad  \text{for all $0<r\le R$}.		
	\end{align}
	Define
	\begin{align*}
		& \wt q(t,x,y)  := \frac{1}{V(x,\phi^{-1}(t))} \wedge \frac{t}{V(x, d(x,y)) \phi(d(x,y))}, \quad t>0, \ x,y\in M.
	\end{align*}
	
Throughout this paper,	we assume that  $Y$ satisfies the following:

	\medskip
	
	\setlength{\leftskip}{0.2in}	
	
	\noindent	{\bf (A)}:  \textit{There exist  $T_0 \in (0, \infty]$ and 
	$A_0\ge 1$ such that
	\begin{align*}
		A_0^{-1} \wt q(t,x,y) \le q(t,x,y) \le A_0 \wt
		q(t,x,y), \qquad (t,x,y) \in (0,T_0)\times M\times M.
	\end{align*}
	When $T_0=\infty$, we also assume that $M$ is unbounded.}

	\begin{remark}\label{r:extend} \rm 
		When {\bf (A)} holds with  $T_0<\infty$, the exact value of $T_0$ is unimportant. In fact, {\bf (A)} implies that for any  $T\in (0,\infty)$, there exists a constant   		$A=A(T)\ge 1$ such that
		\begin{align*}
			A^{-1} \wt q(t,x,y) \le q(t,x,y) \le  A \wt
			q(t,x,y), \qquad (t,x,y) \in (0, T)\times M \times M.
		\end{align*}
		See \cite[Remark 2.3]{CKSV20} for a proof. 
	\end{remark}

	\setlength{\leftskip}{0in}	
	
	As an immediate consequence of {\bf (A)}, $Y$ is Feller and strongly Feller. By general theory, 
	 there exists a  quasi-regular (possibly non-symmetric) Dirichlet form $(\sE^Y, \sF^Y)$ which is properly associated with  $Y$. 
	 Moreover, $Y$ admits a dual Hunt process $\wh Y = (\wh Y, \wh \P_x)$ with respect to $\mu$.
	 That is, the transition semigroup $(\wh Q_t)_{t\ge 0}$ of $\wh Y$ is the dual of  $(Q_t)_{t\ge 0}$ with respect to $\mu$. See \cite[Theorem IV.5.1, Corollary IV.5.2 and Remark IV.5.3(i)]{MR92}. 	 The dual process $\wh Y$ is also Feller and strongly Feller, and has a  transition density given by $\wh q(t,x,y ):=q(t,y,x)$. 
	In the sequel, all objects related to the dual process $\wh{Y}$ will be denoted by a hat.

	 According to  \cite{BJ73},
	  the Hunt process $Y$ admits a \textit{L\'evy system} $(J(x,dy),H_t)$. That is,  there exist a family of  $\sigma$-finite Borel measures  
           $J(x,dy)$  on  $M$ with $J(x,\{x\})=0$ for all $x\in M$, and a positive continuous additive functional $(H_t)_{t\ge 0}$ of $Y$ such that
for  any stopping time $\sigma$ and   any  Borel function $f: M\times M \to [0,\infty]$ vanishing on the diagonal, the following identity holds for all $x\in M$:
	\begin{equation}\label{e:levy-system}
		\E_x\sum_{s\in (0, \sigma]}f(Y_{s-}, Y_s)= \E_x\int^{\sigma}_0\int_{M}f(Y_s, y)J(Y_s, dy)d H_s.
	\end{equation}
	 Since $Y$ has a strictly positive transition density,  by \cite[Proposition 5.1]{Ba79},  $Y$ admits a L\'evy system of the form $(J(x,dy), t)$.  	 Throughout this paper, we assume  without loss of generality that $Y$ has a L\'evy system  $(J(x,dy), t)$.	The kernel $J(x,dy)$ is called the \textit{jump kernel} 	of $Y$.
	
		The Le\'vy system formula \eqref{e:levy-system} implies that for any bounded continuous function $f$ on $M$ with compact support $K$ and $x\in M\setminus K $,
	\begin{equation*}
		\int_{M}f(y)J(x,dy)=\lim_{t\to 0} \frac{Q_t f(x)}{t}.
	\end{equation*} 
	Thus, by {\bf (A)},  the jump kernel has a density $J(x,dy) = J(x,y) \mu(dy)$ and 
	\begin{align}\label{e:J-density}
	\frac{A^{-1}}{V(x, d(x,y)) \phi(d(x,y))} \le J(x,y) \le \frac{A}{V(x, d(x,y)) \phi(d(x,y))},
	\end{align}
	for some $A\ge 1$.
	Moreover, by \cite[Theorem 3.5]{Ge71},   $\wh Y$ has a jump kernel $\wh J(x,dy)=\wh J(x, y)\mu(dy)$ with $\wh J(x, y)=J(y, x)$.

	For an open set $D\subset M$, define the first exit times from $D$ by  $\tau^Y_D:=\inf\{t>0:\, Y_t\notin D\}$ and  $\wh\tau^Y_D:=\inf\{t>0:\, \wh Y_t\notin D\}$. 
	The killed processes $Y^D$ and $\wh Y^D$ on $D$ are defined by
	\begin{align*}
		Y_t^D := \begin{cases} 
			Y_t, &t<\tau^Y_D,\\
			\partial, & t\ge \tau^Y_D,
			\end{cases} \quad \text{and} \quad 	\wh Y_t^D := \begin{cases} 
			\wh Y_t, &t<\wh \tau^Y_D,\\
			\partial, & t\ge \wh \tau^Y_D.
			\end{cases} 
	\end{align*}
	 
	A measure $\kappa$ on an open set $D\subset M$ is 
	said to be a \textit{smooth measure}
	of $Y^D$ with respect to the reference measure $\mu|_D$,
	if there is a positive continuous additive functional
	$A=(A_t)_{t\ge 0}$ of $Y^D$ such that for any
	bounded non-negative Borel function $f$ on $D$,
	\begin{equation*}
		\int_D f(x) \kappa (dx) = \lim_{t\downarrow 0}
		\E_{\mu} \left[ \frac1t \int_0^t f(Y^D_s) 
		dA_s \right].
	\end{equation*}
	For 	a smooth Radon measure $\kappa$ of $Y^D$ and  $a \ge 0$, we define
	\begin{align*}
		N^{D, \kappa}_a(t) := \sup_{x\in M}
		\int_0^t \int_{z\in D: \delta_D(z) >a \phi^{-1}(t)}  \wt q(s,x,z)\kappa(dz)ds, \quad t>0.
	\end{align*}
	Here $\delta_D(x)$ denotes the distance between $x$ and $D^c$.

	\begin{definition}\label{d:KT}
	\rm	Let $D\subset M$ be an open set and $\kappa$ be a non-negative smooth measure for  both $Y^D$ and $\wh Y^D$ with respect to the reference measure $\mu|_D$ and let
		$T\in \{1,\infty\}$.
		The measure $\kappa$ is said to be in the  class
		$\mathbf{K}_T(D)$ if
		
		\smallskip
		
		(1) $\displaystyle \sup_{t \in (0,T)}N^{D, \kappa}_a(t)<\infty$ for all $a \in (0, 1]$;

		(2)
		$\displaystyle\lim_{t \to 0}N^{U, \kappa}_0(t)=0$ for every relatively compact open subset $U$ of $D$.
	\end{definition}
	
	\begin{remark}
	\rm	If $\kappa\in\mathbf{K}_1(D)$, then $ \sup_{t\le T}N^{D, \kappa}_a(t)<\infty$ for all $a>0$ and $T\in (0,\infty)$. See \cite[Remark 2.13]{CKSV20}.
	\end{remark}
	
	We recall some standard examples of elements in 	$\mathbf{K}_T(D)$ from \cite[Examples 2.17 and 2.18]{CKSV20}.
	
	\begin{example}\label{ex:critical-killing}
	\rm	(i)	Let $\kappa(dx) =  \kappa(x)\mu(dx) $
		with $ 0\le \kappa(x) \le C/\phi( \delta_D(x) \wedge 1)$ for some $C>0$.  Then
		$\kappa \in  \mathbf{K}_1(D)$.
		
		\noindent (ii) 	Suppose $T_0=\infty$  and $\kappa(dx) =  \kappa(x)\mu(dx) $
		with $ 0\le \kappa(x) \le C/\phi( \delta_D(x))$ for some $C>0$.  Then
		$\kappa \in  \mathbf{K}_\infty(D)$.
	\end{example}

	Let us fix an open set $D\subset M$ and let $\kappa \in \mathbf K_1(D)$. Denote by
$(A^{\kappa}_t)_{t\ge0}$  the positive continuous additive functional  of $Y^D$
	with Revuz measure $\kappa$ and  by
	$(\wh A^{\kappa}_t)_{t\ge0}$  the
	positive continuous additive functional of
	$\wh Y^D$
	with Revuz measure $\kappa$.
	For any non-negative Borel function $f$ on $D$, we define
	\begin{align*}
		P^{\kappa}_t f(x) = \E_x\left[ \exp (-A^{\kappa}_t) f(Y^D_t)\right], \quad \wh P^{\kappa}_t f(x) = \wh \E_x\left[ \exp (- \wh A^{\kappa}_t) f(\wh Y^D_t)\right],
		\quad t\ge 0 , x\in D.
	\end{align*}
	The semigroup $(P^{\kappa}_t)_{t\ge 0}$ (respectively $(\wh P^{\kappa}_t)_{t\ge 0}$)
	is called the \textit{Feynman-Kac semigroup} of $Y^D$ (respectively $\wh Y^D$)
	associated with $\kappa$.
	By \cite[Theorem 6.10(2)]{Y96},
	$P^{\kappa}_t$ and $\wh P^{\kappa}_t$ are duals of each other
	with respect to the measure $\mu|_D$.

	Let $X^{\kappa}$ ($\wh X^{\kappa}$, respectively) be a Hunt process on $D$ corresponding to the  transition semigroup $(P_t^{\kappa})_{t\ge 0}$ ($(\wh P_t^{\kappa})_{t\ge 0}$, respectively).
	It follows from \cite[Theorem I.3.4]{Y96}
	that the L\'evy system of $X^\kappa$ ($\wh X^{\kappa}$, respectively) 
	restricted to $D\times D$ 
	is the same as that  of $Y$ ($\wh Y$, respectively), hence the following L\'evy system formulas are valid: 	for any $x\in D$,  any stopping times $\sigma$ for  $X^\kappa$ and $\wh \sigma$  for $\wh X^{\kappa}$,  and any Borel function $f: D\times D \to [0,\infty]$ vanishing on the diagonal,  
	\begin{align}
		\E_x\sum_{s\in (0, \sigma]}f(X^{\kappa}_{s-}, X^{\kappa}_s)= \E_x\int^\sigma_0\int_{D}f(X^{\kappa}_s, y)J(X^{\kappa}_s, y)\mu(dy)ds,\label{e:levy_systemX}\\
		\wh\E_x\sum_{s\in (0, \wh\sigma]}f(\wh X^{\kappa}_{s-}, \wh X^{\kappa}_s)= \E_x\int^{\wh \sigma}_0\int_{D}f(\wh X^{\kappa}_s, y)J(y,\wh X^{\kappa}_s)\mu(dy)ds.	\label{e:levy_systemX-dual}
	\end{align}

\subsection{Preliminaries}

By a standard argument, from \eqref{e:VD} and \eqref{e:scale-phi}, we can obtain
that for all $x\in M$ and $r>0$,
\begin{align}\label{e:volumephi-integral}
\int_{B(x,r)^c} \frac{\mu(dy)}{V(y,d(x,y)) \phi(d(x,y))} \asymp 	\int_{B(x,r)^c} \frac{\mu(dy)}{V(x,d(x,y)) \phi(d(x,y))}  	 \le  \frac{C}{\phi(r)}.
\end{align}
 See \cite[Lemma 2.1]{CKW-memo} for a proof. Consequently, by \eqref{e:J-density}, we have
\begin{align}\label{e:jump-integral}
	\int_{B(x,r)^c} J(x,y) \mu(dy)  \vee 	\int_{B(x,r)^c} J(y,x) \mu(dy)\le  \frac{C}{\phi(r)} \quad \text{for all $x\in M$ and $r>0$}.
\end{align}

\begin{lemma}\label{l:EP}
	There exists 	$C=C(\phi,A_0,M)>0$ 
	such that for all $x\in M$, $r>0$ and $t\in (0,T_0)$,
	\begin{align*}
		\P_x\big( \tau^Y_{B(x,r)} < t \wedge \zeta^Y  \big) \vee 	\wh\P_x\big( \wh\tau^Y_{B(x,r)} < t \wedge \wh\zeta^Y  \big) \le \frac{Ct}{\phi(r)}.
	\end{align*}
\end{lemma}
\pf Let $x\in M$ and $r>0$. For any $s\in (0,T_0)$, by the strong Markov property, we have
\begin{align}\label{e:EP-1}
		&\P_x\big( \tau^Y_{B(x,r)} < s \wedge \zeta^Y, \, s/2 < \zeta^Y  \big) \nn\\
		&\le 	\P_x\big( Y_{s/2} \in M\setminus B(x,r/2)  \big)  + 	\P_x\big( Y_{s/2} \in  B(x,r/2) , \, Y_{\tau^Y_{B(x,r)}} \in M \setminus B(x,r); \tau^Y_{B(x,r)}<s  \big)  \nn\\
			&\le 	\P_x\big( Y_{s/2} \in M\setminus B(x,r/2)  \big)+ 	\P_x\big(   Y_{\tau^Y_{B(x,r)}} \in M, \, Y_{s/2} \in  M\setminus B( Y_{\tau^Y_{B(x,r)}},r/2) ; \tau^Y_{B(x,r)}<s/2 \big)     \nn\\
			&\quad  + 	\P_x\big(  Y_{s/2} \in M, \, Y_{\tau^Y_{B(x,r)}} \in M \setminus B(Y_{s/2},r/2) ; s/2<\tau^Y_{B(x,r)}<s  \big)  \nn\\
		&\le 	\P_x\big( Y_{s/2} \in M\setminus B(x,r/2)  \big)  + 2	\sup_{z \in M, \, a \in (0,s/2]} \P_z\big(  Y_{a} \in  M\setminus B(z,r/2)  \big).
\end{align}
Using {\bf (A)} in the first inequality below,  \eqref{e:volumephi-integral} in the second and \eqref{e:scale-phi} in the third, we obtain for all $z\in M$ and $a\in (0,s/2]$,
\begin{align*}
&	\P_z\big( Y_a \in M\setminus B(z,r/2)  \big)  = \int_{M\setminus B(z,r/2)} q(a,z,y)\mu(dy) \\
&\le c_1a \int_{M\setminus B(z,r/2)}  \frac{\mu(dy)}{V(z,d(z,y)) \phi(d(z,y))} \le \frac{c_2 a}{\phi(r/2)} \le \frac{c_3s}{\phi(r)}.
\end{align*}
Combining this with \eqref{e:EP-1}, we obtain
\begin{align*}
\P_x\big( \tau^Y_{B(x,r)} < s \wedge \zeta^Y, \, s/2 < \zeta^Y  \big)  \le \frac{3c_3 s}{\phi(r)} \quad \text{for any $s\in (0,T_0)$}.
\end{align*}
Using this, we conclude that for any $t\in (0,T_0)$,
\begin{align*}
		\P_x\big( \tau^Y_{B(x,r)} < t \wedge \zeta^Y  \big) &= 	\P_x\big( \tau^Y_{B(x,r)} < t < \zeta^Y  \big)+  \sum_{n=1}^\infty 	\P_x\big( \tau^Y_{B(x,r)} <  \zeta^Y, \, \zeta^Y \in (2^{-n} t, 2^{1-n}t]\big)\\
		&\le 	\P_x\big( \tau^Y_{B(x,r)} < t \wedge \zeta^Y, \, t/2<\zeta^Y  \big)+  \sum_{n=1}^\infty 	\P_x\big( \tau^Y_{B(x,r)} <  2^{1-n}t \wedge \zeta^Y, \, 2^{-n}t< \zeta^Y \big)\\
		 &\le  \frac{3c_3t}{\phi(r)}  \bigg( 1+  \sum_{n=1}^\infty 2^{1-n} \bigg) = \frac{9c_3t}{\phi(r)}.
\end{align*}
Similarly, we also obtain $\wh\P_x\big( \wh\tau^Y_{B(x,r)} < t \wedge \wh\zeta^Y  \big)  \le c_4t/\phi(r)$ for any $t\in (0,T_0)$. \qed

	For an open subset $U \subset D$, we denote by $X^{\kappa,U}$  ($\wh X^{\kappa,U}$, respectively)
	the process $X^\kappa$  ($\wh X^\kappa$, respectively) killed upon exiting  $U$.

	 The following two results are established in \cite[Proposition 2.14, the paragraph above Proposition 2.14 and Theorem 2.15]{CKSV20}.
	\begin{prop}\label{p:heatkernel-existence}
		Suppose that $\kappa \in \mathbf{K}_1(D)$. The Hunt process $X^{\kappa}$ on $D$ corresponding to the  transition semigroup
		$(P_t^{\kappa})$ has a transition density $p^\kappa(t,x,y)$
		with respect to $\mu$
		such that
		\begin{align*}
			p^\kappa(t,x,y) \le q(t,x,y) \quad \text{for all $(t,x,y) \in (0,\infty) \times D \times D$}.
		\end{align*}
		Furthermore, if $D$ is relatively compact, then
		$p^\kappa(t,x,y)$ is continuous in $(t,y)$ for each fixed $x$,
		and continuous in $(t,x)$ for each fixed $y$.
	\end{prop}

	\begin{prop}\label{p:heatkernel-interior}
		Suppose that $\kappa \in \mathbf{K}_1(D)$.
		For all $T \in (0, \infty)$ and $a\in (0, 1]$,
		there exists 		$C=C(a,\phi,A_0,M,\sup_{t\le T} N^{D, \kappa}_{2^{-1}a}(t))>0$
		such that for any open $U \subset D$, the process $X^{\kappa,U}$ has a transition density $p^{\kappa, U}(t,x,y)$ with respect to $\mu$ and 
		\begin{align}\label{e:lower_bound}
			p^{\kappa,U}(t,x,y)  \ge C \wt q(t,x,y)
		\end{align}
		for all $t\in (0,T)$ and $x,y\in U$ with $\delta_U(x) \wedge \delta_U(y)\ge a\phi^{-1}(t)$. Moreover, if
		$\kappa \in \mathbf{K}_{\infty}(D)$ and
		$T_0=\infty$,
		then \eqref{e:lower_bound} holds for all $t>0$ and $x,y\in U$ with $\delta_U(x) \wedge \delta_U(y)\ge a\phi^{-1}(t)$.
	\end{prop}

	\section{Survival estimates}\label{s:surv}

	In this and the next section, we
	fix an open set $D\subset M$ and  $\kappa \in \mathbf{K}_1(D)$.  		We will write
	$X, \wh X, P_t, \wh P_t,p(t,x,y)$ and $p^U(t,x,y)$  instead of  $ X^\kappa,   \wh X^\kappa, P^\kappa_t, \wh P^\kappa_t, p^\kappa(t,x,y)$ and $ p^{\kappa,U}(t,x,y)$. 
 For an open subset $U\subset D$, let $\tau_{U}:=\inf\{s> 0: \, X_s\notin U\}$ and $\wh \tau_{U}:=\inf\{s> 0: \, \wh X_s\notin U\}$. 
	 
	 We begin with the following lemma which is a consequence of  the L\'evy system formulas \eqref{e:levy_systemX} and \eqref{e:levy_systemX-dual}. 
		 Since the proofs for the dual processes are the same, 
	 throughout the paper, we present the proofs for $X$ only.

	 \begin{lemma}\label{l:exit-jump}
	 	There exists 
		$C=C(\phi,A_0,M)>0$ 
		such that for any open set $U\subset D$ and any Borel $W \subset D$ with $r:=\text{\rm dist}(U,W)>0$, we have
	 	\begin{align*}
	 		\P_x ( X_{\tau_U} \in W)   \le \frac{C \E_x[\tau_U]}{\phi(r)} \quad \text{and} \quad  
	 		\wh\P_x ( \wh X_{\wh \tau_U} \in W)\le \frac{C \wh\E_x[\wh\tau_U]}{\phi(r)}.
	 	\end{align*}
	 \end{lemma}
	 	 \pf Using \eqref{e:levy_systemX} and \eqref{e:jump-integral}, we obtain
	 	 \begin{align*}
	 	 		\P_x ( X_{\tau_U} \in W)  = \E_x \int_0^{\tau_U} \int_W J(X_s, y) \mu(dy) ds \le \E_x \tau_U \sup_{z\in M} \int_{B(z,r)^c} J(z,y)\mu(dy) \le \frac{c_1\E_x[\tau_U]}{\phi(r)}.
	 	 \end{align*}
	 	\qed

	The following lemma was established in  \cite[Corollary 2.16]{CKSV20}. For completeness, we give the proof below.

	\begin{lemma}\label{l:survival-interior}
		(i)		For any $T>0$ and $a\in (0,1]$, there exists 
		$C=C(a,\phi,A_0,M,\sup_{t\le T} N^{D, \kappa}_{2^{-1}a}(t))>0$  
		such that for all $x \in D$ and $t \in (0, T)$ with $B(x,a\phi^{-1}(t)) \subset D$,
		\begin{equation}\label{e:survival-interior}
			\inf_{z \in B(x,a\phi^{-1}(t)/2)}\P_z(\tau_{B(x,a\phi^{-1}(t))}> t) \wedge \inf_{z \in B(x,a\phi^{-1}(t)/2)}\wh \P_z(\wh \tau_{B(x,a\phi^{-1}(t))}> t) 
			\ge C.
		\end{equation}
		(ii)		If $\kappa \in  \mathbf{K}_{\infty}(D)$ and
		$T_0=\infty$,
		then for any $a\in (0,1]$,  \eqref{e:survival-interior} holds for all $x\in D$ and  $t>0$ with $B(x,a\phi^{-1}(t)) \subset D$.
	\end{lemma}
	\pf For any $z \in B(x,a\phi^{-1}(t)/2)$,  using Proposition \ref{p:heatkernel-interior}, \eqref{e:VD}  and \eqref{e:scale-phi}, we get 
	\begin{align*}
		\P_z(\tau_{B(x,a\phi^{-1}(t))}> t) &\ge \P_z(\tau_{B(z,a\phi^{-1}(t)/6)}> t)   =\int_{B(z,a\phi^{-1}(t)/6)} p^{B(z,a\phi^{-1}(t)/6)}(t,z,y)\mu(dy) \\
		&  \ge c_1 \int_{B(z,a\phi^{-1}(t)/12)}
	\bigg( 	\frac1{V(z, \phi^{-1}(t))}
		\wedge \frac{t}{V(z,a\phi^{-1}(t)/3)\phi(\phi^{-1}(t)/3)}  \bigg)\mu(dy) \\ &\ge  \frac{c_2V(z,a\phi^{-1}(t)/12)}{V(z, \phi^{-1}(t))}  \ge c_3.
	\end{align*}
	\qed 
	
Let $T_0\in (0,\infty]$ be the constant from {\bf (A)}. Define
	\begin{align*}
		R_0:=(\diam(M)/k_0) \wedge \phi^{-1}(T_0),
	\end{align*}
	 where $k_0>1$ is the constant in \eqref{e:uniformlyperfect}. Define 
	$$B_D(x,r)=B(x,r) \cap D, \quad x\in M, \; r>0.$$
	
	\begin{lemma}\label{l:mean-exittime}
		(i) There exists 		$C=C(\phi,A_0,M)>0$ 
		such that   for all  $x\in D$ and $r \in (0,R_0/2)$,
		\begin{equation}\label{e:mean-exittime}
			\sup_{z \in B_D(x,r)}	\E_z[\tau_{B_D(x,r)}] \vee \sup_{z \in B_D(x,r)} \wh \E_z[\wh \tau_{B_D(x,r)}] \le C \phi(r).
		\end{equation}
		(ii)		If $\kappa \in  \mathbf{K}_{\infty}(D)$ and
		$T_0 =\infty$,
		then  \eqref{e:mean-exittime} hold for all $x\in D$ and $r>0$.
	\end{lemma}
	\pf Since the proofs are similar, we only prove  (i).
	
	  Let $x\in D$ and $r\in (0,R_0/2)$. By \eqref{e:uniformlyperfect}, we can choose $y \in B(x,2k_0 r) \setminus B(x,2r)$.  For all $z\in B_D(x,r)$ and $w\in B(y,r)$, by \eqref{e:VD} and \eqref{e:scale-phi}, we have
	  \begin{align*}
	  	\frac{1}{V(z, r)} \wedge \frac{\phi(r)}{V(z,d(z,w)) \phi(d(z,w))}\ge  \frac{\phi(r)}{V(y,(4k_0+3)r) \phi( (2k_0+2)r)} \ge \frac{c_1}{V(y,r)} .
	  \end{align*}
	Thus,  by  Proposition \ref{p:heatkernel-existence} and  {\bf (A)},  we have  for any $z\in B_D(x,r)$,
	\begin{align}\label{e:survival-upper}
		&\P_z(\tau_{B_D(x, r)} > \phi(r))  \le \int_{B_D(x, r)} p(\phi(r),z,w) \mu(dw)   \le \int_{B(x, r)}  q(\phi(r),z,w) \mu(dw) \nn\\
		& \le 1-\int_{B(y, r)}  q(\phi(r),z,w) \mu(dw) \le 1-c_2 \int_{B(y, r)} \wt q(\phi(r),z,w) \mu(dw) \le 1- c_1c_2.
	\end{align}
Set $\delta:=1-c_1c_2\in (0,1)$.		By the Markov property and \eqref{e:survival-upper},  we get for all $n\ge 1$,
	\begin{align*}
		& \sup_{z\in B_{D}(x,r)}	\P_z(\tau_{B_D(x,r)}> n\phi(r)) =  \sup_{z\in B_{D}(x,r)}\P_z\left( \P_{X_{\phi(r)}}\left(\tau_{B_D(x,r)}> (n-1)\phi(r) \right) ; X_{\phi(r)} \in B_D(x,r) \right) \nn\\
		&\le  \sup_{z'\in B_{D}(x,r)}	\P_{z'}(\tau_{B_D(x,r)}> (n-1)\phi(r)) \sup_{z\in B_{D}(x,r)}	\P_z(\tau_{B_D(x,r)}> \phi(r)) \nn\\
		&\le \cdots \le  \bigg( \sup_{z\in B_{D}(x,r)}	\P_z(\tau_{B_D(x,r)}> \phi(r))  \bigg)^n \le \delta^n.
	\end{align*}
	It follows that
	\begin{align*}
		\sup_{z\in B_D(x,r)}\E_z[\tau_{B_D(x,r)}] \le  \phi(r)\sum_{n=1}^\infty n\,\P_x\big(\tau_{B_D(x,r)} \in ((n-1)\phi(r), n\phi(r)]\big) \le \phi(r) \sum_{n=1}^\infty n \delta^{n-1} = c_3 \phi(r).
	\end{align*}	\qed
	
Denote by	 $\zeta$ and $\wh \zeta$  the lifetimes of $X$ and $\wh X$ respectively.  The following survival estimates play a crucial role in establishing the approximate factorization for Dirichlet heat kernels.
\begin{prop}\label{p:survival-main}
	There exists 
	$C=C(\phi,A_0,M)>0$ 
	such that for all  $z\in \overline D$,  $R\in (0,R_0)$ and  $x\in B_D(z,R/2)$,  
	\begin{align}\label{e:survival-main}
	\P_x ( \tau_{B_D(z, R)} <  \zeta  )   \le \frac{C \E_x [ \tau_{B_D(z,R)}]}{\phi(R)}  \quad \text{and} \quad 	\wh\P_x (  \wh\tau_{B_D(z, R)}  <\wh\zeta )  \le \frac{C \wh\E_x [ \wh\tau_{B_D(z,R)}]}{\phi(R)} .
	\end{align}
	Moreover, 	if $\kappa \in  \mathbf{K}_{\infty}(D)$ and
	$T_0 =\infty$, then  \eqref{e:survival-main} holds for all $z\in \overline D$, $R>0$ and $x\in B_D(z,R/2)$.
\end{prop}
	 
	To prove Proposition \ref{p:survival-main}, we adapt  the argument of  \cite[Proposition 3.1]{CC24}, which itself is based on the box method developed in \cite{BB90}. However, due to critical killing, condition  {\bf (EP)}$_{\phi,r_0,\le}$ in \cite{CC24} generally fails, and thus, some non-trivial modifications are required.

	 We begin with the following lemma, cf. \cite[Lemma 3.2]{CC24}.
	\begin{lemma}\label{l:box}
		 	(i)	There exists 			$C=C(\phi,A_0,M)>0$
			such that for any  $z\in \overline D$ and $R\in (0,R_0)$, the following holds  for all $x \in B_D(z, 3R/4)$ and $r \in (0,R/8]$:
		\begin{align}\label{e:box}
		\P_x ( \tau_{B_D(x, r)}  <\zeta )  \le C \sqrt{ \frac{\E_x [ \tau_{B_D(x,2r)}]}{ \phi(r)}} \quad \text{and} \quad 	\wh\P_x ( \wh\tau_{B_D(x, r)} < \wh\zeta  )  \le C \sqrt{ \frac{\wh\E_x [ \wh\tau_{B_D(x,2r)}]}{ \phi(r)}}.
		\end{align}
		
		\noindent (ii) If  $\kappa \in  \mathbf{K}_{\infty}(D)$ and
		$T_0 =\infty$, then  \eqref{e:box} hold for all $z\in \overline D$, $R>0$, $x \in B_D(z, 3R/4)$ and $r \in (0,R/8]$.
	\end{lemma}
	\pf  We  only present the proof for the statement concerning $X$ in (i), as the other statements follow from similar arguments.	We decompose $	\P_x ( \tau_{B_D(x, r)}  < \zeta ) $ as follows:
	\begin{align*}
		&	\P_x ( \tau_{B_D(x, r)}  < \zeta )  =  \P_x \big( X_{\tau_{B_D(x, r)}} \in D \big)\\
		&=\P_x \big( X_{\tau_{B_D(x, r)}} \in B_D(x,2r) \big) +  \P_x \big( X_{\tau_{B_D(x, r)}} \in D\setminus B_D(x,2r) \big)=:I_1+I_2.
	\end{align*}

	By Lemma \ref{l:exit-jump}, $I_2\le c_1\E_x[\tau_{B_D(x,r)}]/\phi(r)$. Thus, since $\E_x[\tau_{B_D(x,r)}] \le c_2\phi(r)$ by Lemma \ref{l:mean-exittime}(i),  we have
	\begin{align*}
		I_2 \le  \frac{c_1\E_x [ \tau_{B_D(x,r)}]}{\phi(r)}  \le  c_3\sqrt{ \frac{\E_x [ \tau_{B_D(x,r)}]}{\phi(r)} } \le  c_3\sqrt{ \frac{\E_x [ \tau_{B_D(x,2r)}]}{\phi(r)} }.
	\end{align*}
	
	For $I_1$, using the Markov inequality, we see that for all $t>0$,
	\begin{align}\label{e:box-I1}
			I_1 = \P_x (  \tau_{B_D(x, r)} <\tau_{B_D(x, 2r)} ) & \le   \P_x ( \tau_{B_D(x, 2r)}>t ) +  \P_x (   \tau_{B_D(x, r)} <t\wedge \tau_{B_D(x, 2r)} )\nn\\
				&\le t^{-1} \E_x  [ \tau_{B_D(x,2r)}] +  \P_x (  \tau_{B_D(x, r)} < t\wedge \zeta  ).
	\end{align}
	Note that  $X_t$ and $Y_t$ has the same distribution for $t<\zeta$, and that $\zeta \le \tau^Y_D$. Thus, we have $	\{ \tau_{B_D(x, r)} < t\wedge \zeta \}\subset \{  \tau^Y_{B_D(x, r)}  <t \wedge \tau^Y_D \}$. Moreover, on the event $\{\tau^Y_{B_D(x, r)}  <t \wedge \tau^Y_D \}$, we have $Y_{  \tau^Y_{B_D(x, r)}} \in D\setminus B(x,r)$ and $ \tau^Y_{B_D(x, r)}  <t$. It follows that $\{\tau^Y_{B_D(x, r)}  <t \wedge \tau^Y_D\} \subset \{ \tau^Y_{B(x,r)}<t \wedge \zeta^Y\}$. Consequently, we have
	 $$\P_x (\tau_{B_D(x, r)} < t\wedge \zeta  ) \le \P_x (  \tau^Y_{B_D(x, r)}  <t \wedge \tau^Y_D) \le \P_x( \tau^Y_{B(x,r)}<t \wedge \zeta^Y).$$ Thus, using Lemma \ref{l:EP}, we get from \eqref{e:box-I1} that for all $t>0$,
\begin{align*}
	I_1 \le \frac{\E_x  [ \tau_{B_D(x,2r)}]}{t} + \P_x( \tau^Y_{B(x,r)}<t \wedge \zeta^Y)   \le   \frac{\E_x  [ \tau_{B_D(x,2r)}]}{t} + \frac{c_4t}{\phi(r)}.
\end{align*}
	Choosing $t= \sqrt{\phi(r)\E_x  [ \tau_{B_D(x,2r)}]}$, we get
	$I_1 \le c_5 \sqrt{ \E_x [ \tau_{B_D(x,2r)}]/ \phi(r)} $. This completes the proof.	\qed

	By following the proof of \cite[Corollary 3.3]{CC24}, we obtain the following corollary 	of Lemma \ref{l:box}. We give the proof for completeness.
	
	\begin{cor}\label{c:box}
	(i)	There exists 
	$C=C(\phi,A_0,M)>0$ such that
		for all $z\in \overline D$, $R\in (0,R_0)$ and $x\in B_D(z,3R/4)$, we have
		\begin{align}
			\frac{\E_x[\tau_{B_D(z,R)}]}{\phi(R) \P_x(\tau_{B_D(z,R)}<\zeta)} \ge C \bigg(\P_x(\tau_{B_D(z,R)}<\zeta) + 	\frac{\E_x[\tau_{B_D(z,R)}]}{\phi(R)} \wedge 1 \bigg),\label{e:box-corollary-1}\\
			\frac{\wh\E_x[\wh\tau_{B_D(z,R)}]}{\phi(R) \wh\P_x(\wh\tau_{B_D(z,R)}<\wh \zeta)} \ge C \bigg(\wh\P_x(\wh\tau_{B_D(z,R)}<\wh \zeta) + 	\frac{\wh\E_x[\wh\tau_{B_D(z,R)}]}{\phi(R)} \wedge 1 \bigg).\label{e:box-corollary-2}
		\end{align}
		
		\noindent (ii) If  $\kappa \in  \mathbf{K}_{\infty}(D)$ and
		$T_0=\infty$, then  \eqref{e:box-corollary-1}--\eqref{e:box-corollary-2} hold for all $z\in \overline D$, $R>0$ and $x \in B_D(z, 3R/4)$.
	\end{cor}
	\pf Since the proofs are similar,  we  only present the proof for  \eqref{e:box-corollary-1} in (i).
	
Let  $z\in \overline D$, $R\in (0,R _0)$ and $x\in B_D(z,3R/4)$. Then $B_D(x,R/4) \subset B_D(z,R)$ so that $\tau_{B_D(x, R/8)}\le \tau_{B_D(x, R/4)} \le \tau_{B_D(z, R)}$. 	Hence,	by applying Lemma \ref{l:box}(i) with $r=R/8$ and using \eqref{e:scale-phi}, we get
		\begin{align}\label{e:box-corollary-key}
				\P_x ( \tau_{B_D(z, R)}  <\zeta ) \le 	\P_x ( \tau_{B_D(x, R/8)}  <\zeta )  \le c_1 \sqrt{ \frac{\E_x [ \tau_{B_D(x,R/4)}]}{ \phi(R/8)}} \le c_2 \sqrt{ \frac{\E_x [ \tau_{B_D(z,R)}]}{ \phi(R)}}.
		\end{align} 
	By applying \eqref{e:box-corollary-key} twice and using the fact that $\sqrt a \ge a\wedge 1$ for all $a\ge 0$, we obtain
	\begin{align*}
	\frac{\E_x[\tau_{B_D(z,R)}]}{\phi(R) \P_x(\tau_{B_D(z,R)}<\zeta)} & \ge  c_1\sqrt{	\frac{\E_x[\tau_{B_D(z,R)}]}{\phi(R)} } = \frac{c_1}{2}\sqrt{	\frac{\E_x[\tau_{B_D(z,R)}]}{\phi(R)} }  +  \frac{c_1}{2}\sqrt{	\frac{\E_x[\tau_{B_D(z,R)}]}{\phi(R)} } \\
	&\ge  c_2 \bigg(\P_x(\tau_{B_D(z,R)}<\zeta) + 	\frac{\E_x[\tau_{B_D(z,R)}]}{\phi(R)} \wedge 1 \bigg).
	\end{align*} \qed

	Now by repeating the proof of \cite[Proposition 3.1]{CC24}, we obtain Proposition \ref{p:survival-main}. 
	
	\medskip
	
	\noindent \textbf{Proof of Proposition \ref{p:survival-main}.}  
	By \eqref{e:jump-integral},
	  condition {\bf (Jt)}$_{\phi, R_0, \le}$ in \cite{CC24} holds. In the proof of \cite[Proposition 3.1]{CC24}, condition  {\bf (EP)}$_{\phi, R_0, \le}$ is used only through  
	  \cite[Lemma 3.2 and Corollary 3.3]{CC24}. These two results
	  can be replaced by Lemma \ref{l:box} and Corollary \ref{c:box} respectively. Indeed, although
	  the range of $r$ in \cite[Lemma 3.2(a)]{CC24} is $(0,R/4]$, while  Lemma \ref{l:box} covers $r\in (0,R/8]$, there is no conflict: \cite[Lemma 3.2]{CC24} is applied only  for $r\le (3/4\pi^2) R$ to get \cite[(3.8)]{CC24}. Hence, since $\zeta= \tau_D$ and $\wh \zeta=\wh \tau_D$, by following the proof of \cite[Proposition 3.1]{CC24}, with \cite[Lemmas 2.8(a), 3.2(a) and Corollary 3.2]{CC24} replaced by   Lemmas \ref{l:exit-jump} and \ref{l:box}, and Corollary \ref{c:box}, respectively, we obtain the  result. \qed

	\section{Approximate factorizations in $\eta$-fat open set}\label{s:fact}
	
	\begin{definition}
		\rm 
		Let $\eta \in (0,1/2]$. We say that an open set $D\subset M$ is $\eta$-fat if there exists
		$\overline R_1\in (0,\diam(D)]$ such that for all
		$x \in \overline{D}$
		and  $r\in(0, \overline R_1)$,
		there is $z_{x,r} \in D$ such that $B(z_{x,r},\eta r) \subset  B_D(x,r)$. The pair $(\overline R_1,\eta)$ is called the characteristics of the $\eta$-fat open set $D$.
	\end{definition}
	
	In this section, we always assume 
	that  $D\subset M$ is $\eta$-fat with characteristics $( \overline R_1,\eta)$.  Set $R_1:=(\overline R_1/k_0) \wedge \phi^{-1}(T_0)$, where $k_0>1$ is the constant in \eqref{e:uniformlyperfect}.

	\subsection{Approximate factorization of the Dirichlet heat kernel}\label{ss:factorization-heatkernel}
	
	Throughout this subsection, we fix a positive constant $T\in (0,T_0)$. Define for $t\in (0,T)$,
	$$r_t=\frac{R_1\phi^{-1}(t)}{3\phi^{-1}(T)} \quad \text{if $R_1<\infty$} \quad \text{and} \quad r_t=\phi^{-1}(t) \quad \text{if $R_1=\infty$}.$$
	In either case, by the monotonicity of $\phi$ and  \eqref{e:scale-phi}, we have $r_t <R_1/3$ and $\phi(r_t) \asymp t$ for all $t\in (0,T)$. 
		For  any $(t,x)\in(0,T)\times D$,	
		we define an open neighborhood $U(x,t)$ of $x$ and an open ball $W(x,t)\subset D\setminus U(x,t)$ as follows (cf. \cite{CKSV20}):

	\smallskip
	
	Since $D$ is $\eta$-fat, there is $z_{x,t} \in D$ such that $B(z_{x,t} ,3\eta r_t) \subset B_D(x,3r_t)$.

	(1) If $\delta_D(x) \le  2\eta r_t$, we define
	\begin{align*} 
	U(x,t)=B_D(x,\frac{\eta r_t}{4k_0}) \quad \text{ and } \quad W(x,t)=B(z_{x,t} ,\frac{\eta r_t}{4k_0}).
	\end{align*}  
	
	(2)	 If  $\delta_D(x) > 2 \eta r_t$, we choose
	$y_{x,t}\in M$
	such that $\eta r_t/k_0 \le d(x,y_{x,t})\le \eta r_t$, where
	$k_0>1$
	is the constant in \eqref{e:uniformlyperfect}.
	Then we define 
	\begin{align*}
		U(x,t)=B(x,\frac{\eta r_t}{4k_0}) \quad \text{ and } \quad  W(x,t)=B(y_{x,t},\frac{\eta r_t}{4k_0}).
	\end{align*}

	Note that in either case,  we have
	\begin{align}
	\frac{\eta r_t}{2k_0} \le  \text{dist}(U(x,t), W(x,t))\le
		4r_t,\label{e:UW}\\
		 \text{dist}(W(x,t),B_D(x,3r_t)^c) \ge \frac{\eta r_t}{2}. \label{e:W-interior}
	\end{align}
Moreover,	by \eqref{e:VD}, there exists $C=C(\eta, \overline R_1, \phi,M)\ge 1$ such that in either case,
	\begin{align}\label{e:W-volume}
	C^{-1}  V(x, r_t) \le 	\mu(W(x,t)) \le C V(x, r_t).
	\end{align} 
	
	\begin{lemma}\label{l:boundary-term-equivalent}
(i)	For any $a> 1$, there are comparison constants depending only on 
$a,T,\overline R_1,\eta, \phi,A_0$ and $M$  such that 
	for any $(t,x) \in (0,T) \times D$,
		\begin{align}\label{e:boundary-term-equivalent-1}
			&\P_x(\zeta>t)  \asymp \P_x(\zeta>at) \asymp 	\P_x(\tau_{B_D(x,3r_t)}>at)  \asymp  t^{-1} \E_x[\tau_{U(x,t)}]  \nn\\[3pt]
			& \asymp  \P_x(X_{\tau_{U(x,t)}} \in D)\asymp
			\P_x(X_{\tau_{U(x,t)}} \in W(x,t))  
		\end{align}
		and
		\begin{align}\label{e:boundary-term-equivalent-2}
			&\wh\P_x(\wh\zeta>t)  \asymp 	\wh\P_x(\wh\zeta>at)  \asymp 	\wh\P_x(\wh\tau_{B_D(x,3r_t)}>at)   \asymp  t^{-1}\wh \E_x[\wh\tau_{U(x,t)}]\nn\\[3pt]
			&\asymp \wh\P_x(\wh X_{\wh \tau_{U(x,t)}} \in D)\asymp
		\wh \P_x(\wh X_{\wh \tau_{U(x,t)}} \in W(x,t)),\end{align}	where $U(x,t)$  and $W(x,t)$ are the open sets defined in the beginning of this subsection.
		
		\noindent (ii) If  $\kappa \in  \mathbf{K}_{\infty}(D)$ and
		$T_0 =\overline R_1=\infty$, then for any $a>1$, there are comparison constants 
		depending on  $a,\eta, \phi,A_0$ and $M$  
		such that \eqref{e:boundary-term-equivalent-1}--\eqref{e:boundary-term-equivalent-2} hold for any $(t,x) \in (0,\infty) \times D$.
	\end{lemma}
	\smallskip
	\pf  We only give the proof for \eqref{e:boundary-term-equivalent-1} in (i), as the other statements can be proved similarly.
	
	Let $t \in (0,T)$ and $x \in D$. 
	By the L\'evy system formula \eqref{e:levy_systemX}, \eqref{e:J-density}, \eqref{e:UW}, \eqref{e:W-volume}, \eqref{e:VD} and \eqref{e:scale-phi},  
	\begin{align*}
		\P_x(X_{\tau_{U(x,t)}} \in D)\ge 	\P_x(X_{\tau_{U(x,t)}} \in W(x,t))&= \E_x\left[\int_0^{\tau_{U(x,t)}} \int_{W(x,t)}J(X_s,w)\mu(dw)ds\right]\nn\\
		&\ge 
		\frac{c_1\mu(W(x,t))\E_x[\tau_{U(x,t)}]}{V(x, 4r_t)\phi(4r_t)} \ge  \frac{c_2\E_x[\tau_{U(x,t)}]}{t}.
	\end{align*}
Moreover, by Proposition \ref{p:survival-main} and \eqref{e:scale-phi}, we have
\begin{align*}
	\P_x(X_{\tau_{U(x,t)}} \in D)& = \P_x(\tau_{U(x,t)}<\zeta) \le \frac{c_3 \E_x [\tau_{U(x,t)} ]}{\phi(\eta r_t/(4k_0))
	} \le \frac{c_4\E_x [\tau_{U(x,t)} ]}{t}.
\end{align*}
Thus, 
\begin{align}\label{e:boundary-term-1}
	\P_x(X_{\tau_{U(x,t)}} \in D) \asymp 	\P_x(X_{\tau_{U(x,t)}} \in W(x,t)) \asymp t^{-1}\E_x[\tau_{U(x,t)}].
\end{align}
Clearly, $		\P_x(\tau_{B_D(x,3r_t)}>at)  \le 	\P_x(\zeta>at)  \le \P_x(\zeta>t) $. Further, by the Markov inequality, 
\begin{equation*}
	\P_x(\zeta>t)  \le \P_x(\tau_{U(x,t)}>t)+\P_x(X_{\tau_{U(x,t)}} \in D) \le   t^{-1}\E_x[\tau_{U(x,t)}]+\P_x(X_{\tau_{U(x,t)}} \in D).
\end{equation*}
 Combining this with \eqref{e:boundary-term-1}, we see that, 
 to obtain \eqref{e:boundary-term-equivalent-1}, it suffices to show that 
\begin{align}\label{e:boundary-term-claim}
		\P_x(\tau_{B_D(x,3r_t)}>at)  \ge c_5 \P_x(X_{\tau_{U(x,t)}} \in W(x,t)) .
\end{align}

By \eqref{e:scale-phi}, we have $\phi^{-1}(at) \le c_6 r_t$. Thus, by \eqref{e:W-interior} and Lemma \ref{l:survival-interior}(i),  $$\inf_{w\in W(x,t)} 	\P_w(\tau_{B_D(x,3r_t)}>at)   \ge \inf_{w\in W(x,t)} \P_w(\tau_{B(w,\eta r_t/2)} >at)  \ge c_7.$$
Hence, using the strong Markov property, we obtain 
\begin{align*}
	&	\P_x(\tau_{B_D(x,3r_t)}>at) \ge \P_x(\tau_{B_D(x,3r_t)}>at + \tau_{U(x,t)}) \\
	&\ge 
	\E_x\left[\inf_{w\in W(x,t)} \P_w(\tau_{B_D(x,3r_t)}>at) : X_{\tau_{U(x,t)}} \in W(x,t)\right]\ge c_7\P_x(X_{\tau_{U(x,t)}} \in W(x,t)),
\end{align*}
proving that \eqref{e:boundary-term-claim} holds. The proof is complete. \qed

For $r>0$, define
\begin{align*}
	D(r):= \left\{ x \in D : \delta_D(x) < r\right\}.
\end{align*}

\begin{cor}\label{c:boundary-term-equivalent}
	
	(i) There are comparison constants depending only on 	$T,\overline R_1, \eta,\phi,A_0$ and $M$  
	such that 	for any $(t,x) \in (0,T) \times D$,
	\begin{align}
		\P_x(\tau_{U(x,t)}>t) &   \asymp 	\P_x(\zeta>t) \asymp  \P_x( X_{\tau_{D(\eta r_t/2)}} \in D),\label{e:boundary-term-equivalent-3}\\
		\wh	\P_x(\wh\tau_{U(x,t)}>t)  &  \asymp 	\wh\P_x(\wh\zeta>t)\asymp  \wh\P_x( \wh X_{\wh\tau_{D(\eta r_t/2)}} \in D).\label{e:boundary-term-equivalent-4}
	\end{align}
	
	\noindent (ii) If  $\kappa \in  \mathbf{K}_{\infty}(D)$ and
	$T_0 =\overline R_1=\infty$, then there are comparison constants depending only on   
	$\eta, \phi,A_0$ and $M$ 
	such that \eqref{e:boundary-term-equivalent-3}--\eqref{e:boundary-term-equivalent-4} hold for any $(t,x) \in (0,\infty) \times D$.
\end{cor}
\pf We prove \eqref{e:boundary-term-equivalent-3} in (i) only, as the other statements can be proved similarly.

Let $t\in (0,T)$ and  $x\in D$. 
By \eqref{e:scale-phi}, there exists $\eps \in (0,1)$ independent of $t$ and $x$ such that  $B_D(x, 3r_{\eps t})\subset B_D(x,  \eta r_t/(4k_0)) \subset U(x,t)$. Hence, applying Lemma \ref{l:boundary-term-equivalent} with $a=1/\eps$, we obtain
\begin{align*}
	\P_x( \zeta >t) \ge 	\P_x(\tau_{U(x,t)}>t)  \ge \P_x(\tau_{B_D(x,3r_{\eps t})}>t)  \ge c_1 \P_x( \zeta >t).
\end{align*}
This proves the first comparison in \eqref{e:boundary-term-equivalent-3}.

By \eqref{e:W-interior}, we have $W(x,t) \subset D\setminus D(\eta r_t/2)$. It follows that
\begin{align*}
	\P_x( X_{\tau_{D(\eta r_t/2)}} \in D) \ge \P_x ( \text{$X_s \in W(x,t)$ for some $s<\zeta$}) \ge \P_x ( X_{\tau_{U(x,t)}} \in  W(x,t)).
\end{align*}
Hence, for the second comparison, by Lemma \ref{l:boundary-term-equivalent}, it suffices to show that 
\begin{align}\label{e:boundary-term-equivalent-claim}
	\P_x( X_{\tau_{D(\eta r_t/2)}} \in D)\le c_2	\P_x(\zeta>t)
\end{align}
for some $c_2>0$ independent of $t$ and $x$.  By \eqref{e:scale-phi}, there exists $\eps' \in (0,1)$ independent of $t$  such that  $6r_{\eps't} \le r_t$.  If $\delta_D(x) > 2\eta r_{\eps' t}$, then using Lemma \ref{l:survival-interior}, we obtain
\begin{align*}
		\P_x(\zeta>t) \ge 	\P_x(\tau_{B(x,2\eta r_{\eps' t})}>t) \ge c_3 \ge c_3\P_x( X_{\tau_{D(\eta r_t/2)}} \in D),
\end{align*}
proving that \eqref{e:boundary-term-equivalent-claim} holds. Assume  $\delta_D(x) \le  2\eta r_{\eps' t}$. Then we get $U(x,\eps' t) = B_D(x, \eta r_{\eps't}/(4k_0)) \subset D(3\eta r_{\eps't}) \subset D(\eta r_t/2)$. Hence, applying Lemma \ref{l:boundary-term-equivalent} with $a=1/\eps'$, we arrive at
\begin{align*}
		\P_x( X_{\tau_{D(\eta r_t/2)}} \in D)  \le 	\P_x( X_{\tau_{U(x,\eps't)}} \in D)  \le c_4 \P_x(\zeta>t).
\end{align*} 
This complete the proof. \qed

By the first comparisons in \eqref{e:boundary-term-equivalent-1} and \eqref{e:boundary-term-equivalent-2}, one sees that  for any $x\in D$, the function $t\mapsto \P_x(\zeta>t)$ decreases  at most polynomially. In the following proposition,  we refine this by showing that the decay rate is strictly less than $1$. This property is crucial for establishing an approximate factorization for the Green function later on.

\begin{prop}\label{p:boundary-term-scale}
	(i)	There exist constants $\theta \in [0,1)$ and $C=C(T)\ge 1$ such that 
	\begin{align}\label{e:boundary-term-scale}
		1\le 	\frac{\P_x(\zeta>s)}{\P_x(\zeta>t)} \le C \bigg( \frac{t}{s}\bigg)^{\theta} \quad \text{and} \quad 	1\le 	\frac{\wh\P_x(\wh\zeta>s)}{\wh\P_x(\wh\zeta>t)} \le C \bigg( \frac{t}{s}\bigg)^{\theta}
	\end{align}
	for all $x\in D$ and  $0<s\le t<T$.
	
	\noindent (ii) If  $\kappa \in  \mathbf{K}_{\infty}(D)$ and
	$T_0 =\overline R_1=\infty$, then there  exist constants $\theta \in [0,1)$ and $C\ge 1$ such that \eqref{e:boundary-term-scale} hold for all $x\in D$ and $0<s\le t$.
\end{prop}
\pf We only prove the first set of inequalities in \eqref{e:boundary-term-scale} in  (i).

Let $x\in D$. It is clear that $\P_x(\zeta>s) \ge \P_x(\zeta>t)$ for all $0<s\le t<T$. For the other inequality, by Lemma \ref{l:boundary-term-equivalent}(i), it suffices to show that
\begin{align}\label{e:boundary-term-scale-claim}
	\frac{\E_x[\tau_{U(x,s)}]}{\E_x[\tau_{U(x,t)}]} \le c_1 \bigg( \frac{s}{t} \bigg)^a \quad \text{for all $0<s\le t<T$},
\end{align}
for some $a \in (0,1)$ and $c_1>0$ independent of $x$. 

Let $t\in (0,T)$. By  \eqref{e:scale-phi}, there exists  $\eps \in (0,1)$ independent of $t$ such that $5r_{\eps t} \le  \eta r_t/(4k_0)$. Define open sets $U(x, \cdot)$ and $W(x,\cdot)$ as in the beginning of this subsection.
By \eqref{e:UW} and \eqref{e:W-interior}, we see that for any $y \in W(x, \eps t)$,
\begin{align*}
	\delta_D(y) \ge  \eta r_{\eps t}/2 \quad \text{and} \quad B(y, \eta r_{\eps t}/2) \subset  B_D(x, 5r_{\eps t}) \subset  B_D(x, \eta r_{t}/(4k_0)) \subset U(x,t).
\end{align*}
Thus, using the strong Markov property, we obtain
\begin{align}\label{e:boundary-term-scale-1}
	\E_x[\tau_{U(x,t)}] &\ge  \E_x\Big[\tau_{U(x,t)}: X_{\tau_{U(x,\eps t)}} \in W(x, \eps t)\Big] \nn\\
	&\ge \E_x[\tau_{U(x,\eps t)}] +  \inf_{y \in W(x,\eps t)} \E_y [ \tau_{B(y, \eta r_{\eps t}/2)} ]\, \P_x \Big(  X_{\tau_{U(x,\eps t)}} \in W(x, \eps t)\Big) .
\end{align}
Since $\delta_D(y) \ge \eta r_{\eps t}/2$ for any $y \in W(x,\eps t)$, by Lemma \ref{l:survival-interior}(i), we get
\begin{align*}
	\inf_{y \in W(x,\eps t)} \E_y [ \tau_{B(y, \eta r_{\eps t}/2)} ] \ge \eps t  \inf_{y \in W(x,\eps t)} \P_y ( \tau_{B(y, \eta r_{\eps t}/2)} > \eps t)  \ge c_2 \eps t.
\end{align*}
Note that, by Lemma \ref{l:boundary-term-equivalent}(i), we have that, $\eps t\,\P_x (  X_{\tau_{U(x,\eps t)}} \in W(x, \eps t))  \ge  c_3\E_x[\tau_{U(x,\eps t)}]$. Hence, we deduce from \eqref{e:boundary-term-scale-1} that
\begin{align}\label{e:boundary-term-scale-2}
	\E_x[\tau_{U(x,t)}] \ge (1+ c_2c_3) \E_x[\tau_{U(x,\eps t)}] \quad \text{for all $t\in (0,T)$}.
\end{align}

We now choose  $0<s\le t <T$ and let $n\ge 1$ be such that $\eps^n t < s \le \eps^{n-1}t$.  Applying \eqref{e:boundary-term-scale-2} $(n-1)$-times, we obtain
\begin{align*}
	\E_x[\tau_{U(x,t)}] &\ge (1+c_2c_3) 	\E_x[\tau_{U(x,\eps t)}]  \ge \cdots \ge  (1+c_2c_3)^{n-1} 	\E_x[\tau_{U(x,\eps^{n-1} t)}]\\
	&\ge (1+c_2c_3)^{n-1} 	\E_x[\tau_{U(x,s)}] \ge (1+c_2c_3)^{-1} (t/s)^{ \log (1+c_2c_3) / |\log \eps|}\,\E_x[\tau_{U(x,s)}] ,
\end{align*}
proving that \eqref{e:boundary-term-scale-claim} holds. The proof is complete. \qed

	Using Lemma \ref{l:boundary-term-equivalent} and Corollary \ref{c:boundary-term-equivalent},  and following the arguments in \cite[Theorems 2.22 and 2.23]{CKSV20}, we obtain  the next theorem. We omit the proof.
	
	\begin{thm}\label{t:factorization-heatkernel} 	Let $D$ be an $\eta$-fat open set with characteristics $(\overline R_1, \eta)$.

\noindent 	(i) 
		If $\kappa \in \mathbf{K}_1(D)$, then
		for every $T>0$, there exists 
		$C=C(T,\overline R_1,\eta, \phi,A_0,M) \ge 1$ 
		such that for any $(t,x,y) \in (0,T) \times D \times D$,
		\begin{equation}\label{e:factorization-heatkernel}
			C^{-1} \P_x(\zeta>t) \wh  \P_y( \wh \zeta>t) \wt q(t,x,y) \le p(t,x,y)\le C \P_x(\zeta>t) \wh  \P_y(\wh \zeta>t) \wt q(t,x,y).
		\end{equation}
			(ii)					If $T_0=\overline R_1=\infty$ and $\kappa \in \mathbf{K}_\infty(D)$, then
			 there exists 
			$C=C(\eta, \phi,A_0,M) \ge 1$ 
			such that \eqref{e:factorization-heatkernel} holds  for any $(t,x,y) \in (0,\infty) \times D \times D$.
	\end{thm}

	 We next establish large time heat kernel estimates in the case when $D$ is relatively compact.
	 For this, we first analyze the spectral properties of the semigroup $(P_t)_{t\ge 0}$.

	   Let $\sA$ and $\wh \sA$ be the infinitesimal generators of $(P_t)_{t\ge 0}$ and $(\wh P_t)_{t\ge 0}$ respectively. Denote by $\sigma(\sA)$ and  $\sigma(P_t)$  the spectra of $\sA$ and $P_t$ respectively. By Theorem \ref{t:factorization-heatkernel},   $\sup_{x,y \in D}p(t,x,y)<\infty$ for any $t>0$. Hence, for any $t>0$, $P_t$ is a compact operator on $L^2(D)$.  By \cite[Corollary 2.3.7 and the paragraph below it]{Pa83}, $\sigma(\sA)$ is countable and  consists only of eigenvalues. Write the eigenvalues of $\sA$ as  $(-\lambda_n)_{n\ge 1}$, repeated according to their multiplicity and ordered so that  $\text{\rm Re\,}\lambda_n$ is non-decreasing.  Note that $\sigma(\wh\sA)$ also consists only of eigenvalues and the eigenvalues of $\wh \sA$ are $(-\overline{\lambda_n})_{n\ge 1}$, where $\overline \lambda_n$ denotes the complex conjugate of $\lambda_n$.
	  Since $P_t$ is compact and is a $L^2(D)$-contraction, $\sigma(P_t)\setminus \{0\}$
	   consists of eigenvalues only and is a subset of $\{z \in \mathbb C: \text{Re\,}z \le 0\}$. Hence,	by \cite[Theorem 2.2.4]{Pa83}, we get $\sigma(P_t)\setminus \{0\}= \{e^{-\lambda_nt}: n \ge 1\}$ and $\text{\rm Re\,}\lambda_1 \ge 0$.

	   Let $\psi_1$ be an eigenfunction of $\sA$ associated with $-\lambda_1$ such that $\lVert \psi_1\rVert_2=1$. Then, for any $t>0$, $\psi_1$ is an eigenfunction of $P_t$ associated with $e^{-\lambda_1t}$. By Proposition \ref{p:heatkernel-interior}, $P_t$ admits a strictly positive transition density. Thus, by Jentzsch's  theorem \cite[Theorem V.6.6]{Sc74},  $\psi_1$ can be chosen to be strictly positive almost everywhere, $\lambda_1\in \R$  and  $\lambda_1<\text{\rm Re\,}\lambda_2$. 
	  	     Applying Jentzsch's  theorem to $\wh P_t$, we can choose an eigenfunction $\wh \psi_1$ of $\wh \sA$ associated with $-\overline{\lambda_1} = -\lambda_1$ which is strictly positive almost everywhere and $\int_D \psi_1(x) \wh \psi_1(x)\mu(dx)=1$. By definitions, we have
	  	     \begin{align*}
	  	     	\psi_1(\cdot) = e^{\lambda_1 } \int_D p(1,\cdot, y) \psi_1(y)\mu(dy) \quad \text{and} \quad 	\wh\psi_1(\cdot) = e^{\lambda_1 } \int_D p(1,y, \cdot) \wh\psi_1(y)\mu(dy) \quad \text{in $L^2(D)$}.
	  	     \end{align*}
	Thus,  since $p(1,\cdot,y)$ and $p(1,y,\cdot)$ are bounded and continuous  by Proposition \ref{p:heatkernel-existence}, we can assume  that  $\psi_1$ and $\wh \psi_1$ are continuous on $D$.

	The following result follows from \cite[Lemma 2.2 and Remark 1.13]{RSZ17}. Although \cite{RSZ17} assumes that $\lambda_1>0$, the proof of  \cite[Lemma 2.2]{RSZ17} remains valid without this assumption.
	\begin{lemma}\label{l:largetime-1}
	For any 	$a\in (0, \text{\rm Re\,}\lambda_2)$,
	there exists  $C=C(a)>0$ such that
	\begin{align*}
		\left| p(t,x,y) - e^{-\lambda_1 t} \psi_1(x) \wh \psi_1(y) \right|  \le C e^{-at} \left( \int_D p(1,x,z)^2 \mu(dz) \right)^{1/2}\left( \int_D p(1,z,y)^2 \mu(dz) \right)^{1/2}
	\end{align*}
	 for any $(t,x,y) \in [1,\infty) \times D \times D$. 
	\end{lemma}

	We will use the following elementary property of 	volume doubling spaces.
	\begin{lemma}\label{l:infimum-volume}
	Let $U \subset M$ be a bounded set. For any  $r>0$, 
	there exists $C=C( r/\diam(U), M)>0$ such that $\inf_{x\in U}	V(x,r) \ge C\sup_{x \in U} V(x,r).$
	Consequently, $\inf_{x\in U} V(x,r)>0$.
	\end{lemma}
	\pf For all $x,y \in U$, using \eqref{e:VD}, we have
	\begin{align*}
		V(x,r) \ge c_1 \bigg( \frac{r}{r + \diam(U)} \bigg)^\alpha V(x,r + \diam(U)) \ge c_1 \bigg( \frac{r}{r + \diam(U)} \bigg)^\alpha V(y,r).
	\end{align*}
	This proves the lemma. \qed 
	\begin{thm}\label{t:largetime-bounded}
		Let $D$ be a relatively compact $\eta$-fat open set. 
		If $\kappa \in \mathbf{K}_1(D)$, then
		there exists $C\ge 1$ such that for any $(t,x,y) \in [3,\infty) \times D \times D$,
		\begin{align*}
		C^{-1} \P_x(\zeta>1) \P_y(\zeta>1)e^{-\lambda_1 t}\le 	p(t,x,y) \le C \P_x(\zeta>1) \P_y(\zeta>1)e^{-\lambda_1 t},
		\end{align*}
		where $-\lambda_1:= \sup \text{\rm Re\,}\sigma(\sA)$.
\end{thm}
	\pf Let $a:= (\lambda_1 +  \text{\rm Re\,} \lambda_2)/2$. By Proposition \ref{p:heatkernel-existence}, {\bf (A)} (with Remark \ref{r:extend}) and Lemma \ref{l:infimum-volume}, 
	\begin{align}\label{e:largetime-1}
		\sup_{x,y\in D} p(1,x,y) \le 	c_1\sup_{x,y\in D} \wt q(1,x,z) \le \sup_{ x\in D} \frac{c_1}{V(x, \phi^{-1}(1))}<\infty.
	\end{align}
Thus, by Lemma \ref{l:largetime-1}, we have
\begin{align}\label{e:largetime-2}
		\left| p(t,x,y) - e^{-\lambda_1 t} \psi_1(x) \wh \psi_1(y) \right|  
		\le c_2e^{-at} 
		\quad \text{for all $t\ge 1$ and $x,y \in D$}.
\end{align}

		\textit{(Upper bound)} We first show that there exists 
		$c_3\ge 1$ such that
	\begin{align}\label{e:largetime-upper-claim}
	  \int_{D\times D} p(t,x,y)\mu(dx) \mu(dy) \le 
	  c_3e^{-\lambda_1 t} \quad \text{for all $t\ge 1$}.
	\end{align}
	Indeed, by using \eqref{e:largetime-2} and the Cauchy-Schwarz inequality, we have for all $t\ge 1$,
	\begin{align*}
		 \int_{D\times D} p(t,x,y)\mu(dx) \mu(dy)&\le e^{-\lambda_1t}  \int_{D} \psi_1(x)\mu(dx) \int_D \wh \psi_1(y) \mu(dy) + 
		 c_2\mu(D)^2 e^{-at} \\
		 &\le \mu(D) \lVert \wh \psi_1 \rVert_{L^2(D)}e^{-\lambda_1t} + 
		 c_2\mu(D)^2 e^{-at}.
	\end{align*}
	Applying Theorem \ref{t:factorization-heatkernel}(i) and using \eqref{e:largetime-1}, we get that for all $x,y \in D$,
\begin{align}\label{e:largetime-upper-3}
	p(1,x,y) \le 
	c_4  \P_x( \zeta> 1) \wh \P_y( \wh \zeta> 1)  \le c_4 (  \P_x( \zeta> 1) \wedge  \wh \P_y( \wh \zeta> 1) ).
\end{align}
 For all $t\ge 3$ and $x,y \in D$, using the semigroup property in the first equality below,  \eqref{e:largetime-upper-3} in the first inequality  and \eqref{e:largetime-upper-claim} in the second, we arrive at
\begin{align*}
	p(t,x,y)& = \int_{D\times D} p(1,x,z) p(t-2,z,w) p(1,w,y) \mu(dz)\mu(dw)\\
	&\le c_4^2 \P_x( \zeta> 1) \, \wh\P_y( \wh \zeta> 1)
	\int_{D\times D}  p(t-2,z,w) \mu(dz)\mu(dw)\\
	& \le  c_5  \P_x( \zeta> 1) \, 
	\wh\P_y( \wh \zeta> 1) e^{-\lambda_1t}.
\end{align*}
	
	\textit{(Lower bound)} Fix a compact set $K \subset D$ and let $r:=\text{dist}(K,D^c)>0$. Since $\psi_1$ and 	$\wh \psi_1$
	 are strictly positive and continuous on $D$, we have 
	 $c_6:=\min_K \psi_1 \wedge \min_K \wh \psi_1>0$.
	 Thus, by \eqref{e:largetime-2}, there exists $T \ge 3$ such that
	\begin{align}\label{e:largetime-lower-1}
		p(t,x,y) \ge  2^{-1}
		c_6^2 e^{-\lambda_1t} 
		\quad \text{for all $(t,x,y) \in [T-2,\infty) \times K \times K$}.
	\end{align}
	For any $t \in [1,T]$, by Lemma \ref{l:infimum-volume},  we have
	\begin{align}\label{e:largetime-lower-2}
			\inf_{x,y\in D} \wt q(t,x,y)&\ge \inf_{x\in D} \frac{1}{V(x, \phi^{-1}(T))} \wedge  \frac{1}{V(x, \diam(D))\phi(\diam(D))}=:
			c_7 >0.
	\end{align}
	Hence, for $t \in [3,T]$, the desired lower bound follows from Theorem \ref{t:factorization-heatkernel}(i).

	For all $z\in K$, we have $\delta_D(z) \ge r$. Thus, by applying Lemma \ref{l:survival-interior}(i) with $a=1 \wedge \frac{r}{\phi^{-1}(1)}$, we have for all $z \in K$,
	\begin{align}\label{e:largetime-lower-3}
  \P_z( \zeta>1) \wedge    \wh\P_z( \wh\zeta>1) \ge  \P_z( \tau_{B(z,r)}>1) \wedge    \wh\P_z(\wh\tau_{B(z, r)}>1)\ge 
  c_8.
	\end{align} 
		For $t \ge T$ and $x,y \in D$, using the semigroup property in the first inequality below, Theorem \ref{t:factorization-heatkernel}(i) 		 and \eqref{e:largetime-lower-1}  in the second, and \eqref{e:largetime-lower-2}  and \eqref{e:largetime-lower-3} in the third,  we obtain
	\begin{align*}
	&p(t,x,y)\ge \int_{K\times K} p(1,x,z) p(t-2,z,w) p(1,w,y) \mu(dz)\mu(dw)\\
		& \ge 
		c_9
		\P_x( \zeta> 1) \, \wh\P_y( \wh \zeta> 1)  e^{-\lambda_1 (t-2)}\int_{K\times K} \wt q(1,x,z)     \wh\P_z( \wh\zeta>1)  \P_w( \zeta>1)   \wt q(1,w,y) \mu(dz)\mu(dw)\\
		&\ge e^{2\lambda_1} 
		c_7^2c_8^2 c_{9}
		\mu(K)^2\P_x( \zeta> 1) \, \wh\P_y( \wh \zeta> 1) e^{-\lambda_1t}.
	\end{align*}
The proof is complete. \qed

	The process $X$ is  called \textit{conservative} if $\P_x(\zeta=\infty)=1$ for all $x\in D$.

	\begin{remark}\label{r:nonconservative}
		\rm We did not assume  $X$ to be non-conservative in Theorem \ref{t:largetime-bounded}.  Under the setting of Theorem \ref{t:largetime-bounded}, it is easy to see that $X$ is non-conservative if and only if  $\lambda_1>0$.
	\end{remark}
	\subsection{Approximate factorization of the Green function}
	
Note that	$X$ is irreducible since, by  \eqref{e:J-density}, we have $J(x,y)>0$ for all $x,y \in D$. Thus, $X$ is either recurrent or transient. When $X$ is transient,  the Green function of $X$ is defined by
	\begin{align*}
		G(x,y):= \int_0^\infty p(t,x,y) dt, \quad x,y \in D.
	\end{align*}

The following hypothesis will be used to establish an approximate factorization for $G(x,y)$. 

	\medskip
	
	\setlength{\leftskip}{0.2in}	
	
	\noindent	{\bf (G)}: \textit{There exists $C\ge 1$ such that 
	\begin{equation}\label{e:Green-condition}
	\int_s^{2\phi(\diam(D))} \frac{dt}{V(x, \phi^{-1}(t))} \le \frac{Cs}{V(x, \phi^{-1}(s))} \quad \text{for all $x\in D$ and $0<s<\phi(\diam(D))$}.
	\end{equation}
	When $D$ is bounded, we also assume that $X$ is not conservative.}

	\setlength{\leftskip}{0in}

	\begin{remark}
		\rm (i) 		
		Let $\beta_1$ and $\alpha_2$  be the constants in \eqref{e:RVD} and \eqref{e:scale-phi}. If $\beta_1>\alpha_2$, 
		then \eqref{e:Green-condition} holds. Indeed,  for all $x \in D$ and $s\in (0,\diam(D))$,  by  \eqref{e:RVD} and \eqref{e:scale-phi}, we have
			\begin{align*}
			\int_s^{2\phi(\diam(D))} \frac{dt}{V(x, \phi^{-1}(t))}  \le \frac{c_1}{V(x, \phi^{-1}(s))}	\int_s^{2\phi(\diam(D))}  
			\left( \frac{s}{t} \right)^{\beta_1/\alpha_2}
			dt \le \frac{c_2s}{V(x, \phi^{-1}(s))}.
		\end{align*}

		\noindent (ii) If $D$ is bounded and $X$ is conservative, then $X$ must be recurrent  since it is irreducible. 
	\end{remark}

	Define
	\begin{align*}
		\wt G(x,y):= \frac{\phi(d(x,y))}{V(x, d(x,y))}, \quad x,y \in D.
	\end{align*}

	\begin{thm}\label{t:factorization-Green}
		Let $D$ be an $\eta$-fat set with characteristics $(\overline R_1, \eta)$.
		Suppose that {\bf (G)} holds, and that either (1)  $D$ is relatively compact  and
		$\kappa \in \mathbf{K}_1(D)$, or  (2) $D$ is unbounded, $\kappa \in \mathbf{K}_\infty(D)$ and $T_0=\overline R_1=\infty$.
	Then there exists $C \ge 1$ such that for all $x,y \in D$,
		\begin{align*}
		&	C^{-1} \P_x\big(\zeta>\phi(d(x,y))\big) \wh  \P_y\big( \wh \zeta>\phi(d(x,y))\big) \wt G(x,y) \\
		&\le G(x,y)\le C \P_x\big(\zeta>\phi(d(x,y))\big) \wh  \P_y\big( \wh \zeta>\phi(d(x,y))\big)\wt G(x,y).
		\end{align*}
	\end{thm}
	\pf Let $x,y \in D$ and define $r=d(x,y)$.  Applying Theorem \ref{t:factorization-heatkernel} (with $T=\phi(2\diam(D))$ if $\diam(D)<\infty$), we obtain for all $t \in [2^{-1}\phi(r),\phi(r)]$,
	\begin{align*}
		p(t,x,y) \ge c_1 \P_x( \zeta>t)\wh\P_y( \wh\zeta>t) \wt q(t,x,y) \ge  \frac{c_1 t\,\P_x( \zeta>\phi(r))\wh\P_y( \wh\zeta>\phi(r))}{V(x,r)\phi(r)}.
	\end{align*}
	It follows that
	\begin{align*}
		G(x,y) &\ge \int_{2^{-1}\phi(r)}^{\phi(r)} p(t,x,y) dt \ge \frac{c_1 \P_x( \zeta> \phi(r))\wh\P_y( \wh\zeta> \phi(r))}{V(x,r) \phi(r) }\int_{2^{-1}\phi(r)}^{\phi(r)}t\, dt\\
		&= 	c_2 \P_x( \zeta> \phi(r))\wh\P_y( \wh\zeta> \phi(r)) \wt G(x,y).
	\end{align*}

	We now prove the upper bound. Applying Theorems \ref{t:factorization-heatkernel}  and \ref{t:largetime-bounded}, we get
		\begin{align*}
		G(x,y) &\le c_3\int_{0}^{\phi(r)}  \frac{t\,  \P_x( \zeta>t)\wh\P_y( \wh\zeta>t)}{V(x,r) \phi(r)}dt  +  c_3\int_{\phi(r)}^{\phi(\diam(D))}   \frac{\P_x( \zeta>t)\wh\P_y( \wh\zeta>t)}{V(x,\phi^{-1}(t)) }dt\\
		&\qquad +  c_3\ind_{(0,\infty)}(\diam(D))\,\P_x( \zeta>1)\wh\P_y( \wh\zeta>1) \int_{\phi(\diam(D))} ^\infty e^{-\lambda_1 t}dt\\
		& =:c_3(I_1+I_2+I_3),
	\end{align*}
	where $-\lambda_1=\sup \text{Re}\, \sigma(\sA)$.	By Proposition \ref{p:boundary-term-scale}, for some $\theta \in [0,1)$, we have
	\begin{align*}
		I_1 &\le c_4 \P_x( \zeta>\phi(r))\wh\P_y( \wh\zeta>\phi(r)) \int_{0}^{\phi(r)}  \frac{t^{1-2\theta}}{V(x,r) \phi(r)^{1-2\theta}}dt =  c_5\P_x( \zeta>\phi(r))\wh\P_y( \wh\zeta>\phi(r)) \wt G(x,y).
	\end{align*}
	For $I_2$, using   \eqref{e:Green-condition},  we obtain
	\begin{align*}
		I_2&\le \P_x( \zeta>\phi(r))\wh\P_y( \wh\zeta>\phi(r)) \int_{\phi(r)}^{\phi(\diam(D))}   \frac{dt}{V(x,\phi^{-1}(t)) }\le c_6\P_x( \zeta>\phi(r))\wh\P_y( \wh\zeta>\phi(r)) \wt G(x,y).
	\end{align*}
	For $I_3$, suppose that $D$ is bounded.  By {\bf (G)},  
	$X$ is not conservative
	 so that $\lambda_1>0$. Moreover, by \eqref{e:Green-condition}, we have
	\begin{align*}
		\wt G(x,y) \ge c_7 
		\int_{\phi(r)}^{2\phi(\diam(D))} \frac{dt}{V(x, \phi^{-1}(t))} \ge \frac{c_7}{\mu(D)}
		\int_{\phi(\diam(D))}^{2\phi(\diam(D))}dt =:c_8>0.
	\end{align*}
	Using Proposition \ref{p:boundary-term-scale}, we see  $
	\P_x(\zeta>\phi(r)) \ge 
	\P_x(\zeta>\phi(\diam(D))) \ge c_9 \P_x(\zeta>1)$. Similarly, we also get $\wh\P_y(\wh\zeta>\phi(r)) \ge c_9 \wh \P_y(\wh\zeta>1)$. Consequently, we arrive at
	\begin{align*}
		I_3& \le c_{10} \P_x( \zeta> \phi(r))\wh\P_y( \wh\zeta> \phi(r)) \int_{\phi(\diam(D))} ^\infty e^{-\lambda_1 t}dt\\
		& = c_{10}\P_x( \zeta> \phi(r))\wh\P_y( \wh\zeta> \phi(r)) \le c_8^{-1}c_{10} \P_x( \zeta> \phi(r))\wh\P_y( \wh\zeta> \phi(r)) \wt G(x,y).
	\end{align*}
	The proof is complete. \qed
	
	We say that a function $h:(0,\infty) \times D\to [0,1]$ is a \textit{regular boundary function} if it satisfies the following two properties: 
	\begin{enumerate}[(i)]
		\item (Doubling property)  
		\emph{For any $x\in D$, $t\mapsto h(t,x)$ is non-increasing and there exists $c_1>1$ such that}
		\begin{align}\label{e:boundary-function-scaling}
			h(t,x) \le c_1h(2t,x) 
			\quad \text{for all $(t, x)\in (0, \infty)\times D$}.
		\end{align}
		
			\item (Interior non-degeneracy)  
			\emph{There exists $c_2\in (0,1]$ such that}
			\begin{align}\label{e:boundary-function-interior}
			\text{$h(t,x) \ge c_2$\quad for all $(t,x) \in (0,\infty) \times D$ with $\delta_D(x) \ge \phi^{-1}(t)$.}
		\end{align}
	\end{enumerate}

	Typical examples of regular boundary functions include $h(t,x) = (1 \wedge (\phi(\delta_D(x))/ t))^p$ for  $p\ge 0$.

	For given regular boundary functions $h_1$ and $h_2$,	consider the following hypotheses:
	
	\smallskip

	\setlength{\leftskip}{0.2in}	
	
	\noindent 	{\rm \bf (UHK$_{h_1,h_2}$)} 
	\emph{There exists $C>0$ such that}
	\begin{align}\label{e:UHK-1}
		p(t,x,y) \le C h_1( t,x)  h_2(t,y) \wt q(t,x,y) \quad \text{for all $(t,x,y) \in (0, 2\phi(\diam(D))) \times D\times D$}.
	\end{align}
	\emph{If $D$ is bounded, 	it further holds that}
	\begin{align}\label{e:UHK-2}
		p(t,x,y) \le C h_1(1,x)  h_2(1,y)e^{-\lambda_1t } \quad \text{for all $(t,x,y) \in [3,\infty) \times D\times D$},
		\end{align}
		\emph{where $-\lambda_1=\sup \text{\rm Re\,}\sigma(\sA)$.}

\vspace{.1in}
	
	\noindent 	{\rm \bf (LHK$_{h_1,h_2}$)} 	\emph{The opposite inequalities to \eqref{e:UHK-1} and \eqref{e:UHK-2} hold.}
	
	\vspace{.1in}

	\noindent 	{\rm \bf (UG$_{h_1,h_2}$)} 
	\emph{There exists $C>0$ such that}
	\begin{align}\label{e:UG-1}
		G(x,y) \le C h_1(\phi(d(x,y)),x) h_2(\phi(d(x,y)),x)  \wt G(x,y)\quad \text{for all $(x,y) \in  D\times D$}.
	\end{align}

	\noindent 	{\rm \bf (LG$_{h_1,h_2}$)} 	
	\emph{The opposite inequality to  \eqref{e:UG-1} holds.}

	\setlength{\leftskip}{0in}
	
	\medskip
	
	As an application of approximate factorizations of heat kernels and Green function, we establish the following equivalence between heat kernel estimates and Green function estimates.

		\begin{cor}\label{c:factorization-equivalence}
		Let $D$ be an $\eta$-fat set with characteristics $(\overline R_1, \eta)$ and let 		$h_1,h_2:(0,\infty) \times D \to [0,1]$ 
		be regular boundary functions.	Suppose that {\bf (G)} holds, and that either (1)  $D$ is relatively compact  and
		$\kappa \in \mathbf{K}_1(D)$, or  (2) $D$ is unbounded, $\kappa \in \mathbf{K}_\infty(D)$ and $T_0=\overline R_1=\infty$.
		Then  we have
		\begin{align*}
			\text{\rm \bf (UHK$_{h_1,h_2}$)} \; \Leftrightarrow \; \text{\rm \bf (UG$_{h_1,h_2}$)}  \quad \text{ and } \quad 	\text{\rm \bf (LHK$_{h_1,h_2}$)} \; \Leftrightarrow \; \text{\rm \bf (LG$_{h_1,h_2}$)} .
		\end{align*}
	\end{cor}
	\pf Since the proofs are similar, we only prove the case where $D$ is bounded and  $\kappa \in \mathbf{K}_1(D)$.
	
	 If {\bf (UHK$_{h_1,h_2}$)} (resp. {\bf (LHK$_{h_1,h_2}$)}) holds, then by Theorem \ref{t:factorization-heatkernel}, we have
	\begin{align*}
		\P_x(\zeta>t) \wh \P_y(\wh \zeta>t) \le c_1h_1(t,x)h_2(t,y)  \quad \left(\text{resp. }	\P_x(\zeta>t) \wh \P_y(\wh \zeta>t) \ge c_1h_1(t,x)h_2(t,y)\right)
	\end{align*}
for all $(t,x,y) \in (0, 2\phi(\diam(D))) \times D\times D$.	Thus, by Theorem \ref{t:factorization-Green}, {\bf (UG$_{h_1,h_2}$)} (resp. {\bf (LG$_{h_1,h_2}$)}) follows.

Now we prove the other direction. Recall that $R_1=(\overline R_1/k_0) \wedge \phi^{-1}(T_0)$. By Proposition \ref{p:boundary-term-scale} and \eqref{e:boundary-function-scaling}, for any $a>1$, there exists $c_2=c_2(a)>1$ such that 	
for all $x,y \in D$ and $t\in (0, 3+ \phi(\text{diam}(D))]$,
	\begin{align}\label{e:equivalence-scaling}
1\le 	\frac{\P_x(\zeta>t) \wh \P_y(\wh \zeta>t)}{\P_x(\zeta>at) \wh \P_y(\wh \zeta>at)} \vee \frac{h_1(t,x)h_2(t,y)}{h_1(at,x) h_2(at,x)}\le c_2.
	\end{align}
 Hence, by Theorems \ref{t:factorization-heatkernel} and  \ref{t:largetime-bounded}, to prove {\bf (UHK$_{h_1,h_2}$)} (resp. {\bf (LHK$_{h_1,h_2}$)}),  it suffices to show that there exists 
 $c_3>0$ such that
 for all  $(t,x) \in (0, \phi(R_1)] \times D$,
	\begin{align}
	&\quad 	\P_x(\zeta>t) \le 
	c_3 h_1(t,x) \quad \text{and} \quad \wh \P_x(\wh \zeta>t) \le c_3h_2(t,x) 
	\label{e:equivalence-claim}\\
	&	 \left(\text{resp. } 	\P_x(\zeta>t) \ge 
	c_3 h_1(t,x) \quad \text{and} \quad \wh \P_x(\wh \zeta>t) \ge c_3h_2(t,x)\right).
	\label{e:equivalence-claim2}
	\end{align}

	Let $(t,x) \in (0, \phi(R_1)] \times D$. 	 Since $D$ is $\eta$-fat, there is $z_0 \in D$ such that $B(z_0, \eta \phi^{-1}(t) )\subset B_D(x, \phi^{-1}(t))$.  Further, by \eqref{e:uniformlyperfect}, there is $z_1 \in B(z_0, \eta  \phi^{-1}(t)/2) \setminus B(z_0, \eta \phi^{-1}(t)/(2k_0))$. Let $w=z_0$  if $d(x,z_0)> \eta \phi^{-1}(t)/(4k_0)$,  and let  $w=z_1$ if  $d(x,z_0)\le \eta \phi^{-1}(t)/(4k_0)$. Then we have 
	\begin{align*}
		 \eta \phi^{-1}(t)/(4k_0) \le d(x,w)  \le \phi^{-1}(t) \quad \text{and} \quad   \delta_D(w) \ge  \eta \phi^{-1}(t)/(4k_0).
	\end{align*} 
	By Lemma \ref{l:survival-interior}(i), we get 
	\begin{align}\label{e:equivalence-inward}
		\P_{w} (\zeta>t) \wedge \wh \P_{w}(\wh \zeta>t) \ge  	\P_{w} (\tau_{B(w,  \eta \phi^{-1}(t)/(4k_0))}>t) \wedge \wh \P_{w}(\wh \tau_{B(w, \eta \phi^{-1}(t)/(4k_0))}>t) \ge 
		c_4.
	\end{align} 
	Note that, by \eqref{e:scale-phi}, it holds that  
	$\phi(\eta \phi^{-1}(t)/(4k_0)) \ge c_5 t$. 
	Hence, by 
\eqref{e:boundary-function-scaling} and	\eqref{e:boundary-function-interior},   
		\begin{align}\label{e:equivalence-inward-2}
		h_1(t,w) \wedge h_2(t,w) \ge 
		c_6\left( 	h_1(c_4t,w) \wedge h_2(c_4t,w)\right)  \ge c_7.
	\end{align} 
	If  {\bf (UG$_{h_1,h_2}$)}   holds, then by \eqref{e:equivalence-inward}, Theorem \ref{t:factorization-Green} and  \eqref{e:equivalence-scaling}, we obtain
		\begin{align*}
	\P_x(\zeta>t) & \le 
	c_4^{-1}\P_x(\zeta>t) \wh\P_{w}(\wh \zeta>t)  \le  c_8h_1(t,x) h_2(t,w)  \le c_8h_1(t,x),\\
	\wh\P_x(\wh \zeta>t) & \le 
	c_4^{-1}\P_w(\zeta>t) \wh\P_{x}(\wh \zeta>t)  \le  c_8h_1(t,w) h_2(t,x)  \le c_8h_2(t,x).
	\end{align*}
	Thus \eqref{e:equivalence-claim} holds.
	
	If {\bf (LG$_{h_1,h_2}$)}   holds, then by Theorem \ref{t:factorization-Green},  \eqref{e:equivalence-scaling} and \eqref{e:equivalence-inward-2}, we get
\begin{align*}
	\P_x(\zeta>t) & \ge \P_x(\zeta>t) \wh\P_{w}(\wh \zeta>t)  \ge 
	c_9h_1(t,x) h_2(t,w) \ge c_7c_9h_1(t,x),\\
	\wh\P_x(\wh \zeta>t) & \ge \P_w(\zeta>t) \wh\P_{x}(\wh \zeta>t)  \ge  
	c_9h_1(t,w) h_2(t,x) \ge c_7c_9 h_2(t,x).
\end{align*}
Thus, \eqref{e:equivalence-claim2} holds and the proof is complete.
\qed 
	
	\section{Applications}\label{s:appl}

A non-decreasing function $\ell:(0,1]\to (0,\infty)$ is said to be a \textit{Dini function} if 
\begin{align*}
	\int_0^1 \frac{\ell(r)}{r} dr <\infty,
\end{align*}
and is 
said to be a \textit{double Dini function} if 
\begin{align*}
\int_0^1 \frac{1}{R}	\int_0^R \frac{\ell(r)}{r} \,dr dR = \int_0^1 \frac{\ell(r)|\log r|}{r}dr  <\infty.
\end{align*}

Denote by $\Dini$  the family of all Dini functions and by $\DDini$ the family of all  double Dini functions.

\begin{defn}
	\rm	Let  $\ell \in \Dini$.  An open set $D\subset \R^d$, $d\ge 2$, is said to be $C^{1,\ell}$ if there exists a localization constant  $r_0\in (0,1]$  such that for any $Q \in \partial D$, there exist a $C^{1}$-function $F_Q:\R^{d-1}\to \R$ satisfying $F_Q(0)=0$, $\nabla F_Q(0)=(0,\cdots,0)$, $ |\nabla F_Q(\wt x) - \nabla F_Q(\wt y) | \le \ell(|\wt x-\wt y|)$ for all $\wt x, \wt y\in \R^{d-1}$, and an orthonormal coordinate system CS$_Q:(\wt y, y_d)= (y_1,\cdots,y_{d-1},y_d)$ with origin at $Q$ such that
	\begin{align*}
	D\cap	B(Q,r_0) = \left\{ y=(\wt y, y_d) \in B(0,r_0) \text{ in CS$_Q$} \, : \, y_d>F_Q(\wt y)\right\}.
	\end{align*}
\end{defn}

\begin{defn}\label{d:Dini}
\rm	An open set $D\subset \R^d$, $d\ge 2$, is said to be 
 $C^{1,{\rm Dini}}$ (resp. $C^{1,2{\text - \rm Dini}}$) if   $D$ is $C^{1,\ell}$ for some $\ell \in \Dini$ (resp. $\ell \in \DDini$) with a localization constant $r_0$.	The pair $(r_0,\ell)$ is called the characteristics of $D$.
\end{defn}

A $C^{1,\ell}$ open set $D$ with $\ell(r) = cr^\eps$ for some  $\eps \in (0,1]$ is referred to as a  $C^{1,\eps}$ open set.
 
 \begin{remark}
 	\rm 
	If $\eps \in (0,1]$, then every $C^{1,\eps}$ set is $C^{1,{\rm Dini}}$. 
 	An example of a $C^{1,{\rm Dini}}$ open set which is not $C^{1,\eps}$ for any $\eps \in (0,1]$ is given by $D:=\{(x_1,x_2): x_2> x_1 (1+|\log x_1|)^{-k}\}$ for $k>1$.
 \end{remark}

 Let $D\subset \R^d$ be a   $C^{1,{\rm Dini}}$ open set with characteristics $(r_0,\ell)$.    By choosing a smaller $r_0$ if necessary, we assume without loss of generality that  $\ell(r_0)\le 1/4$.
  Recall that $\delta_D(x)$ denotes the distance between $x$ and $D^c$. The signed distance to $\partial D$ is defined by
 \begin{align}\label{e:signed-distance}
	 \wt\delta_D(x) := \begin{cases}
 	 	\delta_D(x) &\mbox{ if $x\in \overline D$},\\
 	 	-\delta_{\overline D^c}(x) &\mbox{ if $x\in \overline D^c$}.
 	 \end{cases}
 \end{align}
  Define	$\wt D(r):= \{ x \in \R^d : |\wt \delta_D(x)| < r\}$ 
  for $r>0$. For any $x\in D$,  we fix a $Q_x\in \partial D$  such that $\delta_D(x) = |x-Q_x|$.  For any $Q\in \partial D$, we denote by $n_Q$ the inward unit normal to $\partial D$ at $Q$.

  Following \cite{Li85}, we construct a  \textit{regularized distance}  $\rho$ for $D$ as follows: 
   Fix a non-negative  function  $\vp \in C_c^\infty(B(0,1))$    satisfying $\int_{B(0,1)} \vp(z)dz=1$. Let 
   $\{B(Q^i,r_0/20)\}_{i\in I}$ be a countable open cover of $\partial D$, with
   $Q^i\in \partial D$ for all $i\in I$, such that 
   $B(Q^i,r_0/100)$, $i\in I$, are pairwise disjoint.        For any $i\in I$, define $h_i:B(Q^i,r_0) \to \R$ by 
  $$
  h_i(y) := 
  y_d - F_{Q^i}(\wt y) \quad \text{if $y=(\wt y, y_d)$ in CS$_{Q^i}$}.
  $$ 
  Then we define a function $P_i(x,a)$ by
 \begin{align*}
	P_i(x,a):= 
	\int_{\R^d} h_i(x- (a/4)z)\,\vp(z)dz, \quad x \in B(Q^i,r_0/4), \;   a \in [ -  r_0,   r_0].
 \end{align*}
 Since $\ell(r_0)\le 1/4$, we have  $|h_i(y)-h_i(z)| \le |y_d-z_d| + \ell(r_0) |\wt y - \wt z| \le  2|y-z|$ for all  $y=(\wt y, y_d),z=(\wt z, z_d) \in B(Q^i, r_0/2)$. Thus, by \cite[p. 339]{Li85},  there exists  $\rho_i:B(Q^i,r_0/4)\to [ -  r_0,   r_0]$ satisfying
\begin{align}\label{e:regular-distance-equation}
	\rho_i(x) = P_i(x, \rho_i(x)) \quad \text{for all $x \in B(Q^i,r_0/4)$}.
\end{align}
Further,  by \cite[line 10 on  p. 331,  the second display on p. 339, (1.8), (1.10), and the arguments on p. 333-334]{Li85} (with $L=4$), there exists $C=C(d)>0$  such that the following hold:
\begin{align}
	\text{$2^{-1}\wt \delta_D(x) \le \rho_i(x)\le  2\sqrt 5\,\wt \delta_D(x)$} \quad &\text{ for all $x\in B(Q^i,r_0/4)$,}\label{e:regular-distance-pre-1}\\
 \text{$|\nabla\rho_i(x) - \nabla \rho_i(y)|\le 8 \ell(|x-y|)$} \quad &\text{ for all $x,y\in B(Q^i,r_0/4)$,} \label{e:regular-distance-pre-2}\\
\text{$|D^2 \rho_i(x)| \le C |\wt \delta_D(x)|^{-1} \ell(|\wt \delta_D(x)|)$} \quad &\text{ for all $x\in B(Q^i,r_0/4) \setminus \partial D$.}\label{e:regular-distance-pre-3}
\end{align}

 For any $x \in \wt D(r_0/7)$, there exists $Q^i$ such that $|Q_x-Q^i|\le r_0/20$.  Hence,  $\wt D(r_0/7) \subset \cup_{i\in I} B(Q^i,r_0/5)$.
Let $\{\eta_0\} \cup \{\wt \eta_0\}\cup \{\eta_i\}_{i\in I}$ be a smooth partition of unity subordinate to $(D\setminus \overline{\wt D(r_0/8)}) \cup    (\overline D^c \setminus \overline{\wt D(r_0/8)}) \cup \{B(Q^i,r_0/5)\}_{i\in I}$ with the property that for any multi-index $\alpha$, there exists $C_\alpha>0$ such that  $|\partial^\alpha \eta_i| \le C_\alpha$ for all $i\in I$. The existence of such a  partition of unity follows from \cite[Lemma A.1.3]{Sh92}.  Define 
$$
\rho(x) := 
r_0\eta_{0}(x)- r_0\wt \eta_{0}(x) + \sum_{i\in I} \eta_i(x)\rho_i(x) , \quad x \in \R^d.
$$
By \eqref{e:regular-distance-pre-1}--\eqref{e:regular-distance-pre-3}, there are comparison constants and  $C>0$ depending only on $d$  such that
	\begin{align}
			\text{$\rho(x)\asymp (\wt \delta_D(x) \vee (-r_0)) \wedge r_0$} \quad &\text{for all $x\in \R^d$,}\label{e:regular-distance-1}\\
		\text{$|\nabla\rho(x) - \nabla \rho(y)|\le C \ell(|x-y|)$} \quad &\text{for all $x,y\in \R^d$,} \label{e:regular-distance-2}\\
		\text{$|D^2 \rho(x)| \le C(1+ |\wt \delta_D(x)|^{-1} \ell(|\wt \delta_D(x)|))$} \quad &\text{for all $x\in \R^d \setminus \partial D$.}\label{e:regular-distance-3}
\end{align}

\begin{lemma}\label{l:inward-normal}
Suppose   $D\subset \R^d$ is a   $C^{1,{\rm Dini}}$ open set with characteristics $(r_0,\ell)$
  such that $\ell(r_0) \le 1/4$.  Let $\rho$ be the regularized distance for $D$ defined  above. Then there exists  $r_1=r_1(d,r_0)\in (0, r_0]$   such that
\begin{align}\label{e:inward-normal-result-1}
	\frac14 \le |\nabla \rho(x)| \le 2  \quad \text{for all $x\in \overline{D(r_1)}$}.
\end{align}
Moreover, it holds that
	\begin{align}\label{e:inward-normal-result-2}
		\frac{\nabla \rho(Q)}{|\nabla \rho(Q)|}  = n_Q \quad \text{for all $Q\in \partial D$}.
	\end{align}
\end{lemma}
\pf    
Let $r_1 \in (0, r_0]$ be a constant to be chosen later, 
let $x\in \overline{D(r_1)}$ and set $I_x:=\{i\in I: x\in B(Q^i, r_0/5)\}$. Pick an arbitrary $i \in I_x$ and write $x=(\wt x, x_d)$ in CS$_{Q^i}$. Let $\nabla_z P_i(z,a)$ denote the vector $(\partial_{z_1} P_i(z,a), \cdots, \partial_{z_d} P_i(z,a))$. By \eqref{e:regular-distance-equation}, it holds that 
\begin{align}\label{e:inward-normal-1}
	\nabla \rho_i(x)&= \nabla_z P_i(x, \rho_i(x))   + \partial_a  P_i(x, a)|_{a=\rho_i(x)}	\nabla \rho_i(x).
\end{align} 
 Since  $|\nabla F_{Q^i}(\wt y) | \le \ell(|\wt y|)$ for all $ \wt y\in \R^{d-1}$ and  $\rho_i(x) \le c_1 \delta_D(x) \le c_1r_1$, 
 letting $r_1 \le 16r_0/(5c_1)$, 
 we obtain for all $z=(\wt z, z_d) \in B(0,1)$,
\begin{align*}
&	\big| (\nabla h_i)(x- (\rho_i(x)/4)z) -  (\wt 0,1)\big|   = \big| (\nabla F_{Q^i}) (\wt x - (\rho_i(x)/4)\wt z)\big|\le  \ell (|x-Q^i| + c_1r_1/4)  \le \ell(r_0) \le \frac14.
\end{align*}
Hence, using  $\int_{B(0,1)}\vp(z)dz =1$,  we get
\begin{align}\label{e:inward-normal-2}
	\left| \partial_a  P_i(x, a)|_{a=\rho_i(x)}\right|   = \bigg|\frac{1}{4} \int_{B(0,1)} (\nabla h_i) (x- (\rho_i(x)/4)z) \cdot z \vp(z) dz  \bigg| \le \frac{5}{16}\int_{B(0,1)} |z| \vp(z) dz  < \frac13
\end{align}
and
\begin{align}\label{e:inward-normal-3}
	\left|\nabla_z  P_i(x, \rho_i(x)) - (\wt 0, 1)\right| = 	\int_{B(0,1)} \left| (\nabla h_i)(x- (\rho_i(x)/4)z) - (\wt 0, 1)\right| \vp(z)dz \le \frac14.
\end{align}
By \eqref{e:inward-normal-1}, \eqref{e:inward-normal-2} and \eqref{e:inward-normal-3}, there exists a  constant $a_i\in [3/4,3/2]$ such that 
\begin{align}\label{e:inward-normal-6}
	|\nabla \rho_i(x) - (\wt 0, a_i)|  = a_i 	|\nabla_z P_i(x, \rho_i(x)) -(\wt 0, 1)| \le \frac38.
\end{align}
Using \eqref{e:inward-normal-6} and the fact that $\rho_i(x) \le c_1r_1$ for all $i \in I_x$, and choosing $r_1$ small enough,  we get
\begin{align*}
\left|\nabla \rho (x)-\sum_{i\in I_x}(\wt 0,a_i\eta_i(x))  \right|&=\left|  \sum_{i\in I_x}  \left( \eta_i(x) ( \nabla \rho_i(x)- (\wt 0, a_i)) + \rho_i(x) \nabla \eta_i(x) \right)\right| \\
& \le   \frac38 \sum_{i\in I_x} \eta_i(x) + c_1r_1 \sum_{i\in I_x} |  \nabla \eta_i(x)| \le \frac38 + c_2r_1 \le \frac12.
\end{align*} 
Since $\sum_{i\in I_x}a_i\eta_i(x)   \in [ \min_{i\in I_x} a_i,\max_{i\in I_x} a_i ] \subset [3/4,3/2]$, this proves \eqref{e:inward-normal-result-1}.

For \eqref{e:inward-normal-result-2}, we let $x\in\partial D$. For any $i\in I_x$, we have $\rho_i(x)=0$ and
 $  \nabla_z P_i(x, 0)= \nabla h_i(x) = ( -\partial_{y_1} F_{Q^i}(\wt x), \cdots, -\partial_{y_{d-1}} F_{Q^i}(\wt x),1) = |\nabla h_i(x)|\,n_x$. Thus, using \eqref{e:inward-normal-1} and \eqref{e:inward-normal-2}, we obtain 
\begin{align*}
	\frac{	\nabla \rho_i(x)}{|	\nabla \rho_i(x)|} = \frac{\nabla_z P_i(x, 0)   }{|\nabla_z P_i(x, 0)   |}  = n_x.
\end{align*} 
It follows that
\begin{align*}
	\frac{\nabla \rho(x)}{|\nabla \rho(x)|}  = \frac{\sum_{i\in I_x} \eta_i(x) \nabla \rho_i(x)}{|\sum_{i\in I_x} \eta_i(x) \nabla \rho_i(x)|}  = \frac{\sum_{i\in I_x} \eta_i(x) |\nabla \rho_i(x)| n_x}{\sum_{i\in I_x} \eta_i(x) |\nabla \rho_i(x)|} =n_x.
\end{align*}  
The proof is complete. \qed

 \begin{lemma}\label{l:censored-domain-closedness}
Suppose $D\subset \R^d$ is a   $C^{1,{\rm Dini}}$ open set with characteristics $(r_0,\ell)$ such that $\ell(r_0)\le 1/4$.
  Let $\rho $ be the regularized distance for $D$  defined at the beginning of this section and let $r_1$ be the constant from Lemma \ref{l:inward-normal}. Let $x_0 \in D(r_1)$ and define
 	$$
	E := 
	\big\{ y \in \R^d: \rho(x_0) + \nabla \rho(x_0)\cdot (y-x_0) > 0 \big\}.
 	$$
 	 Suppose that there exist $c_0\ge 1$ and $\eps>0$ such that
 \begin{align}\label{e:ell-scaling}	\frac{\ell(r)}{\ell(s)} \le c_0\bigg( \frac{r}{s}\bigg)^{\eps} \quad \text{for all $0<s\le r\le 1$}.
 \end{align}
 Then the following two assertions hold.
 
 \smallskip
 
 \noindent (i) There exists $C=C(d)\ge 1$ such that
 \begin{align*}
C^{-1} \delta_D(x_0) \le \delta_E(x_0) \le C\delta_D(x_0) \quad \text{and} \quad  	|\delta_D(x_0)-\delta_E(x_0)| \le C \delta_D(x_0) \ell(\delta_D(x_0)).
 \end{align*}
 
 \noindent (ii)
  There exists $C=C(d,r_0,c_0,\ell(1))>0$ such that for any $a>\eps$,
 	\begin{align*}
 		\int_{(D\setminus E) \cup (E\setminus D)} \frac{dy}{|x_0-y|^{d+a}} \le \frac{C\delta_D(x_0)^{-a}\ell(\delta_D(x_0))}{a-\eps}.
 	\end{align*}
 \end{lemma}
 \pf We assume, without loss of generality, that $Q_{x_0}=0$ and 
 use the coordinate system CS$_{Q_{x_0}}$. Set $r:=\delta_D(x_0)$ so that $x_0 = (\wt 0, r)$. Note that $\delta_E(x_0)=\rho(x_0)/|\nabla \rho(x_0)|$. Define
 \begin{align*}
 	\theta(y)&:= \frac{\nabla \rho(y)}{|\nabla \rho(y)|}  \quad \text{for $y\in \overline{D(r_1)}$}, \qquad \theta_0:=\theta(x_0), \qquad z_0:=x_0 -\delta_E(x_0) \theta_0.
 \end{align*} 
 
 (i) 
 By \eqref{e:regular-distance-1} and  \eqref{e:inward-normal-result-1},  we have $\delta_E(x_0) \asymp r$. Moreover, since   $\nabla \rho(0)=  (\wt 0,|\nabla \rho(0)|)$ and $\rho(x_0) \le c\delta_D(x_0)=cr$, using \eqref{e:regular-distance-2} and  the monotone and scaling properties of $\ell$, we obtain
 \begin{align*}
 	\left|\delta_E(x_0)-r \right| &\le 	\left|\frac{\rho(x_0)}{|\nabla \rho(x_0)|} -\frac{\rho(x_0)}{|\nabla \rho(0)|} \right| + 	\left|\frac{\rho(x_0)}{|\nabla \rho(0)|} -r \right|   \\
 	& =  	\rho(x_0)\left|\frac{1}{|\nabla \rho(x_0)|} -\frac{1}{|\nabla \rho(0)|} \right| +  \frac{1}{|\nabla \rho(0)|}\left|  \int_0^{r}  \left( \nabla \rho((\wt 0, s)) \cdot (\wt 0,1) - \nabla \rho(0) \cdot (\wt 0,1) \right) ds \right| \\
 	&\le  c_1 \rho(x_0) \ell(r)  + c_2 \int_0^r \ell(s) ds \le  c_3 r \ell(r).
 \end{align*}
 
 (ii) By Lemma \ref{l:inward-normal}, it holds that $\theta(0)=(\wt 0,1)$. 
 Hence, by \eqref{e:regular-distance-2}, \eqref{e:inward-normal-result-1} and \eqref{e:inward-normal-result-2}, we get 
 \begin{align*}
 	\theta_0 = \frac{\nabla \rho(0)}{|\nabla \rho(0)|}  + \xi_0 = (\wt 0,1) + \xi_0 \quad \text{for some $|\xi_0|\le c_4\ell (r)$}.
 \end{align*}
  Note that $y\in E$ if and only if $\theta_0 \cdot (y-z_0)>0$. Further,  for any $y = (\wt y, y_d) \in B(0,r_0)$, we have $|F_{Q_{0}}(\wt y)| \le |\wt y| \sup_{|\wt z| \le |\wt y|} |\nabla F_{Q_{0}}(\wt z)| \le  
 |\wt y|\ell(|\wt y|) $. Define
 \begin{align*}
 U_1&:=	\left\{ y=(\wt y, y_d) \, : \, |\wt y|<r_0,\,  y_d>-|\wt y| \ell(|\wt y|),\; \theta_0\cdot (y-z_0)\le 0\right\},\\
 	U_2&:=	\left\{ y=(\wt y, y_d) \,  : \, |\wt y|<r_0,\,  y_d\le |\wt y| \ell(|\wt y|) , \theta_0\cdot (y-z_0)> 0\right\}.
 \end{align*}
Suppose $y=(\wt y,y_d)\in U_1$ with  $y_d \ge 0$. It follows from (i) that 
 \begin{align*}
 	|z_0| \le |x_0 - r\theta_0| + 	\left|\delta_E(x_0)-r \right|  = r |\xi_0| + c_3r\ell(r) \le c_5 r \ell(r).
 \end{align*} Thus
we have
\begin{align}\label{e:censored-domain-closedness-1}
0 &\ge \theta_0 \cdot (y-z_0) \ge y_d - |\xi_0||y|  - |z_0|  \ge y_d - c_4 \ell(r) y_d -  c_4 \ell(r)  |\wt y|   - c_5 r \ell(r).
\end{align}
Similarly, for any $y=(\wt y,y_d)\in U_2$ with $y_d\le 0$,  we get
\begin{align}\label{e:censored-domain-closedness-2}
	0 &< \theta_0 \cdot (y-z_0) \le   y_d - c_4 \ell(r) y_d +  c_4 \ell(r)  |\wt y|   + c_5 r \ell(r).
\end{align}
Let $r_*:=\sup\{s \in (0,1]: 2(c_4+c_5)\ell(s)\le  1/2\}$.  If $r\ge r_*$, then using $\delta_E(x_0) \wedge r \ge c_6r$, we obtain
\begin{align*}
	\int_{(D\setminus E) \cup (E\setminus D)} \frac{dy}{|x_0-y|^{d+a}} \le 	\int_{B(x_0,c_6r)^c} \frac{dy}{|x_0-y|^{d+a}} =c_7r^{-a} \le 4(c_4+c_5)c_7  r^{-a} \ell(r).
\end{align*} 

Suppose $r<r_*$. Observe that
\begin{align*}
	\int_{(D\setminus E) \cup (E\setminus D)}\frac{dy}{|x_0-y|^{d+a}}& \le \int_{y=(\wt y, y_d) \in B(0,r_0) \setminus E \,:\, y_d>F_0(\wt y)} \frac{dy}{|x_0-y|^{d+a}} \\
	&\quad + \int_{y=(\wt y, y_d) \in B(0,r_0) \cap E \,:\, y_d\le F_0(\wt y)} \frac{dy}{|x_0-y|^{d+a}} + \int_{B(0,r_0)^c} \frac{dy}{|x_0-y|^{d+a}}\\
	&\le \int_{U_1} \frac{dy}{|x_0-y|^{d+a}} + \int_{U_2} \frac{dy}{|x_0-y|^{d+a}}+\int_{B(0,r_0)^c} \frac{dy}{|x_0-y|^{d+a}}\\
	&=:I_1+I_2+I_3.
\end{align*}  Since $ 2(c_4+c_5) r\ell(r) \le r/2$, using \eqref{e:censored-domain-closedness-1}, we obtain
\begin{align*}
	I_1 &\le c_8\int_0^{r_0} s^{d-2} \int_{-s \ell(s)}^{ 2c_4 \ell(r) s + 2c_5 r\ell(r) } \frac{dy_d \, ds}{(s^2 +(r-y_d)^2)^{(d+a)/2}}\\
	&\le  c_8\int_0^{r} s^{d-2} \int_{-r \ell (r)}^{ 2(c_4+ c_5) r \ell(r)} \frac{dy_d \, ds}{(r/2)^{d+a}} + c_8\int_{r}^{r_0} s^{d-2} \int_{-s\ell(s)}^{  2(c_4+c_5)s\ell(s)} \frac{dy_d \, ds}{s^{d+a}} \\
	&\le c_9r^{-a} \ell(r) + c_{10}\int_{r}^{r_0} \frac{\ell(s)}{s^{1+a}}ds.
\end{align*}
Similarly, we get from \eqref{e:censored-domain-closedness-2}  that
\begin{align*}
	I_2 		&\le  c_{11}\int_0^{r} s^{d-2} \int^{r\ell(r)}_{-2(c_4+c_5)r\ell(r)} \frac{dy_d \, ds}{(r/2)^{d+a}} + c_{11}\int_{r}^{r_0} s^{d-2} \int^{s\ell(s)}_{-2(c_4+c_5)s\ell(s)} \frac{dy_d \, ds}{s^{d+a}} \\
		&\le c_{12}r^{-a} \ell(r) + c_{13}\int_{r}^{r_0} \frac{\ell(s)}{s^{1+a}}ds.
\end{align*}
Clearly, $I_3=c_{14}$. Using \eqref{e:ell-scaling}, we see $r^{-a} \ell(r)\ge c_0^{-1} \ell(1) = c_{15} \ell(1)I_3$  and
\begin{align*}
	\int_{r}^{r_0} \frac{\ell(s)}{s^{1+a}}ds \le 
\frac{c_{0}\ell(r)}{r^{\eps}}	\int_{r}^{r_0} \frac{ds}{s^{1+a-\eps}} \le \frac{c_{16} r^{-a} \ell(r)}{a-\eps}.
\end{align*}
Combining the estimates for $I_1,I_2$ and $I_3$, we get (ii). The proof is complete. \qed

 \subsection{Stable-like processes with critical killings in $C^{1,{\rm Dini}}$ open sets}\label{ss:critical-killing}\label{ss:5.1}
 
 Let $\alpha\in (0,2)$ and 	$d\ge 2$. Before we present the general setup of this subsection, we first introduce stable-like processes and their killed processes in open sets. Killed stable-like processes will  serve as special cases of our general framework.

  Let $K$ be a symmetric Borel function on $\R^d\times \R^d$ satisfying 
 \begin{align}\label{e:pre-censored-1}
 	C \le 	K(x,y) \le C' \quad 
	\text{for all $x,y \in \R^d$},
 \end{align}
for some constants $C' \ge C>0$.  Consider the bilinear form $(\overline \sE^{ \R^d}, \overline \sF^{ \R^d})$ on  $L^2(\R^d,dx)$ defined by
 \begin{align*}
 	\overline \sE^{\R^d}(u,v) &:= \int_{\R^d\times \R^d\setminus \diag} (u(x)-u(y))(v(x)-v(y)) \frac{K(x,y)}{|x-y|^{d+\alpha}} dxdy, \nn\\
\overline	\sF^{\R^d}&:= \bigg\{ u \in L^2(\R^d,dx): \int_{\R^d\times \R^d\setminus \diag}  \frac{(u(x)-u(y))^2}{|x-y|^{d+\alpha}} dxdy<\infty\bigg\}.
 \end{align*}
 Note that $
 \overline	\sF^{\R^d}$ is the fractional Sobolev space of order $\alpha/2$.
 By  \eqref{e:pre-censored-1},
 $(\overline \sE^{\R^d},\overline \sF^{\R^d})$ is a regular  Dirichlet form on $L^2(\R^d,dx)$.  Let $Z$ be a symmetric Hunt process associated with $(\overline \sE^{\R^d},\overline \sF^{\R^d})$. The  process $Z$ is referred to as an 
 \emph{$\alpha$-stable-like process} on $\R^d$.
 
 Let $D\subset \R^d$ be a proper open subset and let   $\sF^D$ be the closure of $C_c^\infty(D)$ 
 with respect to the 
 $\overline \sE^{\R^d}_1$-norm, where  $\overline \sE^{\R^d}_1(u,u):=\overline \sE^{\R^d}(u,u)+\lVert u \rVert_{2}^2$. Define
  \begin{align}\label{e:killingfunction}
  	\kappa_D(x):= \int_{
  	\R^d \setminus D}  \frac{K(x,y)}{|x-y|^{d+\alpha}}dy, \quad x \in D.  
	\end{align}
	Let $Z^D$ be the killed  process of  $Z$ on $D$.  Then $Z^D$ is a Hunt process on $D$ associated with the regular Dirichlet form  $(\sE^D, \sF^D)$, where
  \begin{align*}
  	\sE^D(u,v)  :=\int_{\R^d\times \R^d\setminus \diag} (u(x)-u(y))(v(x)-v(y)) \frac{K(x,y)}{|x-y|^{d+\alpha}} dxdy + \int_D u(x) v(x) \kappa_D(x) dx. 
  \end{align*}
 Consider the following condition: There exist constants $\theta >(\alpha-1)_+$ and  $C>0$ such that
  \begin{align}\label{e:pre-censored-2}
  	|K(x,x) - K(x,y) | \le 
  	C |x-y|^\theta \quad \text{ for all 
  		$x \in D$ and $y \in \R^d$ with $|x-y|\le 1$.}
  	\end{align}
  	
 For $q\in [0,\alpha)$, define
 \begin{align}\label{e:def-constant-Cq}
 	\sC(d,\alpha,q):=\frac{\omega_{d-1}}{2} \sB\left(\frac{\alpha+1}{2}, \frac{d-1}{2}\right) \int_0^1 \frac{(t^q -1)(1-t^{\alpha-q-1})}{(1-t)^{1+\alpha}} dt,
 \end{align} 
 where $\omega_{d-1}$ is the $(d-2)$-dimensional Lebesgue measure of the unit sphere in $\R^{d-1}$ and $\sB(\cdot,\cdot)$ is the beta function.

\begin{lemma}\label{l:critical-killing-1}  
 Let  $D$ be a $C^{1,{\rm Dini}}$ open set with characteristics $(r_0,\ell_0)$, with $r^{-\alpha/2}\ell_0(r)$ being non-increasing.  Assume \eqref{e:pre-censored-1} and  \eqref{e:pre-censored-2}.    Then  there exists $C>0$ such that for any $x_0\in D$,
	\begin{align*}
		\left|\kappa_D(x_0) - \sC(d,\alpha, \alpha/2) K(x_0,x_0)\delta_D(x_0)^{-\alpha}\right|\le C\left( \delta_D(x_0)^{-\alpha} \big( \delta_D(x_0)^{\theta 
			\wedge (\alpha/2)} + \ell_0(\delta_D(x_0))\big) + 1 \right).
	\end{align*}
\end{lemma}
\pf  
By choosing a smaller $r_0$ if necessary, we assume  that $\ell_0(r_0)\le 1/4$. Let $\rho $ be the regularized distance for $D$  defined at the beginning of this section (with $\ell$ replaced by $\ell_0$) and let $r_1$ be the constant from Lemma \ref{l:inward-normal}. Let $x_0\in D$ and set $r:=\delta_D(x_0)$.  If $r\ge r_1$, then by \eqref{e:pre-censored-1}, 
$$
	\left|\kappa_D(x_0) - \sC(d,\alpha, \alpha/2) K(x_0,x_0)r^{-\alpha}\right| \le c_1 \int_{B(x_0,r)^c} \frac{dy}{|x_0-y|^{d+\alpha}}  + c_2r^{-\alpha}= c_3 r^{-\alpha} \le c_3r_1^{-\alpha}.
$$

Suppose $r<r_1$.  Let 
$	E:= \{ y \in \R^d: \rho(x_0) + \nabla \rho(x_0)\cdot (y-x_0) > 0 \}.$
Using  the fact that $h(y):= (y_d)_+^{\alpha/2}$  satisfies $\Delta^{\alpha/2}h=0$ in $\R^d_+$, we get from \cite[(5.4)]{BBC03} that
\begin{align*}
	\int_{E^c} \frac{dy}{|x_0-y|^{d+\alpha}} = \sC(d,\alpha, \alpha/2)\delta_E(x_0)^{-\alpha}.
\end{align*}
Using this, we obtain
\begin{align*}
		&\left|\kappa_D(x_0) - \sC(d,\alpha, \alpha/2) K(x_0,x_0)r^{-\alpha}\right|\\
		 &\le \left|\int_{D^c} \frac{K(x_0,y)}{|x_0-y|^{d+\alpha}}dy - \int_{E^c}\frac{K(x_0,x_0)}{|x_0-y|^{d+\alpha}}dy\right| + \sC(d,\alpha, \alpha/2) K(x_0,x_0)\left| \delta_E(x_0)^{-\alpha} -r^{-\alpha}\right|\\
	&\le \int_{D^c} \frac{|K(x_0,x_0)-K(x_0,y)|}{|x_0-y|^{d+\alpha}}dy  +  K(x_0,x_0) \int_{(D \setminus E) \cup (E \setminus D)}\frac{dy}{|x_0-y|^{d+\alpha}} \\
	&\quad   + \sC(d,\alpha, \alpha/2) K(x_0,x_0)\left| \delta_E(x_0)^{-\alpha} - r^{-\alpha}\right|\\
	&=:I_1+I_2+I_3.
\end{align*}
By \eqref{e:pre-censored-1} and \eqref{e:pre-censored-2}, we have $|K(x_0,x_0) - K(x_0,y)|\le c_4(|x_0-y|\wedge 1)^\theta
\le c_4(|x_0-y|\wedge 1)^{\theta \wedge (\alpha/2)}
$. It follows that
\begin{align*}
	I_1 \le c_4\int_{B(x_0,r)^c} \frac{(|x_0-y|\wedge 1)^{\theta \wedge (\alpha/2)}}{|x_0-y|^{d+\alpha}}dy = c_5\int_r^1 \frac{ds}{s^{1+\alpha-\theta \wedge (\alpha/2)}}ds + c_5 \int_1^\infty \frac{ds}{s^{1+\alpha}}  \le c_6 r^{\theta \wedge (\alpha/2)-\alpha}. 
\end{align*} 
For $I_2$, by \eqref{e:pre-censored-1} and Lemma \ref{l:censored-domain-closedness}(ii), we get $
I_2 \le  c_7 r^{-\alpha} \ell_0(r)$. For $I_3$, using \eqref{e:pre-censored-1}, 
$\delta_E(x_0) \ge c_8r$
and Lemma \ref{l:censored-domain-closedness}(i), we obtain $I_3\le  
c_9 r^{-\alpha-1} |\delta_E(x_0)-r| \le c_{10} r^{-\alpha}\ell_0(r).$ 
The proof is complete.\qed 

	For every $\ell_0 \in \Dini$, by  Lemma \ref{l:Dini-regularizing}, there exists $\wt \ell_0\in \Dini$ such that $\ell_0 \le \wt \ell_0$ and  $r^{-\alpha/2}\wt\ell_0(r)$ is non-increasing. Thus,  the condition that  $r^{-\alpha/2}\ell_0(r)$ is  non-increasing in Lemma \ref{l:critical-killing-1}   is not restrictive.

\medskip

We now introduce our general framework, which encompasses    killed  stable-like processes.  Let  $D\subset \R^d$ be a  $C^{1,{\rm Dini}}$ open set.  Let $K$ be a symmetric Borel function on $D\times D$ satisfying 
\begin{align}\label{e:censored-1}
	A_1 \le 	K(x,y) \le A_2 \quad 
	\text{for all $x,y \in D$},
\end{align} 
for some  $A_2\ge A_1>0$. 
We assume that there exist constants $\theta >(\alpha-1)_+$ and  
$A_3>0$ such that
\begin{align}\label{e:censored-2}
|K(x,x) - K(x,y) | \le A_3 |x-y|^\theta \quad \text{ for all 
	$x,y \in D$.}
	\end{align}
Define $\sC(d,\alpha,q)$ as \eqref{e:def-constant-Cq}. Then $q \mapsto \sC(d,\alpha,q)$ is strictly increasing on $[(\alpha-1)/2,\alpha)$, 
\begin{align}\label{e:criticall-killing-range}
	\sC(d,\alpha,\alpha-1) = \sC(d,\alpha,0) =0 \quad \text{and} \quad \lim_{q\uparrow \alpha} \sC(d,\alpha,q)=\infty.
\end{align} 
Motivated by killed stable-like processes and Lemma \ref{l:critical-killing-1}, we let
$\sK_\alpha$ be the family of all non-negative Borel functions $\kappa$ on $D$ with the following properties: 
\emph{There exist constants $a_0\ge 0$,
$A_4>0$ and $\ell_1 \in \Dini$ such that $\kappa(x) \le A_4$ for all $x\in D$ 
with $\delta_D(x)\ge 1$, and}
\begin{align}\label{e:critical-class}
	\big|\kappa(x) - a_0K(x,x) \delta_D(x)^{-\alpha}\big| \le 
	A_4 \delta_D(x)^{-\alpha} \ell_1(\delta_D(x))  
	\quad \text{for all $ x\in D$ with $\delta_D(x)<1$.}
\end{align}   
\emph{When $\alpha\le 1$, we further assume that $a_0>0$.}
For any $q \in [\alpha-1, \alpha) \cap (0,\alpha)$, let
   \begin{align*}
   	\sK_\alpha(q):=\left\{ \kappa \in \sK_\alpha: \text{ the constant $a_0$ in \eqref{e:critical-class} is $\sC(d,\alpha,q)$}\right\}.
   \end{align*}
 By \eqref{e:criticall-killing-range}, we get that $\sK_\alpha(q)$, $q\in [\alpha-1,\alpha) \cap (0,\alpha)$, are disjoint and  $\sK_\alpha = \cup_{q \in [\alpha-1, \alpha) \cap (0,\alpha)} 	\sK_\alpha(q)$.

  Define    the bilinear form $(\overline \sE^{D}, \overline \sF^{ D})$ on  $L^2(D,dx)$ defined by
  \begin{align*}
  	\overline \sE^{D}(u,v) &:= \int_{D\times D\setminus \diag} (u(x)-u(y))(v(x)-v(y)) \frac{K(x,y)}{|x-y|^{d+\alpha}} dxdy, \nn\\
  	\overline	\sF^{D}&:= \bigg\{ u \in L^2(D,dx): \int_{D\times D\setminus \diag}  \frac{(u(x)-u(y))^2}{|x-y|^{d+\alpha}} dxdy<\infty\bigg\}.
  \end{align*}
 By \cite[Remark 2.1(1)]{BBC03}, 
 $(\overline \sE^D,\overline \sF^D)$ is a regular symmetric Dirichlet form on $L^2(D,dx)$.   Let   $\overline Y^D$ be the Hunt process associated with $(\overline \sE^D, \overline \sF^D)$ 
  and let $Y^D$ be the killed process of $\overline Y^D$ on $D$.
  The  process $\overline Y^D$ is referred to as an \textit{(actively) reflected 
  $\alpha$-stable-like process} on $\overline D$, and $Y^D$  a \textit{censored $\alpha$-stable-like process} on $D$.
 According to \cite[Theorem 1.1]{CK03}, $\overline Y^D$  satisfies   {\bf (A)} with $\phi(r)=r^\alpha$ and $T_0=1$.
  
    We fix a $q\in [\alpha-1,\alpha) \cap (0,\alpha)$ and a $\kappa \in \sK_\alpha(q)$, and let $X$ be a Hunt process on $D$ corresponding to the  Feynman-Kac semigroup
 \begin{align}\label{e:Feynman-Kac}
 	\E_x[f(X_t)] = \E_x \left[ e^{-\int_0^t \kappa(Y^D_s)ds} f(Y^D_t)\right], \quad t\ge 0, \, x\in D.
 \end{align}  
  Let $\sA$ be the  infinitesimal generator of $X$ and let
 \begin{align*}
 	C_c^2(D;\R^d)=\left\{ u:D\to \R \, : \, \text{there exists $f \in C_c^2(\R^d)$ such that $f=u$ on $D$}\right\}.
 \end{align*}
 Under \eqref{e:censored-1} and \eqref{e:censored-2},   by  \cite[Proposition 5.1]{CKSV25},  originally due to \cite[Corollary 4.5]{KSV23},	the  generator $\sA$  restricted to $C_c^2(D;\R^d)$ coincides with the following non-local operator  $\sL^{\kappa}$:
 \begin{align*}
 	\sL^\kappa u (x) = p.v. \int_D (u(y)-u(x)) \frac{K(x,y)}{|x-y|^{d+\alpha}}dy - \kappa(x) u(x).
 \end{align*} 
 
 The following  theorem is the main result of this subsection.

 \begin{thm}\label{t:critical}
 	Let $\alpha \in (0,2)$, $d\ge 2$ and $D\subset \R^d$ be a   $C^{1,{\rm Dini}}$ 
	open set with characteristics $(r_0,\ell_0)$.
	Assume \eqref{e:censored-1}  and \eqref{e:censored-2}.  Then for all 
 	$q\in [\alpha-1,\alpha) \cap (0,\alpha)$ and $\kappa \in \sK_\alpha(q)$,  	the Hunt  process $X$ corresponding to the 
 	Feynman-Kac semigroup \eqref{e:Feynman-Kac} 
 	has a  transition density $p(t,x,y)$ on $(0,\infty) \times D \times D$. Moreover, the following hold.

 	\smallskip
 	
 	\noindent 		(i) For every $T>0$, there are comparison constants such that for  all $(t,x,y) \in (0,T]\times D\times D$,
 	\begin{align*}
 		p(t,x,y) \asymp \bigg(1 \wedge\frac{\delta_D(x)}{t^{1/\alpha}} \bigg)^{q}\bigg(1 \wedge\frac{\delta_D(y)}{t^{1/\alpha}} \bigg)^{q} \bigg( t^{-d/\alpha} \wedge \frac{t}{|x-y|^{d+\alpha}}\bigg).
 	\end{align*}

 	\noindent 		(ii) Suppose in addition that $D$ is bounded. Then there exists $C\ge 1$ such that for all $(t,x,y) \in [3,\infty)\times D\times D$,
 	\begin{align*}
 		C^{-1}\delta_D(x)^{q}\delta_D(y)^{q}e^{-\lambda_1t}\le 		p(t,x,y) \le C\delta_D(x)^{q}\delta_D(y)^{q}e^{-\lambda_1t},
 	\end{align*}
 	where 
	$-\lambda_1:=\sup \text{\rm Re\,} \sigma(\sA)$  
	and $\sigma(\sA)$ is the spectrum of the generator $\sA$ of $X$. Moreover, $X$  
	admits 
	a Green function $G(x,y)$ on $D\times D$  such that for  all $(x,y) \in D\times D$,
 	\begin{align*}
 		G(x,y) \asymp\bigg(1 \wedge\frac{\delta_D(x)}{|x-y|} \bigg)^{q}\bigg(1 \wedge\frac{\delta_D(y)}{|x-y|} \bigg)^{q}  |x-y|^{\alpha-d} .
 	\end{align*}
 \end{thm}
 
 \smallskip

In the case of $\alpha$-stable-like processes, under conditions \eqref{e:pre-censored-1} and \eqref{e:pre-censored-2}, the function $\kappa_D$  defined in \eqref{e:killingfunction} belongs to $\sK_\alpha(\alpha/2)$ for all $\alpha\in (0, 2)$.
When $\alpha\in (1, 2)$, the constant 0 belongs to $\sK_\alpha(\alpha-1)$. So as  consequences of the theorem above, we get the following generalizations of the main results of \cite{CKS10a, CKS10}.  

 \begin{cor}\label{c:cricital}
		Let  $D\subset \R^d$, $d\ge 2$, be a   $C^{1,{\rm Dini}}$ open set.

				\noindent  (i) Assume \eqref{e:pre-censored-1}  and \eqref{e:pre-censored-2}. For any $\alpha \in (0,2)$, the assertions of Theorem \ref{t:critical}(i)-(ii) are valid for the killed stable-like process $Z^D$ with $q=\alpha/2$.
		
		\noindent	(ii) Assume \eqref{e:censored-1}  and \eqref{e:censored-2}. For $\alpha>1$, 
		the assertions of Theorem \ref{t:critical}(i)-(ii) are valid for the censored stable-like process $Y^D$ with $q=\alpha-1$.
	
 \end{cor}
 
 To prove Theorem \ref{t:critical}, we first establish several preliminary results. The proof of Theorem \ref{t:critical} will be given at the end of this subsection.

 In the  remainder of this subsection,  we always assume that we are in the setup of of Theorem \ref{t:critical}, that is, $\alpha \in (0,2)$, $d\ge 2$,  $D\subset \R^d$ is a   $C^{1,{\rm Dini}}$ open set 
 with characteristics $(r_0,\ell_0)$, and  \eqref{e:censored-1}  and \eqref{e:censored-2} hold.
 
 We begin by recalling the following result from \cite[Corollary 5.4]{CKSV25}. 
 \begin{prop}\label{p:maximum-censored}
 	Let $U\subset D$ be a
 	bounded 
 	open set. Suppose that    $u$ is a bounded Borel function on $D$   
	such that $u$ is continuous on $\overline{U}\cap D$ and vanishes continuously at  $\partial U\cap \partial D$,
	$u|_{U} \in C^2(U)$, and $\sL^\kappa u \ge 0$  on $U$ {\rm (}resp. $\sL^\kappa u \le 0$  on $U${\rm )}. Then we have $$u(x) \le  \E_x\left[ u(X_{\tau_U})\right] \;\; \text{{\rm (}resp. $u(x) \ge  \E_x\left[ u(X_{\tau_U})\right]${\rm )}} \quad \text{ for all $x\in U$.}$$
 \end{prop}

  By \cite[Theorem 2.7]{BBC03},  $Y^D$ is non-conservative for $\alpha>1$.
 Since the constant $a_0$ in \eqref{e:critical-class} is  strictly positive when $\alpha\le 1$, it follows that  $X$ is always non-conservative for all $\alpha \in (0,2)$.

 Recall that $D(r)=\{x\in D: \delta_D(x)<r\}$ for $r>0$.

 \begin{cor}\label{c:maximum-censored}
 	Let $r>0$. Suppose that    $u$ is a bounded Borel function on $D$  
	such that $u$ vanishes continuously on $\partial D$ and is continuous on $\overline{D(r)}$, 
	$u|_{D(r)} \in C^2(D(r))$, and $\sL^\kappa u \ge 0$  on $D(r)$ {\rm (}resp. $\sL^\kappa u \le 0$  on $D(r)${\rm )}. Then we have $$u(x) \le  \E_x\big[ u(X_{\tau_{D(r)}})\big] \;\;\; \text{{\rm \big(}resp. $u(x) \ge  \E_x\big[ u(X_{\tau_{D(r)}})\big]${\rm \big)}} \quad \text{ for all $x\in D(r)$.}$$
 \end{cor}
 \pf Since the proofs are similar, we only prove  the case where $\sL^\kappa u \ge 0$ on $D(r)$.
 
  Let $x\in D(r)$. 
  Define $D_n(r):=D(r) \cap B(x,2^n)$ for $n\ge 1$. By Proposition \ref{p:maximum-censored}, we have $u(x) \le \liminf_{n\to \infty}  \E_x[ u(X_{\tau_{D_n(r)}})] $. Thus, to obtain the result,  it suffices to show that
  \begin{align}\label{e:maximum-censored-claim}
  	\lim_{n\to \infty}  \E_x[ u(X_{\tau_{D_n(r)}})] = \E_x[ u(X_{\tau_{D(r)}})].
  \end{align}
  Let $a:=\sup_D |u|$.   For all $n \ge 1$, we have
  \begin{align*}
 &\left| \E_x[ u(X_{\tau_{D(r)}})]- \E_x[ u(X_{\tau_{D_n(r)}})]\right| = \left| \E_x\big[ u(X_{\tau_{D(r)}}) -  u(X_{\tau_{D_n(r)}}) : \tau_{D_n(r)} < \tau_{D(r)}\big] \right|\\
  &\le 2a \P_x ( \tau_{D_n(r)} < \tau_{D(r)} ) \le 2a \P_x ( \tau_{D_n(r)} < \tau_{D(r)} \wedge n )  + 2a \P_x (  \tau_{D(r)} \ge n ) .
  \end{align*}
  Since $X$ is non-conservative,  $\lim_{n\to \infty}  \P_x (  \tau_{D(r)} \ge n )  \le \lim_{n\to \infty}  \P_x (  \zeta \ge n ) =0$. Besides, we observe that for all $n\ge 1$,
   \begin{align*}
  	\P_x ( \tau_{D_n(r)} < \tau_{D(r)} \wedge n ) & =  \P_x ( \tau_{D_n(r)}  = \tau_{B(x,2^n)}< \tau_{D(r)} \wedge n )\le \P_x (  \tau_{B(x,2^n)}< n ) \le    \P_x (\tau^{\overline Y^D}_{B(x,2^n)}< n ),
  \end{align*}
  where $\tau^{\overline Y^D}_{B(x,2^n)}:=\inf\{t>0: \overline Y^D_t \notin B(x,2^n)\}$. By \cite{MU11},  $\overline Y^D$ is conservative.  Hence, using the Markov property and Lemma \ref{l:EP}, we get for all $n\ge 1$,
 \begin{align*}
 	\P_x (\tau^{\overline Y^D}_{B(x,2^n)}< n ) &\le  \P_x \left(\tau^{\overline Y^D}_{B(\overline Y^D_m,2^n/n)}< 1 \quad \text{for some $0\le m\le n-1$} \right)\\
 	&\le n \sup_{z\in \overline D}  \P_z \left( \tau^{\overline Y^D}_{B(z,2^n/n)}< 1\right) \le \frac{c_1n}{(2^n/n)^\alpha},
 \end{align*}
 implying that $\lim_{n\to \infty} 	\P_x (\tau^{\overline Y^D}_{B(x,2^n)}< n ) =0$. Hence, we obtain \eqref{e:maximum-censored-claim} and the proof is complete.  \qed

  By choosing a smaller $r_0$ and a larger   $A_4$ if necessary, 
  we assume without loss of generality that $ \ell_0(r) = \ell_0(r_0) \le 1/16$ for all $r\in [r_0,1]$, and $\ell_1(1) \le 1/16$, where $\ell_1\in \Dini$ is the function from \eqref{e:critical-class}. 
  Define $\ell_2(r):=\ell_0(r) \vee \ell_1(r)$ for $r\in (0,1]$. Clearly, $ \ell_2 \in \Dini$ and $\ell_2(1) \le 1/16$.   Let 
   $$\eps_0:=\frac{ \theta \wedge q \wedge(\alpha-q)  \wedge 1}{4} .$$ By Lemma \ref{l:Dini-regularizing}, there exists $\ell \in \Dini\cap C^2((0,\infty))$ such that $\ell_2(r) \le \ell(r)$, $r\ell'(r) \le  2\eps_0\ell(r)$, $|r^2\ell''(r)|\le 6\eps_0\ell(r)$ and   $\ell(r) \le \eps_0 \int_0^r s^{-1}\ell(s)ds$ for all $r\in (0,1]$,     $r^{-\eps_0}\ell(r)$ is non-increasing on $(0,1]$, and  
    $\ell(1) \le  4\ell_2(1) \le 1/4$.  We extend  $\ell$ to $(0,\infty)$ by setting $\ell(r)=\ell(1)$ for $r\ge1$.   Define
 \begin{align*}
 	F_\ell(r):= \int_0^r \frac{\ell(s)}{s}ds, \quad r>0.
 \end{align*}
Note that $F_\ell$ is slowly varying at $\infty$ and that for all $0<r\le R$,
 \begin{align}\label{e:F-ell-scaling}
 	F_\ell(R) =  \int_0^r \frac{\ell(R u/r)}{u}du \le \bigg(\frac{R}{r}\bigg)^{\eps_0}\int_0^r \frac{\ell( u)}{u}du = \bigg(\frac{R}{r}\bigg)^{\eps_0} F_\ell(r).
 \end{align}

  According to \cite[(5.4)]{BBC03}, it holds that\begin{align}\label{e:critical-harmonic}	p.v.	\int_{\R^d_+} \frac{y_d^q  - x_d^q }{|x-y|^{d+\alpha}}dy=\sC(d,\alpha,q) x_d^{q-\alpha} \quad \text{ for all $x=(\wt x,x_d) \in \R^d_+$} . \end{align}

 \begin{lemma}\label{l:critical-superharmonic} There exists $C>0$ such that for all $r\in (0,1]$ and $x=(\wt x,x_d) \in \R^d_+$,
 	\begin{align*}
 	p.v.	\int_{\R^d_+} \frac{y_d^q F_\ell(y_d/r) - x_d^q F_\ell(x_d/r)}{|x-y|^{d+\alpha}}dy \ge    \sC(d,\alpha,q)x_d^{q-\alpha}   F_\ell(x_d/r)  + Cx_d^{q-\alpha}\ell(x_d/r).
 	\end{align*}
 \end{lemma}
 \pf  Since $F_\ell\in C^2((0,\infty))$, $\lim_{R\to \infty} F_\ell(R)/\log R<\infty$ and $q<\alpha$, the principal value integral  is well-defined.
 Using \eqref{e:critical-harmonic}, we obtain
 \begin{align}\label{e:critical-superharmonic-1}
 	&p.v.	\int_{\R^d_+} \frac{y_d^q F_\ell(y_d/r) - x_d^q F_\ell(x_d/r)}{|x-y|^{d+\alpha}}dy \nn\\
 	&=p.v.	\int_{\R^d_+} \frac{( F_\ell(x_d/r) + q^{-1}\ell(x_d/r))(y_d^q-x_d^q) + y_d^q (F_\ell(y_d/r)-F_\ell(x_d/r)) - q^{-1}\ell(x_d/r)(y_d^q-x_d^q)}{|x-y|^{d+\alpha}}dy \nn\\	&= \sC(d,\alpha,q)x_d^{q-\alpha} (  F_\ell(x_d/r)  + q^{-1}\ell(x_d/r)) + p.v. 	\int_{ \R^d_+} \frac{H(y_d)}{|x-y|^{d+\alpha}}dy,
  \end{align}
 where $H(y_d):=y_d^q (F_\ell(y_d/r)-F_\ell(x_d/r)) - q^{-1}\ell(x_d/r)(y_d^q-x_d^q).$ Observe that for all $0<s<x_d$, 
\begin{align*}
	H'(s) = qs^{q-1} (F_\ell(s/r)-F_\ell(x_d/r)) + s^{q-1} \ell(s/r) - s^{q-1}\ell(x_d/r)  \le 0,\end{align*}
	and that for all $s>x_d$,
	\begin{align*}
		H'(s) &= qs^{q-1} (F_\ell(s/r)-F_\ell(x_d/r)) + s^{q-1} \ell(s/r) - s^{q-1}\ell(x_d/r)  \\
		&\ge qs^{q-1} (F_\ell(s/r)-F_\ell(x_d/r))  \ge \frac{qs^{q-1}(s-x_d)}{r}  \inf_{u \in (x_d,s)} \frac{\ell(u/r)}{u/r} \ge qs^{q-2} \ell(x_d/r)(s-x_d) \ge 0.
\end{align*}
Hence, we deduce that $H(s) \ge H(x_d)=0$ for all $s>0$ and 
\begin{align*}
	H(s) \ge \int_{2x_d}^{3x_d} H'(u)du \ge q (3x_d)^{q-2}   \ell(x_d/r) \int_{2x_d}^{3x_d} (s-x_d) ds = c_1 x_d^q \ell(x_d/r) \quad \text{for all $s\ge 3x_d$}.
\end{align*} 
It follows that
\begin{align*}
	 p.v. 	\int_{ \R^d_+} \frac{H(y_d)}{|x-y|^{d+\alpha}}dy &\ge c_1x_d^q \ell(x_d/r) \int_{y=(\wt y, y_d)\in  \R^d_+:y_d \ge 3x_d} \frac{ dy}{|x-y|^{d+\alpha}}\\
	 & =  c_1x_d^{q-\alpha} \ell(x_d/r) \int_{z=(\wt z, z_d)\in  \R^d_+:z_d \ge 3} \frac{ dz}{|(\wt 0,1)-z|^{d+\alpha}} = c_2 x_d^{q-\alpha} \ell(x_d/r),
\end{align*}
where we used the change of the variables $y_d = x_dz_d$ in the first equality. Combining this with \eqref{e:critical-superharmonic-1}, since $\sC(d,\alpha,q)\ge 0$, we get the result.
  \qed

 Let $\rho$ be the regularized distance for $D$  defined at the beginning of this section and  $r_1>0$ be the constant from Lemma \ref{l:inward-normal}. 
 Set 
 $\chi:=(d,r_0,\ell,\alpha, q,A_1,A_2,A_3,A_4,\theta)$, where $A_i$, 
 $1\le i\le 4$, and $\theta$ are the constants in \eqref{e:censored-1}, \eqref{e:censored-2} and \eqref{e:critical-class}.

 \begin{lemma}\label{l:critical-power-functions}
Let $r\in (0,r_1/4]$ and  $x_0 \in D(r)$, and 
define for $y\in \R^d$, 
 	\begin{align*} 
 	f(y) &:=\frac{ \left(\rho(x_0) + \nabla \rho(x_0) \cdot (y-x_0) \right)_+}{|\nabla \rho(x_0)| }, \\
 	f_q(y)&:=f(y)^q \quad \text{ and } \quad g_q (y)= f(y)^q F_\ell(f(y)/r).
 	\end{align*}
 	There exist positive constants 
	$C_1,C_2,C_3$ 
	depending only on $\chi$ such that
 	\begin{align}
 		|	\sL^{\kappa} f_q(x_0) |&\le 
		C_1\delta_D(x_0)^{q-\alpha} \ell(\delta_D(x_0)),
		\label{l:critical-power-functions(i)}\\
 			\sL^{\kappa} g_q(x_0) &\ge 
			\left( C_2   - C_{3} F_\ell(\delta_D(x_0)/r)  \right) 
			\delta_D(x_0)^{q-\alpha} \ell(\delta_D(x_0)/r) .\label{l:critical-power-functions(ii)}
 	\end{align}
 \end{lemma}
 \pf Set $s:=\delta_D(x_0)$.  Define 
 $E: = \{ y \in \R^d: \rho(x_0) + \nabla \rho(x_0)\cdot (y-x_0) > 0 \}.$
 By Lemma \ref{l:censored-domain-closedness}(i), 
 \begin{align}\label{e:censored-power-functions-0}
 c_1^{-1}s \le 	\delta_E(x_0) \le c_1 s \quad \text{and}  \quad |\delta_E(x_0)-s| \le c_2 s\ell(s).
 \end{align}
 
 (i) Observe that
 \begin{align*}
 	\sL^\kappa f_q(x_0)   &= p.v. \int_{D\cap E} (f_q(y)-f_q(x_0)) \frac{K(x_0,y)}{|x_0-y|^{d+\alpha}}dy  - \int_{D\setminus E}  \frac{f_q(x_0)K(x_0,y)}{|x_0-y|^{d+\alpha}}dy - \kappa(x_0) f_q(x_0) \\
 	&= \bigg( p.v. \int_{E} (f_q(y)-f_q(x_0)) \frac{K(x_0,x_0)}{|x_0-y|^{d+\alpha}}dy - \kappa(x_0) f_q(x_0) \bigg) \\
 	&\quad\;\;   +   \int_{D\cap E} (f_q(y)-f_q(x_0)) \frac{(K(x_0,y)-K(x_0,x_0))}{|x_0-y|^{d+\alpha}}dy \\
 		&\quad \;\;-  \int_{E\setminus D} (f_q(y)-f_q(x_0)) \frac{K(x_0,x_0)}{|x_0-y|^{d+\alpha}}dy   - \int_{D\setminus E}  \frac{f_q(x_0)K(x_0,y)}{|x_0-y|^{d+\alpha}}dy \\
 	&=:I_1+I_2-I_3-I_4.
 \end{align*}

 For $I_1$, by \eqref{e:critical-harmonic}, we have 
 \begin{align*}
 	&|I_1|  =\delta_E(x_0)^q \left|\sC(d,\alpha,q)K(x_0,x_0) \delta_E(x_0)^{-\alpha} - \kappa(x_0) \right|  \\
 	&\le \delta_E(x_0)^q \Big[\left|\sC(d,\alpha,q)K(x_0,x_0) \delta_D(x_0)^{-\alpha} - \kappa(x_0)  \right| + \sC(d,\alpha,q)K(x_0,x_0)\left|  \delta_E(x_0)^{-\alpha} -  \delta_D(x_0)^{-\alpha}  \right| \Big].
 \end{align*}
 Thus, using  \eqref{e:critical-class} (with $a_0=\sC(d,\alpha,q)$), \eqref{e:censored-1} and \eqref{e:censored-power-functions-0}, we obtain
 \begin{align*}
 	|I_1| \le c_3s^{q-\alpha} \ell(s)  + c_4 s^{q-\alpha-1} | \delta_E(x_0) -  \delta_D(x_0)|  \le c_5 s^{q-\alpha}\ell(s).
 \end{align*}

It follows from \eqref{e:censored-power-functions-0} that for any $y\in \R^d$,
\begin{align}\label{e:censored-power-functions-1}
	f_q(y) \le (\delta_E(x_0)+ |x_0-y|)^q \le c_6(s^q  +  |x_0-y|^q).
\end{align} Hence, using \eqref{e:censored-1} and   Lemma \ref{l:censored-domain-closedness}, since $\ell$ satisfies \eqref{e:ell-scaling} with $\eps=\eps_0< (\alpha-q)/2$, we get that
 \begin{align*}
 	|I_3| + |I_4| \le c_7\int_{(D\setminus E) \cup (E\setminus D)} \bigg(  \frac{s^q }{|x_0-y|^{d+\alpha}} + \frac{1}{|x_0-y|^{d+\alpha-q}}\bigg)dy \le c_8 s^{q-\alpha}\ell(s).
 \end{align*}

  For $I_2$, by \eqref{e:censored-1} and \eqref{e:censored-2}, we have
 \begin{align*}
|I_2| &\le c_{9} \bigg( \int_{B(x_0, s/(2c_1))} + \int_{(D\cap E) \setminus B(x_0, s/(2c_1))} \bigg)   | f_q(y)-f_q(x_0)|\frac{ (|x_0-y| \wedge 1)^{\theta}}{|x_0-y|^{d+\alpha}}dy=:I_{2,1} + I_{2,2}.
 \end{align*}
 For all $y \in B(x_0,s/(2c_1))$, since $ f(y) \asymp f(x_0) \asymp  s$, by the mean value theorem, 
 \begin{align}\label{e:censored-power-functions-2}
 	|f_q(y) - f_q(x_0)| \le c_{10}s^{q-1} |y-x_0|.
 \end{align}
Set $\theta_0:=\theta \wedge ((\alpha-q)/2)$. By \eqref{e:censored-power-functions-1} and  \eqref{e:censored-power-functions-2}, 
 \begin{align*}
 	I_{2,1}+I_{2,2} &\le c_{11} s^{q-1}\int_{B(x_0,s/(2c_1))} \frac{dy}{|x_0-y|^{d+\alpha-\theta-1}}  + c_{11} \int_{B(x_0,s/(2c_1))^c} \frac{s^q\, dy}{|x_0-y|^{d+\alpha-\theta_0}}\\
 	&\quad  + c_{11} \int_{B(x_0,s/(2c_1))^c} \frac{ dy}{|x_0-y|^{d+\alpha-\theta_0-q}}  \\
 	&\le c_{12} (s^{q+\theta-\alpha} +  s^{q+\theta_0-\alpha} ) \le 2c_{12} s^{q+\theta_0-\alpha}.
 \end{align*}
 Since $
 \ell(s) \ge s^{\eps_0} \ell(1) \ge  s^{\theta_0}\ell(1)$, the proof of \eqref{l:critical-power-functions(i)} is complete. 
 
(ii) To establish \eqref{l:critical-power-functions(ii)}, we follow the arguments  for \eqref{l:critical-power-functions(i)}. Observe that
 \begin{align*}
 	\sL^\kappa g_q(x_0)    	&= \bigg( p.v. \int_{E} (g_q(y)-g_q(x_0)) \frac{K(x_0,x_0)}{|x_0-y|^{d+\alpha}}dy - \kappa(x_0) g_q(x_0) \bigg) \\
 	&\quad \;\;  +   \int_{D\cap E} (g_q(y)-g_q(x_0)) \frac{(K(x_0,y)-K(x_0,x_0))}{|x_0-y|^{d+\alpha}}dy \\
 	&\quad \;\;-  \int_{E\setminus D} (g_q(y)-g_q(x_0)) \frac{K(x_0,x_0)}{|x_0-y|^{d+\alpha}}dy   - \int_{D\setminus E}  \frac{g_q(x_0)K(x_0,y)}{|x_0-y|^{d+\alpha}}dy \\
 	&=:I_1'+I_2'-I_3'-I_4'.
 \end{align*}
 Using Lemma \ref{l:critical-superharmonic}, \eqref{e:censored-1},  \eqref{e:critical-class},  \eqref{e:censored-power-functions-0}, and the scaling properties of $\ell$ and $F_\ell$, we get
 \begin{align*}
 	&I_1'  =  p.v. \int_{E} (\delta_E(y)^qF_\ell(\delta_E(y)/r)-\delta_E(x_0)^qF_\ell(\delta_E(x_0)/r)) \frac{K(x_0,x_0)}{|x_0-y|^{d+\alpha}}dy - \kappa(x_0) \delta_E(x_0)^q F_\ell(\delta_E(x_0)/r)  \\
 	& \ge  c_{13} K(x_0,x_0) \delta_E(x_0)^{q-\alpha} \ell(\delta_E(x_0)/r)  + \delta_E(x_0)^q F_\ell(\delta_E(x_0)/r)	\left(\sC(d,\alpha,q)K(x_0,x_0) \delta_E(x_0)^{-\alpha} - \kappa(x_0) \right) \\ 
 		& \ge  c_{13} K(x_0,x_0) \delta_E(x_0)^{q-\alpha} \ell(\delta_E(x_0)/r)  - \delta_E(x_0)^q F_\ell(\delta_E(x_0)/r)	\left|\sC(d,\alpha,q)K(x_0,x_0) \delta_D(x_0)^{-\alpha} - \kappa(x_0) \right| \\ &\quad  - \sC(d,\alpha,q)K(x_0,x_0) \delta_E(x_0)^q F_\ell(\delta_E(x_0)/r)\left| \delta_D(x_0)^{-\alpha}- \delta_E(x_0)^{-\alpha}  \right|  \\
 	&\ge c_{14} s^{q-\alpha} \ell(s/r)  - c_{15} s^{q-\alpha} F_\ell(s/r) \ell(s).
 \end{align*}
 By \eqref{e:F-ell-scaling}, $F_\ell(a+b) \le (1+ b/a)^{\eps_0} F_\ell(a) $ for all $a,b>0$. Since $c_1^{-1}s\le f(x_0)\le c_1s$ and $F_\ell(f(x_0)/r)\le c F_\ell(s/r)$, we have  for all $y \in \R^d$,
 \begin{align}\label{e:censored-power-functions-3}
 g_q(y)&\le (f(x_0) + |x_0-y|)^q F_\ell((f(x_0) + |x_0-y|)/r)   \le f(x_0)^q(1+ |x_0-y|/f(x_0))^{q+\eps_0}F_\ell(f(x_0)/r)\nn\\
 &\le c_{16}s^q ( 1 + |x_0-y|^{q+\eps_0}/s^{q+\eps_0}) F_\ell(s/r). 
 \end{align} 
 Thus,  using $g_q(x_0)\le c_{17} s^q F_\ell(s/r)$, 
 \eqref{e:censored-1} and applying Lemma \ref{l:censored-domain-closedness} (since $\eps_0<(\alpha-q)/4$), we get
 \begin{align*}
 	|I_3'| + |I_4'| &\le 
	c_{18} 
	s^qF_\ell(s/r)\int_{(D\setminus E) \cup (E\setminus D)} \frac{dy}{|x_0-y|^{d+\alpha}} +  
	c_{18} 
	s^{-\eps_0} F_\ell(s/r)\int_{(D\setminus E) \cup (E\setminus D)} \frac{dy}{|x_0-y|^{d+\alpha-q-\eps_0}}   \\
	&\le c_{19} 
	s^{q-\alpha}F_\ell(s/r)\ell(s).
 \end{align*}
Since $F_\ell(a) \ge\ell(a)$ for any $a>0$, we have
$$( s^q F_\ell(s/r))'=  qs^{q-1} F_\ell(s/r) + s^{q-1}\ell(s/r) \le (q+1)s^{q-1}F_\ell(s/r) \quad \text{for all $s>0$}.$$
 Thus,  for all $y \in B(x_0,s/(2c_1))$, by the mean value theorem and the scaling of $F_\ell$,
 \begin{align*}
 	|g_q(y) - g_q(x_0)| \le 
	c_{20}s^{q-1} F_\ell(s/r) |y-x_0|.
 \end{align*}
 Using this,  \eqref{e:censored-power-functions-3}, and the fact  $\theta_0+q+\eps_0< q + 3(\alpha-q)/4 <\alpha$,  we get from \eqref{e:censored-1} and \eqref{e:censored-2} that 
 \begin{align*}
 	|I_2'| &\le c_{21} \bigg( \int_{B(x_0, s/(2c_1))} + \int_{(D\cap E) \setminus B(x_0, s/(2c_1))} \bigg)   | g_q(y)-g_q(x_0)|\frac{ (|x_0-y| \wedge 1)^{\theta}}{|x_0-y|^{d+\alpha}}dy\\
 	&\le c_{22} s^{q-1}F_\ell(s/r)  \int_{B(x_0, s/(2c_1))}  \frac{dy}{|x_0-y|^{d+\alpha-\theta-1}}  +  c_{22} s^q F_\ell(s/r) \int_{B(x_0,s/(2c_1))^c} \frac{dy}{|x_0-y|^{d+\alpha-\theta_0}} \\
 	&\quad + c_{21} s^{-\eps_0} F_\ell(s/r) \int_{B(x_0,s/(2c_1))^c} \frac{dy}{|x_0-y|^{d+\alpha-\theta_0-q-\eps_0}}\\
 	&\le c_{23} s^{q+\theta_0-\alpha} F_\ell(s/r) \le c_{24} s^{q-\alpha} F_\ell(s/r)\ell(s).
 \end{align*}
Since $\ell$ is non-decreasing and $r\le 1$,  we have  $\ell(s)\le \ell(s/r)$,  the proof is complete. 	\qed
 
 For $r\in (0,r_1]$ and $\sigma \in (0,1/4]$, define $h_{r,\sigma},\psi_{r,\sigma}:\R^d  \to [0,\infty)$ by
 		\begin{align}\label{e:def-h-r-sigma-critical}
 		h_{r,\sigma}(x) := \begin{cases}
 			\displaystyle	\bigg(\frac{\rho(x)}{r}\bigg)^{q}    &\mbox{ if $x\in D(2\sigma r)$},\\[3pt]
 			\sigma^q  &\mbox{ if $x\in D\setminus D(2\sigma r)$},\\
 			0 &\mbox{ if $x\in D^c$},
 		\end{cases} 
 	\end{align} 
 		\begin{align}\label{e:def-psi-r-sigma-critical}
 		\psi_{r,\sigma}(x) := \begin{cases}
 			\displaystyle	\bigg(\frac{\rho(x)}{r}\bigg)^{q}   F_\ell\left(\frac{\rho(x)}{|\nabla \rho(x)|\,r}\right)    &\mbox{ if $x\in D(2\sigma r)$},\\[3pt]
 			\sigma^q F_\ell(\sigma)  &\mbox{ if $x\in D\setminus D(2\sigma r)$},\\
 			0 &\mbox{ if $x\in D^c$}.
 		\end{cases} 
 	\end{align}

 \begin{prop}\label{p:barriers-critical}
	There 
	exist positive constants 
	$C_1,C_2,C_3,C_4$  
	depending only on  $\chi$ such that for all $x\in D(\sigma r)$,
 	 \begin{align}
 	 	|\sL^\kappa h_{r,\sigma} (x)|&\le 
		C_1 r^{-q} \delta_D(x)^{q-\alpha} \ell(\delta_D(x)) + C_2
		\sigma^{q-\alpha} r^{-\alpha}, \label{e:barriers-critical-result-1}\\
 	 		\sL^\kappa  \psi_{r,\sigma}(x)& \ge\left( 
			C_3  - \frac{C_4
			F_\ell(\delta_D(x)/r)}{\ell(\sigma)}  \right)   r^{-q} \delta_D(x)^{q-\alpha} \ell(\delta_D(x)/r).\label{e:barriers-critical-result-2}
 	 \end{align}
 \end{prop}
 \pf 
   Pick an arbitrary  $x_0 \in D(\sigma r)$. Let   $s:=\delta_D(x_0)$ and $a:=|\nabla \rho(x_0)|^{-1}$. Note that, by \eqref{e:inward-normal-result-1},  $a\in [1/2,4]$.    
   Define for $y\in \R^d$, 
 	\begin{align*} 
 	f(y) &:=a(\rho(x_0) + \nabla \rho(x_0) \cdot (y-x_0) )_+,   \\
	 H_1 (y)&:= (f(y)/r)^q \qquad \text{and} \qquad H_2(y)=(f(y)/r)^q F_\ell(f(y)/r) .
 \end{align*}
We have $f(x_0)=a\rho(x_0)$ and $\nabla f(x_0) = a \nabla \rho(x_0)$.   By \eqref{e:regular-distance-1} and \eqref{e:inward-normal-result-1}, there exists $b\in (0,1/2]$ such that $f(y)/(a\rho(x_0)) \in [1/2,2]$ for all $y \in B(x_0, bs)$, and   
 \begin{align}\label{e:barriers-critical-1}
 |\rho(y)|  + 	|f(y)| \le  c( s+ |y-x_0|) \quad \text{ for all $y \in \R^d$.}
 \end{align}
 Further, using \eqref{e:regular-distance-2}, we get for any $y \in \R^d$,
 \begin{align}\label{e:barriers-critical-2}
 	\left|a\rho(y)_+-f(y)\right| \le 4| \rho(y) - \rho(x_0) - \nabla \rho(x_0) \cdot (y-x_0) | \le c  |y-x_0|\ell(|y-x_0|).
 \end{align}

 (i) 
 We first prove \eqref{e:barriers-critical-result-2}. Since $a^q \psi_{r,\sigma}(x_0) = H_2(x_0)$, we have
 \begin{align*}
 	& 	a^q\sL^\kappa \psi_{r,\sigma} (x_0) -  	\sL^\kappa H_2 (x_0) = p.v.\int_{D} \left(a^q\psi_{r,\sigma} (y) - H_2(y)  \right) \frac{K(x_0,y)}{|x_0-y|^{d+\alpha}}dy \nn\\
 	&=   \int_{B(x_0,bs)} \left( a^q \psi_{r,\sigma} (y) - H_2(y)  \right) \frac{(K(x_0,y)- K(x_0,x_0))}{|x_0-y|^{d+\alpha}}dy   \nn\\ 
 	&\quad + \int_{B(x_0,bs)} \left(a^q \psi_{r,\sigma} (y)  - H_2(y)  - (y-x_0) \cdot \nabla \left(a^q \psi_{r,\sigma}     - H_2\right)(x_0) \right) \frac{K(x_0,x_0)}{|x_0-y|^{d+\alpha}}dy    \nn\\
 	&\quad +   \int_{D \cap (B(x_0,\sigma r)\setminus B(x_0,bs))} \left( a^q \psi_{r,\sigma} (y) - H_2(y)  \right) \frac{K(x_0,y)}{|x_0-y|^{d+\alpha}}dy  \nn\\
 	&\quad +   \int_{D\setminus B(x_0,\sigma r)} \left( a^q \psi_{r,\sigma} (y) - H_2(y)  \right) \frac{K(x_0,y)}{|x_0-y|^{d+\alpha}}dy  \nn\\
 	&=:I_1+I_2+I_3+I_4.
 \end{align*}

 For any $y \in B(x_0,bs)$, we have $\rho(y) \asymp f(y) \asymp s$,
 \begin{align*}
 	\partial_{i} \left(  a^q \psi_{r,\sigma}   - H_2 \right)(y) &= \frac{qa}{r^{q} }    \left(a^{q-1}\rho(y)^{q-1} F_\ell(a\rho(y)/r)\partial_{i} \rho(y)  -f(y )^{q-1} F_\ell(f(y)/r) \partial_i \rho(x_0)  \right) \\
 	&\quad + \frac{a}{r^{q} }    \left(a^{q-1}\rho(y)^{q-1} \ell(a\rho(y)/r)\partial_{i} \rho(y)  -f(y )^{q-1} \ell(f(y)/r) \partial_i \rho(x_0)  \right) 
 \end{align*}
 and 
 \begin{align*}
 &	\partial_{ij}\left(  a^q \psi_{r,\sigma}   - H_2 \right)(y)= \frac{qa^q}{r^{q} }   \rho(y)^{q-1} F_\ell(a\rho(y)/r)\partial_{ij} \rho(y)  + \frac{a}{r^q}\rho(y)^{q-1} \ell(a\rho(y)/r)\partial_{ij} \rho(y) \\
 	&\qquad\qquad    +  \frac{ q(q-1)a^2}{r^{q} }    \left(a^{q-2}\rho(y)^{q-2} F_\ell(a\rho(y)/r)\partial_{i} \rho(y) \partial_{j} \rho(y)  -f(y )^{q-2} F_\ell(f(y)/r) \partial_{i} \rho(x_0) \partial_{j} \rho(x_0)  \right) \\
 	&\qquad\qquad   + \frac{(2q-1)a^2}{r^{q} }    \left(a^{q-2}\rho(y)^{q-2} \ell(a\rho(y)/r)\partial_{i} \rho(y) \partial_{j} \rho(y)  -f(y )^{q-2} \ell(f(y)/r) \partial_i \rho(x_0) \partial_j \rho(x_0)  \right) \\
 		&\qquad\qquad   + \frac{a^2}{r^{q+1} }    \left(a^{q-1}\rho(y)^{q-1} \ell'(a\rho(y)/r)\partial_{i} \rho(y)  \partial_{j} \rho(y)  -f(y )^{q-1} \ell'(f(y)/r) \partial_i \rho(x_0) \partial_j \rho(x_0)  \right) .
 \end{align*}
 By \eqref{e:regular-distance-2} and \eqref{e:regular-distance-3}, for any $y \in B(x_0,bs)$, we have 
 $| \nabla \rho(y) - \nabla \rho(x_0)| \le  c_1\ell(s)$ and  $| D^2 \rho(y)| \le  c_2s^{-1}\ell(s)$. 
 Thus, using \eqref{e:barriers-critical-2}, the facts that 
 $u\ell'(u) \le  c_3\ell(u)\le c_4F_\ell(u)$ and  $|u^2\ell''(u)|\le c_5\ell(u)$ for $u>0$, 
 and the scaling properties of $F_\ell$ and $\ell$, we get that  for any $y \in B(x_0,bs)$,
 \begin{align*}
 	\left|D^k ( a^q \psi_{r,\sigma}   - H_2 ) (y)\right|
	\le  \frac{c_6 s^{q-k}F_\ell(s/r)\ell(s) }{r^q} \quad \text{for $k\in \{1,2\}$}.
 \end{align*}
 It follows that for any $y \in B(x_0,bs)$,
 \begin{align}\label{e:barriers-critical-3}
 	\left| a^q \psi_{r,\sigma}(y)   - H_2(y) \right| \le \frac{c_7  s^{q-1}F_\ell(s/r)\ell(s) |y-x_0| }{r^q}  
 \end{align}
 and
 \begin{align}\label{e:barriers-critical-4}
 	\left| a^q \psi_{r,\sigma}(y)  - H_2(y)  - (y-x_0) \cdot \nabla \left(a^q \psi_{r,\sigma}    - H_2\right)(x_0)  \right| \le \frac{c_8  s^{q-2} F_\ell(s/r)\ell(s) |y-x_0|^2 }{r^{q}} .
 \end{align}
 By \eqref{e:censored-2} and  \eqref{e:barriers-critical-3}, we have
 \begin{align*}
 	|I_1| \le \frac{c_9 s^{q-1} F_\ell(s/r)\ell(s)}{r^{q}} \int_{B(x_0,bs)} \frac{dy}{|x_0-y|^{d+\alpha-\theta-1}} = \frac{c_{10} s^{q+\theta-\alpha}F_\ell(s/r)\ell(s)}{r^{q}}.
 \end{align*}
 Further, by  \eqref{e:censored-1} and \eqref{e:barriers-critical-4}, it holds that 
 \begin{align*}
 	|I_2| \le \frac{c_{11}  s^{q-2}F_\ell(s/r)\ell(s)  }{r^{q}}\int_{B(x_0,bs)} \frac{dy}{|x_0-y|^{d+\alpha-2}}  = \frac{c_{12}s^{q-\alpha}F_\ell(s/r)\ell(s)}{r^{q}}.
 \end{align*}

 Observe that, by the scaling of $F_\ell$,  $\psi_{r,\sigma} \le c_{13} \sigma^q F_\ell(\sigma)$ on $D\setminus B(x_0,\sigma r)$. Using this,
  \eqref{e:censored-1}, \eqref{e:barriers-critical-1} and the scaling of $F_\ell$, since   $s\le \sigma r$ and  $u^{q-\alpha}F_\ell(u/r)\ell(u/r)$ is decreasing   by  \eqref{e:F-ell-scaling}, we obtain
 \begin{align*}
 	|I_4|&\le c_{14} \int_{B(x_0,\sigma r)^c} \frac{( \sigma^q F_\ell(\sigma) +  (|x_0-y|/r)^q F_\ell(|x_0-y|/r))}{|x_0-y|^{d+\alpha}} dy  \\
 	& \le \frac{c_{15}\sigma^q F_\ell(\sigma)}{(\sigma r)^\alpha}  +  \frac{c_{16}  F_\ell(\sigma)}{\sigma^{(\alpha-q)/2} r^{(\alpha+q)/2}} \int_{B(x_0,\sigma r)^c} \frac{dy}{|x_0-y|^{d+(\alpha-q)/2}} \\
 	&=  \frac{c_{17}\sigma^{q-\alpha} F_\ell(\sigma)}{r^\alpha} \le \frac{c_{17}F_\ell(s/r)\ell(s/r)}{\ell(\sigma)r^\alpha } \bigg(\frac{s}{r}\bigg)^{q-\alpha}.
 \end{align*}
 
 By \eqref{e:F-ell-scaling}, since $\eps_0<q$,  $u^q F_\ell(v) \ge v^q F_\ell(u)$ for any $0<v\le u$. Thus, for all $y \in D\cap (B(x_0,\sigma r)\setminus  B(x_0,bs))$,  using \eqref{e:barriers-critical-1}, $s+|x_0-y|\le (1+b^{-1})|x_0-y|$ and the scaling of $F_\ell$, we see that   
 \begin{align}\label{e:critical-I3-diff}
 		\left|  a^q \psi_{r,\sigma}(y) - H_2(y)  \right| &\le \frac{| a^q \rho(y)^q  - f(y)^q|}{r^q}\bigg|  F_\ell\left(\frac{a\rho(y)}{r}\right) +  F_\ell\left(\frac{f(y)}{r}\right) \bigg| \nn\\
 		&\le \frac{c_{18} F_\ell( |x_0-y|/r)| a^q \rho(y)^q  - f(y)^q|}{r^q}.
 \end{align}
 
 For $I_3$, we consider the cases $q\ge 1$ and $q<1$ separately.
We first assume   $q\ge 1$. For any $y \in D\cap (B(x_0,\sigma r)\setminus  B(x_0,bs))$,  using \eqref{e:critical-I3-diff}, \eqref{e:barriers-critical-1} and \eqref{e:barriers-critical-2},  we obtain
 \begin{align*}
	\left|  a^q \psi_{r,\sigma}(y) - H_2(y) \right| 	 &\le  \frac{c_{19}F_\ell(|x_0-y|/r  )((a^q \rho(y)) \vee  f(y))^{q-1} | a\rho(y) - f(y)| }{r^{q}}  \nn\\
 	 & \le \frac{c_{20} |x_0-y|^{q} F_\ell ( |x_0-y|/r ) \ell(|x_0-y|)}{r^{q}}.
 \end{align*}
 Using this and \eqref{e:censored-1}, since $u^{-(\alpha-q)/2}F_\ell(u/r)\ell(u)$ is non-increasing, we obtain
 \begin{align}\label{e:barriers-I3-case1}
 	|I_3|& \le \frac{c_{21}}{r^q} \int_{B(x_0,\sigma r)\setminus B(x_0,bs)} \frac{F_\ell(|x_0-y|/r)\ell(|x_0-y|)}{|x_0-y|^{d+\alpha-q}} dy \nn\\
 	&\le \frac{c_{21}F_\ell(bs/r)\ell(bs)}{r^q (bs)^{(\alpha-q)/2}} \int_{B(x_0,\sigma r)\setminus B(x_0,bs)} \frac{dy}{|x_0-y|^{d+(\alpha-q)/2}} \le \frac{c_{22} s^{q-\alpha}F_\ell(s/r)\ell(s/r)}{r^q }.
 \end{align}
 
 Suppose $q<1$. Using \eqref{e:critical-I3-diff}, the inequality $|u^q-v^q|(u^{1-q} + v^{1-q}) \le 2|u-v|$ for all $u,v\ge 0$, \eqref{e:barriers-critical-2} and \eqref{e:regular-distance-1}, we obtain  for any $y \in D\cap (B(x_0,\sigma r)\setminus  B(x_0,bs))$,
 \begin{align*}
 	\left|  a^q \psi_{r,\sigma}(y) - H_2(y)  \right| &\le  \frac{2c_{18} F_\ell( |x_0-y|/r)| a^q \rho(y)  - f(y)|}{r^q( (a\rho(y))^{1-q} + f(y)^{1-q})} \\
 	&\le  \frac{c_{23}|x_0-y| F_\ell(|x_0-y|/r)\ell(|x_0-y|)}{\rho(y)^{1-q}r^{q} }   \le  \frac{c_{24}|x_0-y| F_\ell(|x_0-y|/r)\ell(|x_0-y|)}{\delta_D(y)^{1-q}r^{q} } .
 \end{align*}
Hence, using the coordinate system CS$_{Q_{x_0}}$, and the scaling of $F_\ell$ and  $\ell$, we get that
 \begin{align}\label{e:barriers-I3-case1+}
 	|I_3|& \le \frac{c_{25}}{r^q} \int_{D\cap (B(x_0,\sigma r)\setminus B(x_0,bs))} \frac{F_\ell(|x_0-y|/r)\ell(|x_0-y|)}{\delta_D(y)^{1-q}|x_0-y|^{d+\alpha-1}} dy \nn\\
 	& \le \frac{c_{26}}{r^q}  \int_0^{r} \int_0^{r}  \frac{\ind_{\{l^2+|s-y_d|^2>(bs)^2\}} F_\ell((l^2+|s-y_d|^2)^{1/2}/r ) \ell((l^2+|s-y_d|^2)^{1/2} )}{y_d^{1-q} (l^2+|s-y_d|^2)^{(d+\alpha-1)/2}} \, l^{d-2 } dldy_d \nn\\
 	& \le \frac{c_{26}}{r^q}  \int_0^{r}  \frac{\ind_{\{|s-y_d|^2>(bs)^2/2\}} F_\ell(2^{1/2} |s-y_d| /r) \ell(2^{1/2} |s-y_d| )}{y_d^{1-q} |s-y_d|^{d+\alpha-1}} \int_0^{|s-y_d|}  l^{d-2 } dldy_d  \nn\\
 	&\quad +\frac{c_{26}}{r^q}\int_0^{r} \frac{\ind_{\{l^2>(bs)^2/2\}}F_\ell(2^{1/2}l/r)\ell(2^{1/2}l )}{ l^{\alpha+1 }}    \,\int_{0}^{r}   \frac{\ind_{\{|s-y_d|<l\}}}{y_d^{1-q}} dy_d dl \nn\\
 	& \le \frac{c_{27}}{r^q}  \int_0^{(1-2^{-1/2}b)s}  \frac{F_\ell(|s-y_d|/r) \ell( |s-y_d| )}{y_d^{1-q} |s-y_d|^{\alpha}} dy_d   + \frac{c_{21}}{r^q}  \int_{(1+2^{-1/2}b)s} ^{r}  \frac{ F_\ell(|s-y_d|/r)\ell( |s-y_d| )}{y_d^{1-q} |s-y_d|^{\alpha}} dy_d  \nn\\
 	&\quad +\frac{c_{27}}{r^q}\int_{2^{-1/2}bs}^{r} \frac{F_\ell(l/r)\ell(l )}{ l^{\alpha+1 }}    \,\int_{0}^{l+s}   \frac{dy_d}{y_d^{1-q}}  dl\nn\\
 	& =:\frac{c_{27}}{r^q} \left( I_{3,1}+I_{3,2}+I_{3,3}\right) .
 \end{align} 
 Since $|s-y_d|\asymp s$ for $y_d \in (0,(1-2^{-1/2}b)s)$,  we have
 \begin{align*}
 	I_{3,1} \le c_{28}s^{-\alpha}F_\ell(s/r)\ell(s)\int_0^{(1-2^{-1/2}b)s}  \frac{dy_d}{y_d^{1-q}}  = c_{29}s^{q-\alpha}F_\ell(s/r)\ell(s).
 \end{align*}
 Further, since $|s-y_d| \asymp y_d$ for $y>(1+2^{-1/2}b)s$  $u^{-(\alpha-q)/2}F(u/r)\ell(u)$ is non-increasing, we get
 \begin{align}\label{e:barriers-I32}
 	I_{3,2}& \le c_{30} \int_{(1+2^{-1/2}b)s} ^{r}  \frac{ F_\ell(y_d/r)\ell( y_d )}{y_d^{1+\alpha-q}} dy_d \nn\\
 	&\le \frac{c_{30}F_\ell(s)\ell(s )}{s^{(\alpha-q)/2}} \int_{(1+2^{-1/2}b)s} ^{r}  \frac{ dy_d}{y_d^{1+(\alpha-q)/2}}  \le c_{31} s^{q-\alpha}F_\ell(s/r)\ell(s).
 \end{align} 
 For $I_{3,3}$, using \eqref{e:barriers-I32}, we obtain
 \begin{align}\label{e:barriers-I33}
 	I_{3,3} & = \frac{1}{q} \left( \int_{2^{-1/2}bs}^{(1+2^{-1/2}b)s} \frac{(l+s)^qF_\ell(l/r)\ell(l )}{ l^{\alpha+1 }}     dl +  \int_{(1+2^{-1/2}b)s}^{r} \frac{(l+s)^qF_\ell(l/r)\ell(l )}{ l^{\alpha+1 }}     dl \right) \nn\\
 	&\le c_{32} s^{q-\alpha-1}F_\ell(s/r) \ell(s) \int_{2^{-1/2}bs}^{(1+2^{-1/2}b)s}   dl  + c_{32} \int_{(1+2^{-1/2}b)s}^{r} \frac{F_\ell(l/r)\ell(l)}{l^{1+\alpha-q}}dl \nn\\
 	&\le c_{33} s^{q-\alpha}F_\ell(s/r)\ell(s).
 \end{align}

 Since $\ell(s)\le \ell(s/r)$, $\ell(\sigma)\le \ell(1)\le 1/4$ and  $\sL^\kappa H_2(x_0)\ge  \left( c_{34}   - c_{35} F_\ell(s/r)  \right)r^{-q} s^{q-\alpha} \ell(s/r)  $ by Lemma \ref{l:critical-power-functions}(ii), 
combining the estimates for $|I_1|$, $|I_2|$, $|I_3|$ and $|I_4|$, we conclude that \eqref{e:barriers-critical-result-2} holds.

 \smallskip

 (ii) For \eqref{e:barriers-critical-result-1},  we  follow the arguments in (i) line by line. Since $a^qh_{r,\sigma}(x_0) = H_1(x_0)$, we have
 \begin{align*}
 	& 	a^q\sL^\kappa h_{r,\sigma} (x_0) -  	\sL^\kappa H_1 (x_0) = p.v.\int_{D} \left(a^qh_{r,\sigma}(y) - H_1(y) \right) \frac{K(x_0,y)}{|x_0-y|^{d+\alpha}}dy \nn\\
 	&=   \int_{B(x_0,bs)}  \left(a^qh_{r,\sigma}(y) - H_1(y) \right) \frac{(K(x_0,y)- K(x_0,x_0))}{|x_0-y|^{d+\alpha}}dy   \nn\\ 
 	&\quad + \int_{B(x_0,bs)} \left(a^q h_{r,\sigma}(y)  - H_1(y) - (y-x_0) \cdot \nabla \left(a^q h_{r,\sigma}    - H_1\right)(x_0) \right) \frac{K(x_0,x_0)}{|x_0-y|^{d+\alpha}}dy    \nn\\
 	&\quad +   \int_{D \cap (B(x_0,\sigma r)\setminus B(x_0,bs))} \left( a^q h_{r,\sigma}(y) - H_1(y)  \right) \frac{K(x_0,y)}{|x_0-y|^{d+\alpha}}dy  \nn\\
 	&\quad +   \int_{D\setminus B(x_0,\sigma r)} \left( a^q h_{r,\sigma}(y) - H_1(y) \right) \frac{K(x_0,y)}{|x_0-y|^{d+\alpha}}dy  \nn\\
 	&=:I_1+I_2+I_3+I_4.
 \end{align*}

 For any $y \in B(x_0,bs)$, we have
 \begin{align*}
 	& 	\partial_{i} \left(  a^q h_{r,\sigma}   - H_1 \right)(y)= \frac{aq}{r^{q} }    \left(a^{q-1}\rho(y)^{q-1}\partial_{i} \rho(y)  -f(y )^{q-1} \partial_i \rho(x_0)  \right) 
 \end{align*}
 and 
 \begin{align*}
 	\partial_{ij}\left(  a^q h_{r,\sigma}  - H_1 \right)(y)&	 =  \frac{aq}{r^{q} }    \left(a^{q-1}\rho(y)^{q-1}\partial_{ij} \rho(y)  -f(y )^{q-1} \partial_{ij} \rho(x_0)  \right)  \\
 	&\quad  + \frac{a^2q(q-1)}{r^{q}}    \left(a^{q-2}\rho(y)^{q-2}\partial_{i} \rho(y)  \partial_{j} \rho(y)  -f(y )^{q-2} \partial_i \rho(x_0) \partial_j \rho(x_0)  \right).
 \end{align*} 
 Repeating the arguments leading to \eqref{e:barriers-critical-3} and \eqref{e:barriers-critical-4}, 
 we get that  for any $y \in B(x_0,bs)$,
 \begin{gather*}
 	\left| a^q h_{r,\sigma}(y)   - H_1(y)\right| \le \frac{c_1  s^{q-1}\ell(s) |y-x_0| }{r^q} ,\nn\\
 	 	\left| a^q h_{r,\sigma}(y)  - H_1(y) - (y-x_0) \cdot \nabla \left(a^q h_{r,\sigma}    - H_1\right)(x_0)  \right| \le \frac{c_{2}  s^{q-2} \ell(s) |y-x_0|^2 }{r^{q}} .
 \end{gather*}
 By \eqref{e:censored-1} and \eqref{e:censored-2}, we have
 \begin{align*}
 	|I_1| &\le \frac{c_3 s^{q-1}\ell(s)}{r^{q}} \int_{B(x_0,bs)} \frac{dy}{|x_0-y|^{d+\alpha-\theta-1}} = \frac{c_4 s^{q+\theta-\alpha}\ell(s)}{r^{q}},\\
 		|I_2| &\le \frac{c_{5}  s^{q-2}\ell(s)  }{r^{q}}\int_{B(x_0,bs)} \frac{dy}{|x_0-y|^{d+\alpha-2}}  = \frac{c_{6}s^{q-\alpha}\ell(s)}{r^{q}}.
 \end{align*}
 For $I_4$, by \eqref{e:censored-1} and \eqref{e:barriers-critical-1},  we have
  \begin{align*}
 	|I_4|&\le c_{7} \int_{B(x_0,\sigma r)^c} \frac{( \sigma^q+  (2|x_0-y|/r)^q)}{|x_0-y|^{d+\alpha}} dy  = \frac{c_8 \sigma^{q-\alpha}}{r^\alpha} .
 \end{align*}
 For $I_3$, by repeating the arguments for \eqref{e:barriers-I3-case1}--\eqref{e:barriers-I33}, we get that if $q\ge 1$, then
 \begin{align*}
 |I_3|& \le \frac{c_{9}}{r^q} \int_{B(x_0,\sigma r)\setminus B(x_0,bs)} \frac{\ell(|x_0-y|)}{|x_0-y|^{d+\alpha-q}} dy \nn\\
 &\le \frac{c_{9}\ell(bs)}{r^q (bs)^{(\alpha-q)/2}} \int_{B(x_0,\sigma r)\setminus B(x_0,bs)} \frac{dy}{|x_0-y|^{d+(\alpha-q)/2}} \le \frac{c_{10} s^{q-\alpha}\ell(s)}{r^q }
 \end{align*}
 and if $q<1$, then
 \begin{align*}
 	|I_3|& \le \frac{c_{10}}{r^q}  \int_0^{r} \int_0^{r}  \frac{\ind_{\{l^2+|s-y_d|^2>(bs)^2\}} \ell((l^2+|s-y_d|^2)^{1/2} )}{y_d^{1-q} (l^2+|s-y_d|^2)^{(d+\alpha-1)/2}} \, l^{d-2 } dldy_d \nn\\
 		& \le \frac{c_{10}}{r^q}  \int_0^{r}  \frac{\ind_{\{|s-y_d|^2>(bs)^2/2\}} \ell(2^{1/2} |s-y_d| )}{y_d^{1-q} |s-y_d|^{d+\alpha-1}} \int_0^{|s-y_d|}  l^{d-2 } dldy_d  \nn\\
 		&\quad +\frac{c_{10}}{r^q}\int_0^{r} \frac{\ind_{\{l^2>(bs)^2/2\}}\ell(2^{1/2}l )}{ l^{\alpha+1 }}    \,\int_{0}^{r}   \frac{\ind_{\{|s-y_d|<l\}}}{y_d^{1-q}} dy_d dl \nn\\
 			& \le \frac{c_{11}}{r^q}  \int_0^{(1-2^{-1/2}b)s}  \frac{ \ell( |s-y_d| )}{y_d^{1-q} |s-y_d|^{\alpha}} dy_d   + \frac{c_{11}}{r^q}  \int_{(1+2^{-1/2}b)s} ^{r}  \frac{ \ell( |s-y_d| )}{y_d^{1-q} |s-y_d|^{\alpha}} dy_d  \nn\\
 		&\quad +\frac{c_{11}}{r^q}\int_{2^{-1/2}bs}^{r} \frac{\ell(l )}{ l^{\alpha+1 }}    \,\int_{0}^{l+s}   \frac{dy_d}{y_d^{1-q}}  dl\nn\\
 			& \le \frac{c_{12}s^{q-\alpha}\ell(s)}{r^q} .
 \end{align*}

By Lemma \ref{l:critical-power-functions}(i), $|\sL^\kappa H_1(x_0)| \le c_{13} r^{-q} s^{q-\alpha}\ell(s)$.
 Hence,  combining the estimates for $|I_1|$, $|I_2|$, $|I_3|$ and $|I_4|$, we obtain \eqref{e:barriers-critical-result-2}. The proof is complete.\qed

 \begin{cor}\label{c:critical-barrier}
	There exist $r_2 \in (0,r_1/4]$ and  comparison constants depending only on $\chi$ such that 
 	\begin{align*}
 		\P_{x} (X_{\tau_{D(r)}} \in D)  \asymp ( \delta_D(x) /r)^{q} \quad \text{for  $r\in (0,r_2]$ and $x\in D(r)$.} 
 	\end{align*}
 \end{cor}
 \pf  
  Let $R\in (0,r_1]$ and let $\sigma \in (0,1/4]$ be a constant to be determined later. Define $h_{R,\sigma}$ and $\psi_{R, \sigma}$ by \eqref{e:def-h-r-sigma-critical} and \eqref{e:def-psi-r-sigma-critical}.
  By Proposition  \ref{p:barriers-critical}, for all $\eta\in (0,\sigma)$ and  $x\in D(\eta R)$, we have
 \begin{align}\label{e:critical-barrier-generator}
 	 |	\sL^\kappa h_{R,  \sigma}(x)| &\le  c_1 R^{-q} \delta_D(x)^{q-\alpha}  \ell(\delta_D(x)/R)  + c_2 \sigma^{q-\alpha} R^{-\alpha},\nn\\
 	\sL^\kappa \psi_{R, \sigma }(x) &\ge  \left( c_3  - \frac{c_4 F_\ell(\eta)}{\ell( \sigma)}  \right)   R^{-q} \delta_D(x)^{q-\alpha} \ell(\delta_D(x)/R).
 \end{align} By  \eqref{e:regular-distance-1}, \eqref{e:inward-normal-result-1} and the scaling of $F_\ell$,  $\psi_{R, \sigma}(x)/h_{R, \sigma}(x) \le c_5 F_\ell( \sigma)$ for all $x \in D$. Thus,  by choosing $ \sigma$ small enough,  we obtain  $8c_1c_3^{-1}\psi_{R, \sigma}(x) \le h_{R,\sigma}(x)$ for all $x \in D$.

 Note that $s^{q-\alpha}\ell(s)$ is decreasing and  $\lim_{s\to 0} s^{q-\alpha}\ell(s) \ge \lim_{s\to 0} s^{(q-\alpha)/2}\ell(1) = 0$.
 Thus, there exists $\sigma_0 \in (0,\sigma)$ such that   $c_2\sigma^{q-\alpha} \le c_1 \inf_{s\in (0,\sigma_0]} s^{q-\alpha}\ell(s)$ and  $c_4F_\ell(\sigma_0) \le c_3\ell(\sigma)/2$. By \eqref{e:critical-barrier-generator}, we obtain for all  $x \in D(\sigma_0R)$,
 \begin{align}\label{e:critical-barrier-1}
 	|	\sL^\kappa h_{R,  \sigma}(x)| &\le   c_1 R^{-q} \delta_D(x)^{q-\alpha}  \ell(\delta_D(x)/R)  + c_1R^{-\alpha} \inf_{s\in (0,\sigma_0R]} (s/R)^{q-\alpha} \ell(s/R) \nn\\
 	&\le 2  c_1 R^{-q} \delta_D(x)^{q-\alpha}  \ell(\delta_D(x)/R),\nn\\
 			\sL^\kappa \psi_{R, \sigma }(x) &\ge 2^{-1}c_3 R^{-q} \delta_D(x)^{q-\alpha} \ell(\delta_D(x)/R).
 \end{align}
Define $u_1(x):=h_{R,\sigma}(y) + 4c_1c_3^{-1}\psi_{R, \sigma}(y)$ and $u_2(y):=h_{R,\sigma}(y) -4c_1c_3^{-1}\psi_{R, \sigma}(y)$. Then we have
 $$(1/2)h_{R,\sigma}(y) \le  u_2(y) \le u_1(y) \le (3/2) h_{R,\sigma}(y) \quad \text{for all $y\in D$}.$$
 Hence, by  \eqref{e:regular-distance-1}, there exists $c_6>1$ such that $c_6^{-1} \le u_2(y) \le u_1(y) \le c_6$ for all $y \in D\setminus D(\sigma_0 R)$. Moreover, by \eqref{e:critical-barrier-1},  $\sL^\kappa u_1 \ge 0$ in $D(\sigma_0R)$ and $\sL^\kappa u_2 \le 0$ in $D(\sigma_0R)$. 
 Thus, 
 by Corollary \ref{c:maximum-censored}, we obtain for all $x\in D(\sigma_0R)$,
 \begin{align*}
 	u_1(x) \le \E_x\big[ u_1(X_{\tau_{D(\sigma_0R)}}) \big] = \E_x\big[ u_1(X_{\tau_{D(\sigma_0R)}}): X_{\tau_{D(\sigma_0R)}} \in D \setminus D(\sigma_0R) \big] \le c_6 \P_x(X_{\tau_{D(\sigma_0r)}} \in D ) 
 \end{align*}
 and
 \begin{align*}
 	u_2(x) \ge \E_x\big[ u_2(X_{\tau_{D(\sigma_0R)}}) \big] = \E_x\big[ u_2(X_{\tau_{D(\sigma_0R)}}): X_{\tau_{D(\sigma_0R)}} \in D \setminus D(a_0R) \big] \ge c_6^{-1} \P_x(X_{\tau_{D(\sigma_0R)}} \in D ) .
 \end{align*}
 Since $u_1(x) \asymp u_2(x) \asymp 	(\delta_D(x)/R )^{q}$ for  $x\in D(\sigma_0R)$, this implies the  result with $r_2=\sigma_0r_1$.\qed

 \medskip

 \noindent \textbf{Proof of Theorem \ref{t:critical}.}  
 Recall that  $\overline Y^D$ satisfies   {\bf (A)} with $\phi(r)=r^\alpha$ and $T_0=1$.  
   Thus, by
  Theorems \ref{t:factorization-heatkernel}, \ref{t:largetime-bounded} and \ref{t:factorization-Green}, and Proposition \ref{p:boundary-term-scale}, 
 since $X$ is symmetric, 
 it suffices to show that
 \begin{align}\label{e:critical-boundary}
 	\P_{x} (\zeta> t)  \asymp  \bigg(1 \wedge\frac{\delta_D(x)}{t^{1/\alpha}} \bigg)^{q},
	\quad (t,x) \in (0,1]\times D. 
 \end{align}	

Let $r_2$ be the constant in Corollary \ref{c:critical-barrier}.  By Proposition \ref{p:boundary-term-scale} and  Corollary \ref{c:boundary-term-equivalent},  we have \begin{align}\label{e:critical-boundary-1}	\P_{x} (\zeta> t) \asymp  \P_{x} \big(X_{\tau_{D(r_2t^{1/\alpha})}} \in D\big)  \quad \text{ for  $(t,x) \in (0,1]\times D$.} \end{align}  	 For any  $(t,x) \in (0,1]\times D$, if  $\delta_D(x)\ge r_2t^{1/\alpha}$, then by Lemma \ref{l:survival-interior}(i), $1\ge 	\P_x(\zeta>t) \ge \P_x(\tau_{B(x,r_2t^{1/\alpha})}>t) \ge c_1$. Thus,  $\P_x(\zeta>t) \asymp 1 \asymp (1 \wedge (\delta_D(x)/t^{1/\alpha})^{q}$.  If $\delta_D(x)<r_2t^{1/\alpha}$, then by \eqref{e:critical-boundary-1} and Corollary \ref{c:critical-barrier}, we get $	\P_{x} (\zeta> t) \asymp (\delta_D(x)/t^{1/\alpha})^{q}$.  Hence, \eqref{e:critical-boundary} holds and the proof is complete.\qed

	\subsection{Non-symmetric stable processes in $C^{1,2{\text - \rm Dini}}$ open sets}\label{ss:5.2}
	
	Let  $\alpha \in (0,2)$, $d\ge 2$ and let $Y$ be a  strictly  $\alpha$-stable L\'evy process on $\R^d$, that is,  for all $a>0$, $(Y_{at})_{t\ge 0}$ and $(a^{1/\alpha}Y_{t})_{t\ge 0}$ have the same distribution. 	By \cite[Theorems 14.3 and 14.7]{Sa13}, the infinitesimal generator of $Y$ restricted to $C_c^2(\R^d)$ coincides with a pseudo-differential operator  $\sL$  of the following form: 
	\begin{align*}
		\sL u(x)= \begin{cases}
			\int_{\R^d} \left(u(x+y)-u(x)\right) \nu(dy) &\mbox{ if $\alpha \in (0,1)$},\\[3pt]
			p.v.  \int_{\R^d} \left(u(x+y)-u(x)\right) \nu(dy) + \mu \cdot \nabla u(x)&\mbox{ if $\alpha=1$},\\[3pt] \int_{\R^d} \left(u(x+y) -u(x) - \nabla u(x) \cdot y\right) \nu(dy) &\mbox{ if $\alpha\in (1,2)$},
		\end{cases}
	\end{align*}
	where $\mu \in \R^d$ and $\nu$ is a non-negative homogeneous measure of degree $-d-\alpha$. That is,  there exists a non-negative finite measure $m(d\xi)$ on the unit sphere $\bS^{d-1}:=\{x\in \R^d:|x|=1\}$ such that for all Borel set $B\subset \R^d$,
	\begin{align*}
		\nu(B)=\int_{\bS^{d-1}}m(d\xi) \int_0^\infty \ind_B(r\xi) \frac{dr}{r^{1+\alpha}}.
	\end{align*}
	Additionally, when $\alpha=1$, the following cancellation property also holds:
	\begin{align*}
		\int_{\bS^{d-1}} \xi \, m(d\xi)=0.
	\end{align*}
	Throughout this subsection, we assume that $m(d\xi)$ has a density $m(\xi) $ with respect to the surface measure on $\bS^{d-1}$ and there exist 
	$A_6\ge A_5>0$ such that
\begin{align}\label{e:nonsymmetric-density}
	 A_5\le m(\xi)\le A_6 \quad \text{for all $\xi \in \bS^{d-1}$}.
\end{align}
As an immediate consequence of \eqref{e:nonsymmetric-density}, we see that the measure $\nu(dx)$ has a density $\nu(x)$ with respect to the Lebesgue measure on $\R^d$ and, for some $C\ge 1$,
\begin{align}\label{e:nonsymmetric-nu-density}
C^{-1}|x|^{-d-\alpha}\le \nu(x)\le C|x|^{-d-\alpha} \;\; \text{for all $x \in \R^{d}$}.
\end{align}
Moreover,  by  \cite[Theorem 1.5 and Remark 1.6 and (2.3)]{Wa07},  $Y$ has a smooth transition density $q(t,x,y)$ with respect to the Lebesgue measure such that
\begin{align*}
	q(t,x,y)  \asymp t^{-d/\alpha} \wedge \frac{t}{|x-y|^{d+\alpha}} \quad \text{for all $(t,x,y) \in (0,\infty) \times \R^d \times \R^d$}.
\end{align*}
Thus, $Y$ satisfies the hypothesis {\bf (A)} with $\phi(r)=r^\alpha$ and $T_0=\infty$. 

 Let $\Psi$ be the L\'evy exponent of $Y$ defined by $	\E[ e^{i\xi \cdot  Y_t }] = e^{-t \Psi(\xi)}$ for $\xi \in \R^d$, $t>0$. By \cite[Lemma 2.2]{DRSV22}, we have
\begin{align*}
	\text{Re}\,\Psi(\xi) &= |\Gamma(-\alpha)| \cos (\pi \alpha/2) \int_{\bS^{d-1}} | \theta \cdot \xi  |^ \alpha \, \nu_e(\theta)d\theta,\\
	\text{Im}\,\Psi(\xi) &= \begin{cases}
		- |\Gamma(-\alpha)| \sin(\pi \alpha/2)\int_{\bS^{d-1}} | \theta \cdot \xi  |^ {\alpha-1} ( \theta \cdot \xi)\, \nu_o(\theta)d\theta &\mbox{ if $\alpha \neq 1$},\\[3pt]
		\int_{\bS^{d-1}} ( \theta \cdot \xi) \, \log | \theta \cdot  \xi |\, \nu_o(\theta)d\theta  -  \mu \cdot \xi  &\mbox{ if $\alpha = 1$},
	\end{cases} 
\end{align*}
where $\Gamma(z)$ denotes the Gamma function, and  $\nu_e(x)$ and $\nu_o(x)$ are the even and odd parts of $\nu(x)$, namely,
\begin{align*}
	\nu_e(x):= \frac{\nu(x) +\nu(-x)}{2} \quad \text{and} \quad 	\nu_o(x):= \frac{\nu(x) -\nu(-x)}{2}.
\end{align*}
 Define functions $\gamma$ and $\wh \gamma$ on $\bS^{d-1}$ by
\begin{align}\label{e:def-gamma}
	\gamma(\theta)&= \frac{\alpha}{2} + \frac{1}{\pi} \arctan \left( \frac{	\text{Im}\,\Psi(\theta)}{	\text{Re}\,\Psi(\theta)}\right) ,\qquad 
	\wh \gamma(\theta)=\alpha - \gamma(\theta).
\end{align}
 Note that if $\nu$ is symmetric, namely, $\nu(x)=\nu(-x)$ for all $x\in \R^d$, then  $\gamma(\theta)=\wh\gamma(\theta)=\alpha/2$ for all $\theta \in \bS^{d-1}$. By \cite[Lemma 2.3 and Propositions 4.3 and 4.4]{DRSV22}, we have $\gamma(\theta), \wh \gamma(\theta) \in ((\alpha-1)_+, \alpha \wedge 1)$ for all $\theta \in \bS^{d-1}$.

Let $D\subset \R^d$ be a $C^{1,2{\text - \rm Dini}}$ open set with characteristics $(r_0,\ell_0)$.
Recall that  for  $x\in D$, $Q_x\in \partial D$ is chosen so that  $\delta_D(x) = |x-Q_x|$, and   for  $Q\in \partial D$, $n_Q$ is  the inward unit normal to $\partial D$ at $Q$.  
 Let $X:=Y^D$ be  the killed process of $Y$ on $D$ and    let  $\sA$ be the infinitesimal generator of $X$.  In this subsection, we establish the following two-sided estimates for the transition density and the Green function of $X$. The proof of the next theorem will be presented at the end of this subsection.
	\begin{thm}\label{t:nonsymmetricLevy}
 Let $\alpha \in (0,2)$, $d\ge 2$ and  $D\subset \R^d$  be a  $C^{1,2{\text - \rm Dini}}$  open set with characteristics $(r_0,\ell_0)$.
 Suppose that  \eqref{e:nonsymmetric-density} holds. Then the process $X$ has a  transition density $p(t,x,y)$ on $(0,\infty) \times D \times D$. Moreover, the following two assertions hold.
	
	\smallskip
		
\noindent 		(i) For every $T>0$, there are comparison constants such that for  all $(t,x,y) \in (0,T]\times D\times D$,
		\begin{align*}
			p(t,x,y) \asymp \bigg(1 \wedge\frac{\delta_D(x)}{t^{1/\alpha}} \bigg)^{\gamma(n_{Q_x})}\bigg(1 \wedge\frac{\delta_D(y)}{t^{1/\alpha}} \bigg)^{\wh\gamma(n_{Q_y})} \bigg( t^{-d/\alpha} \wedge \frac{t}{|x-y|^{d+\alpha}}\bigg).
		\end{align*}

		\noindent 		(ii) Suppose in addition that $D$ is bounded. Then there exists $C\ge 1$ such that for all $(t,x,y) \in [3,\infty)\times D\times D$,
		\begin{align*}
	C^{-1}\delta_D(x)^{\gamma(n_{Q_x})}\delta_D(y)^{\wh\gamma(n_{Q_y})}e^{-\lambda_1t}\le 		p(t,x,y) \le C\delta_D(x)^{\gamma(n_{Q_x})}\delta_D(y)^{\wh\gamma(n_{Q_y})}e^{-\lambda_1t},
		\end{align*}
		where $-\lambda_1:=\sup \text{\rm Re\,} \sigma(\sA)<0$  and $\sigma(\sA)$ is the spectrum of  the   generator  $\sA$ of $X$.
			 Moreover, $X$ admits a Green function $G(x,y)$ on $D\times D$  satisfying the following estimate: For  all $(x,y) \in D\times D$,
		\begin{align*}
			G(x,y) \asymp \bigg(1 \wedge\frac{\delta_D(x)}{|x-y|} \bigg)^{\gamma(n_{Q_x})}\bigg(1 \wedge\frac{\delta_D(y)}{|x-y|} \bigg)^{\wh\gamma(n_{Q_y})}  |x-y|^{\alpha-d}.
		\end{align*}
	\end{thm}
	
	\smallskip

	In the remainder of this section, we always assume that we are in the setup of Theorem \ref{t:nonsymmetricLevy}, that is, $\alpha \in (0,2)$, $d\ge 2$, $D$ is a $C^{1,2{\text - \rm Dini}}$ open set with characteristics $(r_0,\ell_0)$, and \eqref{e:nonsymmetric-density} holds.

	 \begin{lemma}\label{l:maximum-nonsymmetric}
	 	Let $U\subset D$ be a 
	 bounded 
	 open set. Suppose that    $u$ is a bounded Borel function on $D$ 
	 such that $u$ is continuous on $\overline{U}\cap D$ and vanishes continuously on $\partial U\cap \partial D$, 
	 $u|_{U} \in C^2(U)$, and $\sL u \ge 0$  on $U$ {\rm (}resp. $\sL u \le 0$  on $U${\rm )}. Then we have $$u(x) \le  \E_x\left[ u(X_{\tau_U})\right] \;\; \text{{\rm (}resp. $u(x) \ge  \E_x\left[ u(X_{\tau_U})\right]${\rm )}} \quad \text{ for all $x\in U$.}$$
	 \end{lemma}
	 \pf  Since the proofs are similar, we only prove  the case where $\sL u \ge 0$ on $D(r)$.
	 
	 We follow the arguments in  \cite[Proposition 5.3 and Corollary 5.4]{CKSV25}.	 For $j \ge 1$, define $U_j:=\{y \in U: \delta_D(y) >2^{-j} \}$. Let $x \in U$ and let $j_0\ge 1$ be such that $\delta_D(x)>2^{-j_0}$.
	 	  Pick an arbitrary $j\ge j_0$  and let $f \in C^2_c(D)$ be such that $f = u$ on $U_{j+1}$. By Ito's formula,
	  \begin{align}\label{e:maximum-nonsymmetric-1}
	  	 	\E_x \big[ f(X_{\tau_{U_j}})\big] 
	  	= u(x) + \E_x \bigg[\int_0^{\tau_{U_j}} \sL f(X_s)ds \bigg].	
	  \end{align}
	 	Define $h:=u-f$. Then $h=0$ on $U_{j+1}$ and $h$ is bounded. Thus, by  \eqref{e:nonsymmetric-nu-density}, 
	 	\begin{align*}
	 		\sup_{x \in \overline U_j}	 \int_{D\setminus U_{j+1}}  h(y)\nu(x-y)dy\le c_1 \lVert h \rVert_{L^\infty(D)} \int_{B(0, 2^{-j-1})^c} \frac{dy}{|y|^{d+\alpha}} <\infty.
	 	\end{align*}  Using the L\'evy system formula \eqref{e:levy_systemX} and \eqref{e:nonsymmetric-nu-density}, we obtain
	 \begin{align}\label{e:maximum-nonsymmetric-2}
	 		&\E_x \big[ h(X_{\tau_{U_j}})  \big]= 	\E_x \big[ h(X_{\tau_{U_j}}) : X_{\tau_{U_j}} \in D\setminus U_{j+1} \big]\nn\\
	 		&= \E_x\bigg[\int_0^{\tau_{U_j}} \int_{D\setminus U_{j+1}}  h(y)\nu(X_s-y)dyds \bigg]	= \E_x\bigg[\int_0^{\tau_{U_j}} \sL h (X_s)ds \bigg].
	 \end{align}
	 Adding \eqref{e:maximum-nonsymmetric-1} and \eqref{e:maximum-nonsymmetric-2},  we arrive at
	 \begin{align}\label{e:maximum-nonsymmetric-conclusion}
	 	\E_x \big[ u(X_{\tau_{U_j}})  \big] = u(x) +  \E_x\bigg[\int_0^{\tau_{U_j}} \sL u(X_s)ds \bigg] \ge u(x).
	 \end{align}
Since $u$ vanishes continuously on $\partial U \cap \partial D$, $$\lim_{j\to \infty} 	\left| \E_x [ u(X_{\tau_{U_j}})  : X_{\tau_{U_j}} \in U \setminus U_j]  \right| \le \lim_{j\to \infty} \sup_{y \in U: \delta_D(y) \le 2^{-j}} |u(y)| =0. $$
Thus, since $U_j \uparrow U$ and $u$ is bounded,  by the bounded convergence theorem,  we get
 $$\lim_{j\to \infty}	\E_x [ u(X_{\tau_{U_j}})  ] = \lim_{j\to \infty}	\E_x [ u(X_{\tau_{U_j}}) : X_{\tau_{U_j}} \in U^c   ] =	\E_x [ u(X_{\tau_{U}})  ].$$  Hence, by taking the limit in \eqref{e:maximum-nonsymmetric-conclusion} as $j\to \infty $,  we obtain the result. \qed

By following the argument of Corollary \ref{c:maximum-censored}, we deduce the result below  from Lemma \ref{l:maximum-nonsymmetric}.
	 \begin{cor}\label{c:maximum-nonsymmetric}
	 	Let $r>0$. Suppose that    $u$ is a bounded Borel function on $D$   		such that $u$ vanishes continuously on $\partial D$ and is continuous on $\overline{D(r)}$, 
		$u|_{D(r)} \in C^2(D(r))$, and $\sL u \ge 0$  on $D(r)$ {\rm (}resp. $\sL u \le 0$  on $D(r)${\rm )}. Then we have $$u(x) \le  \E_x\big[ u(X_{\tau_{D(r)}})\big] \;\;\; \text{{\rm \big(}resp. $u(x) \ge  \E_x\big[ u(X_{\tau_{D(r)}})\big]${\rm \big)}} \quad \text{ for all $x\in D(r)$.}$$
	 \end{cor}
	 
		Recall that  $\gamma(\theta)$ is defined in \eqref{e:def-gamma}.  Set 
	 $$\text{$\gamma_0=\min_{\theta \in \bS^{d-1}}\gamma(\theta)$ \;\; and \;\;  $	\gamma_1=\max_{\theta \in \bS^{d-1}}\gamma(\theta)$}.$$
	 Note that $(\alpha-1)_+< \gamma_0 \le \gamma_1 < \alpha \wedge 1$.
	  By choosing a smaller $r_0$ if necessary, we assume without loss of generality that $\ell_0(r) = \ell_0(r_0) \le 1/16$ for all $r\in [r_0,1]$. Let 
	 $$\eps_0:= \frac{  \gamma_0 \wedge (\alpha-\gamma_1)\wedge  1}{8}.$$ 
	 By Lemma \ref{l:Dini-regularizing}, there exists $\ell \in \DDini \cap C^2((0,\infty))$ satisfying  properties (a)--(d) in Lemma \ref{l:Dini-regularizing} with $f$ replaced by $\ell_0$ and  $\eps$  by $\eps_0$. In particular, 
	  $\ell(1) \le  4\ell_0(1) \le 1/4$.  We extend $\ell$ to $(0,\infty)$ by setting $\ell(r)=\ell(1)$ for $r\ge 1$.

	Define
	\begin{align*}
		L_{\ell}(r):= \int_0^r \frac{\ell(s) (2\gamma_0^{-1}+|\log s|)}{s}ds, \quad r>0.
	\end{align*} 
	Note that $L_{\ell}$ is slowly varying at $\infty$ since $\ell(s)=1$ for $s\ge 1$. 	Since $u^{-\eps_0}\ell(u)$ is decreasing and  $\sup_{a\ge 1, u >0}(2\gamma_0^{-1}+|\log au|)/(a^{\eps_0}(2\gamma_0^{-1}+ |\log u|))<\infty$, we have
	\begin{align*}
		L_\ell(R) =  \int_0^r \frac{\ell(R u/r) (2\gamma_0^{-1}+ |\log Ru/r|)}{u}du \le  c\bigg(\frac{R}{r}\bigg)^{2\eps_0} L_\ell(r) \quad \text{ for all $0<r\le R$.} 
	\end{align*}
	Moreover, 	using $\gamma_0<1$,
	we see that  $L_\ell(r) \ge \int_0^r s^{-1}\ell(s) ds \ge \ell(r)$ for all $r>0$ and
	\begin{align*}
	L_\ell(r) \ge |\log r|  \int_0^r \frac{\ell(s)}{s}ds \ge \ell(r) |\log r|  \quad \text{for all $r\in (0,1]$}.
	\end{align*}

	\begin{prop}\label{p:nonsymmetric-superhamonic}
		Let  $\theta\in \bS^{d-1}$ and 		$\gamma(\cdot)$
		be given by \eqref{e:def-gamma}.  For $q\in (0,\alpha)$ and  $r\in (0,1]$, define 
		$$f_{\theta, q}(x):=(x\cdot \theta)_+^q \quad \text{ and } \quad g_{\theta, q,r}(x):=(x\cdot \theta)_+^q L_{\ell}((x\cdot\theta)_+/r).$$
	\noindent 	(i) There exists $C=C(d, \Psi)>0$ such that 
	\begin{align*}
		|\sL f_{\theta,q}(x) | \le C |q - \gamma(\theta)|  (x\cdot \theta)^{q-\alpha} \quad \text{ for all $q\in (0,\alpha)$ and $x\in \R^d$ with $x\cdot \theta>0$.} 
	\end{align*}
	
		\noindent 	(ii) There exist $\delta_0 =\delta_0(d,\Psi,\ell)\in (0,(\gamma_0 \wedge (\alpha- \gamma_1))/2)$ and  $C=C(d,\Psi,\ell)>0$  such that 
	\begin{align*}
		\sL g_{\theta, q,r} (x)  - L_{\ell}((x\cdot \theta)/r) \sL f_{\theta,q}(x)  \ge C(x\cdot \theta)^{q-\alpha} \ell( (x\cdot \theta)/r)(1+|\log ((x\cdot \theta)/r)| )
	\end{align*}
	 for all $q\in (\gamma(\theta) - \delta_0,\gamma(\theta) + \delta_0)$ and $x\in \R^d$ with $x\cdot \theta\in (0,r]$.
	\end{prop}
	\pf (i) The result follows from \cite[Corollary 4.7]{DRSV22}.
	
	\noindent (ii)  We use a coordinate system where 
	 $\theta = (\wt 0,1)$.  By \cite[Lemma 2.3]{DRSV22}, there exist  $c_1=c_1(\theta),c_2=c_2(\theta)>0$ and $b=b(\theta) \in \R$ such that for any function $u:\R^d\to \R$ of the form $u =u(x_d)$,
		\begin{align}\label{e:operator-nonsymmetric}
			\sL u(x)   &=  \begin{cases}
			\displaystyle \int_{\R}   \left( u(y_d)-  u(x_d) \right) \frac{c_1 \ind_{y_d\ge x_d} + c_2 \ind_{y_d<x_d} }{|x_d-y_d|^{1+\alpha}} dy_d  &\mbox{if $\alpha\in (0,1)$},\\[8pt]
			\displaystyle p.v. \int_{\R}  (u(y_d) - u(x_d)) \frac{c_2}{|x_d-y_d|^{2}} dy_d  +  b \partial_{x_d} u(x) &\mbox{if $\alpha=1$},\\[8pt]
			\displaystyle \int_{\R}   \left( u(y_d) - u(x_d) - (y_d-x_d) \partial_{x_d}u(x) \right) \frac{c_1 \ind_{y_d\ge x_d} + c_2 \ind_{y_d<x_d} }{|x_d-y_d|^{1+\alpha}} dy_d &\mbox{if $\alpha\in (1,2)$},
		\end{cases}
	\end{align}
	whenever the right-hand side is well-defined.  Moreover, $c_3:=\inf_{\xi \in \bS^{d-1}} c_2(\xi)>0$.
	
	Let $q\in (\gamma(\theta) - \delta_0,\gamma(\theta) + \delta_0)$, $r\in (0,1]$ and $x=(\wt x, x_d)\in \R^d$ with $x_d \in (0,r]$.  We write  $f_q=f_{(\wt 0,1),q}$ and  $g_q=g_{(\wt 0,1),q,r}$. Since  $f_q, g_q \in C^2(\R^d_+)$, 
	$\limsup_{r\to \infty} f_q(r)/r^q<\infty$, 
	$q<\alpha$, and   
	$L_\ell$ is slowly varying at $\infty$, 
	the right-hand-side of 
	\eqref{e:operator-nonsymmetric} is well-defined for  $f_q$ and $g_q$. Define
	\begin{align*}
		I_0:=  \begin{cases}
			\displaystyle \int_{0}^\infty    \left( y_d^q-x_d^q \right)\frac{c_1 \ind_{y_d\ge x_d} + c_2 \ind_{y_d<x_d} }{|x_d-y_d|^{1+\alpha}} dy_d  &\mbox{if $\alpha\in (0,1)$},\\[8pt]
			\displaystyle p.v. \int_{0}^\infty  (y_d^q- x_d^q) \frac{c_2}{|x_d-y_d|^{2}} dy_d  +  qbx_d^{q-1}  &\mbox{if $\alpha=1$},\\[8pt]
			\displaystyle \int_{0}^\infty   \left( y_d^q - x_d^q - qx_d^{q-1}(y_d-x_d) \right) \frac{c_1 \ind_{y_d\ge x_d} + c_2 \ind_{y_d<x_d} }{|x_d-y_d|^{1+\alpha}} dy_d  +  \frac{c_2qx_d^{q-\alpha}}{\alpha-1} &\mbox{if $\alpha\in (1,2)$}.
		\end{cases}
	\end{align*}
	Using (i),	by choosing $\delta_0$ sufficiently small,  we obtain
	\begin{align}\label{e:nonsymmetric-superhamonic-0}
	I_0		& = \sL f_q (x_d) + x_d^q \int_{-\infty}^{0}  \frac{c_2}{|x_d-y_d|^{1+\alpha}}dy_d   \ge -	c_4 \delta_0 x_d^{q-\alpha} +  \frac{c_2 x_d^{q-\alpha}}{\alpha}  \ge  \frac{c_3 x_d^{q-\alpha}}{2\alpha}.
	\end{align}
	Moreover, using the fact that
	$$ \partial_{y_d} \left(  y_d^q L_\ell(y_d/r)\right) |_{y_d=x_d} - qx_d^{q-1} L_\ell(x_d/r) = x_d^{q-1} \ell(x_d/r) (2\gamma_0^{-1} + |\log (x_d/r)|), $$
	we see that
	\begin{align*}
		&\sL g_q (x)  - L_\ell(x_d/r) \sL f_q(x) - q^{-1}\ell(x_d/r) (2\gamma_0^{-1} + |\log (x_d/r)|) I_0 \\
		&= \begin{cases}
			\displaystyle  c_1\int_{x_d}^\infty  \frac{H(y_d)  }{|x_d-y_d|^{1+\alpha}} dy_d+ c_2 \int_0^{x_d}  \frac{H(y_d)}{|x_d-y_d|^{1+\alpha}} dy_d &\mbox{ if $\alpha \in (0,1)\cup (1,2)$},\\[8pt]
				\displaystyle c_2\,  p.v.\int_{0}^\infty  \frac{H(y_d)  }{|x_d-y_d|^{2}} dy_d &\mbox{ if $\alpha =1$},
		\end{cases}
	\end{align*} 
	where 
	$$
	H(y_d):= y_d^q  \left( L(y_d/r)-  L(x_d/r) \right) - q^{-1} \ell(x_d/r)(2\gamma_0^{-1}+|\log (x_d/r)|)(y_d^q-x_d^q).
	$$
	For all $s>0$, we have
	\begin{align*}
	\frac{H'(s)}{s^{q-1}} & = q\left( L(s/r)-  L(x_d/r) \right) + \ell(s/r) (2\gamma_0^{-1}+ |\log (s/r)|) -  \ell(x_d/r)(2\gamma_0^{-1}+ |\log (x_d/r)|) \\
		&= \int_{x_d/r}^{s/r} \left(\frac{q\ell (u) (2\gamma_0^{-1}+|\log u|)}{u}  + \ell'(u) (2\gamma_0^{-1}+ |\log u|) -  \frac{\ell(u)}{u}\ind_{u<1}   + \frac{\ell(u)}{u}\ind_{u\ge 1}  \right)  du.
	\end{align*}
	Since $q>\gamma_1-\theta \ge \gamma_0/2$, the integrand in the integral above is  positive for all $u>0$. Thus, $H$ is increasing on $(x_d,\infty)$ and is decreasing on $(0,x_d)$, implying that   $H(s) \ge H(x_d)= 0$ for all $s>0$. It follows that  $	\sL g_q (x)  - L_\ell(x_d/r) \sL f_q(x) \ge q^{-1}\ell(x_d/r) (2\gamma_0^{-1} + |\log (x_d/r)|) I_0$. 
	Now the the desired result follows from \eqref{e:nonsymmetric-superhamonic-0}.
	\qed

  Let $\rho $ be the regularized distance for $D$  defined at the beginning of this section, and let $r_1$ be the constant from Lemma \ref{l:inward-normal}. Since  $\Psi$ is $C^{1}(\bS^{d-1})$, we have $|\gamma(n_{Q})- \gamma(n_{Q'})| \le c \ell(|Q-Q'|)$ for all $Q,Q'\in \partial D$. Thus, by an extension theorem (see \cite[VI.2.2]{St70}),     there exist  a non-negative function $\Phi \in C(\R^d)$ and $C=C(d,\Psi,\ell)>0$
  such that 
  \begin{align}
  &	\Phi(Q) = \gamma(n_Q) \quad \text{for all $Q\in \partial D$}, \label{e:exponent-Phi-1}\\
  & |\Phi(x) - \Phi(y) | \le C \ell(|x-y|)\quad \text{for all $x,y \in \R^d$}, \label{e:exponent-Phi-2}\\
  &|D^k \Phi (x)| \le C |\wt\delta_D(x)|^{-k}\ell (|\wt \delta_D(x)|) \quad \text{for all $k=1,2$ and $x \in \R^d\setminus \partial D$},\label{e:exponent-Phi-3}
  \end{align}
  where $\wt \delta_D(x)$ is the signed distance to $\partial D$ defined by \eqref{e:signed-distance}.

For $r\in (0,r_2]$ and $\sigma \in (0,1/4]$, define
$h_{r,\sigma}, \psi_{r,\sigma}:\R^d \to [0,\infty)$ by
		\begin{align}\label{e:def-h-r-sigma-nonsymmetric}
		h_{r,\sigma}(x) := \begin{cases}
			\displaystyle	\bigg(\frac{\rho(x)}{r}\bigg)^{\Phi(x)}    &\mbox{ if $x\in D(r)$},\\[3pt]
		1  &\mbox{ if $x\in D\setminus D(r)$},\\
			0 &\mbox{ if $x\in D^c$},
		\end{cases} 
	\end{align} 
	\begin{align}\label{e:def-psi-r-sigma-nonsymmetric}
		\psi_{r,\sigma}(x) := \begin{cases}
			\displaystyle	\bigg(\frac{\rho(x)}{r}\bigg)^{\Phi(x)}   L_\ell\left(\frac{\rho(x)}{|\nabla \rho(x)|\,r}\right)    &\mbox{ if $x\in D(2\sigma r)$},\\[3pt]
			\displaystyle	\bigg(\frac{\rho(x)}{r}\bigg)^{\Phi(x)} L_\ell(\sigma) &\mbox{ if $x\in D(r)\setminus D(2\sigma r)$},\\
			L_\ell(\sigma)  &\mbox{ if $x\in D\setminus D( r)$},\\
			0 &\mbox{ if $x\in D^c$}.
		\end{cases} 
	\end{align}
	
	The next proposition is motivated by \cite[Proposition 4.10]{DRSV22}. The proof is similar to that of Proposition \ref{p:barriers-critical}.

\begin{prop}\label{p:barriers-nonsymmetric}
There exist constants $r_2 \in (0,r_1]$ and $C_1,C_2,C_3,C_4>0$  depending only on $d,r_0,\Psi,\ell$   such that for any $r\in (0,r_2]$,  $\sigma \in (0,1/4]$ and $x\in D(\sigma r)$,
	\begin{align}
		|\sL h_{r,\sigma} (x)|&\le
		C_1
		r^{-\Phi(x)} \delta_D(x)^{\Phi(x)-\alpha} \ell(\delta_D(x)) |\log (\delta_D(x)/r)| +  
		C_2
		\sigma^{\gamma_0-\alpha} r^{-\alpha},\label{e:barriers-nonsymmetric-result-1}\\
		\sL  \psi_{r,\sigma}(x)& \ge\left( 
		C_3- \frac{C_4 
		L_\ell(\delta_D(x)/r)}{\sigma^{\gamma_1-\gamma_0}\ell(\sigma)}  \right)   r^{-\Phi(x)} \delta_D(x)^{\Phi(x)-\alpha} \ell(\delta_D(x)/r)  |\log (\delta_D(x)/r)| . \label{e:barriers-nonsymmetric-result-2}
	\end{align}
\end{prop}
\pf  Let $r_2\in (0,r_1]$ be a constant to be determined later, and let $r\in (0,r_2]$ and $\sigma \in (0,1/4]$.  Pick an arbitrary  $x_0 \in D(\sigma r)$. Let   $s:=\delta_D(x_0)$, $a:=|\nabla \rho(x_0)|^{-1}$ and $\theta_0:= a\nabla \rho(x_0)$.  Define 
\begin{align*}
	\wt \rho(y):= \frac{\rho(y)}{|\nabla \rho(y)|} \quad \text{ for $y \in D(r_1)$},
\end{align*}
and
\begin{align*} 
	f(y) &:=a(\rho(x_0) + \nabla \rho(x_0) \cdot (y-x_0) )_+,   \\
	H_1(y)&:=(f(y)/r)^{\Phi(x_0)} \qquad \text{and} \qquad H_2(y):= (f(y)/r)^{\Phi(x_0)}L_\ell(f(y)/r)  \quad \text{ for $y \in \R^d$}.
\end{align*} 
Recall that,  by \eqref{e:regular-distance-1} and \eqref{e:inward-normal-result-1}, we have  $a\in [1/2,4]$ and there exists $b\in (0,1/2]$ such that $f(y)/(a\rho(x_0)) \in [1/2,2]$ for all $y \in B(x_0, bs)$. Further,   \eqref{e:barriers-critical-1} and \eqref{e:barriers-critical-2} hold.

Let $\delta_0$ be the constant in Proposition \ref{p:nonsymmetric-superhamonic}(ii). By Lemma \ref{l:inward-normal} and \eqref{e:exponent-Phi-1},   $\Phi(Q)  =\gamma(n_Q)= \gamma(\nabla \rho(Q)/|\nabla \rho(Q)|)$ for all $Q \in \partial D$. Hence, since  $|x_0-Q_{x_0}| =s \le \sigma r_2$,  by choosing  $r_2$ sufficiently small and using \eqref{e:regular-distance-2}, \eqref{e:exponent-Phi-2} and the fact that $\psi \in C^1(\bS^{d-1})$, we obtain
\begin{align*}
	| \Phi(x_0)  - \gamma( \theta_0) | \le 	| \Phi(x_0)  -\Phi(Q_{x_0})| + \left| \gamma (\nabla \rho(Q)/|\nabla \rho(Q)| )- \gamma( \theta_0) \right| & \le  c_1 \ell(s) <\delta_0. 
\end{align*}
Applying Proposition \ref{p:nonsymmetric-superhamonic}, we get that $	\left| 	\sL H_1 (x_0) \right| \le c_2r^{-\Phi(x_0)}   f(x_0)^{\Phi(x_0)-\alpha}\ell(s)$ and 
 \begin{align*}
 			\sL H_2(x_0)  -L_{\ell}(f(x_0)/r) \sL H_1(x_0) &\ge  c_3 r^{-\Phi(x_0)} f(x_0)^{\Phi(x_0)-\alpha} \ell( f(x_0)/r)  (1+ |\log (f(x_0)/r)|).
 \end{align*}
 Since $f(x_0) = a\rho(x_0) \asymp s$ and $\ell(s)\le \ell(s/r)$, by the scaling properties of $L_\ell$ and $\ell$, we have
 \begin{align}
 	\left| 	\sL H_1 (x_0) \right| &\le c_4r^{-\Phi(x_0)} s^{\Phi(x_0)-\alpha}  \ell(s),\label{e:barriers-nonsymmetric-H2}\\
 	\sL H_2(x_0) & \ge  \left(  c_5 |\log (s/r)|  - c_6 L_{\ell}(s/r)\right) r^{-\Phi(x_0)}s^{\Phi(x_0)-\alpha}  \ell(s/r).\label{e:barriers-nonsymmetric-H1}
 \end{align}

 (i)  We first prove \eqref{e:barriers-nonsymmetric-result-2}. Since $a^{\Phi(x_0)}\psi_{r,\sigma}(x_0)=H_2(x_0)$, we have
 \begin{align*}
 &	\sL \left( a^{\Phi(x_0)} \psi_{r,\sigma}  -  H_2 \right)(x_0)\\
  &= \int_{B(x_0,bs)} \left(a^{\Phi(x_0)} \psi_{r,\sigma} (y) - H_2(y) -\ind_{\alpha \ge 1} (y-x_0) \cdot \nabla \left(a^{\Phi(x_0)} \psi_{r,\sigma}   - H_2\right)(x_0)  \right) \nu ( y-x_0)dy\\
  &\quad + \int_{D \cap (B(x_0,\sigma r) \setminus B(x_0,bs))}  \left(a^{\Phi(x_0)} \psi_{r,\sigma} (y) - H_2(y) \right) \nu ( y-x_0)dy \\
   &\quad + \int_{D\setminus B(x_0,\sigma r)} \left(a^{\Phi(x_0)} \psi_{r,\sigma} (y) - H_2(y) \right) \nu ( y-x_0)dy  + \ind_{\alpha =1} \mu \cdot \nabla \left(a^{\Phi(x_0)} \psi_{r,\sigma}   - H_2\right)(x_0)\\
  &\quad   -\ind_{\alpha>1} \int_{B(x_0,bs)^c} (y-x_0) \cdot  \nabla \left(a^{\Phi(x_0)} \psi_{r,\sigma}   - H_2\right)(x_0)  \, \nu(y-x_0)dy\\
 	& =:I_1+I_2+I_3 + \ind_{\alpha=1} I_4 - \ind_{\alpha>1}I_5.
 \end{align*}
  For $y \in B(x_0,bs)$, we have $f(y) \asymp  \delta_D(y) \asymp \rho(y) \asymp \wt \rho(y)  \asymp s$,
\begin{align*}
	& 	\partial_{i} \left( a^{\Phi(x_0)} \psi_{r,\sigma}  -  H_2 \right)(y)= 	\partial_{i} \bigg( a^{\Phi(x_0)} \bigg(\frac{\rho(\cdot)}{r}\bigg)^{\Phi(\cdot)}   L_\ell\left(\frac{\wt \rho(\cdot )}{r}\right)   - \frac{f(\cdot )^{\Phi(x_0)}}{r^{\Phi(x_0)}} L_\ell  \left(\frac{f(\cdot)}{r}\right)  \bigg)(y) \\
	&=  \frac{a^{\Phi(x_0)}\Phi(y)\rho(y)^{\Phi(y)-1} L_\ell(\wt \rho(y)/r)}{r^{\Phi(y)}}    \partial_i \rho(y) -  \frac{a\Phi(x_0)f(y )^{\Phi(x_0)-1}  L_\ell(f(y)/r)}{r^{\Phi(x_0)}}   \partial_i \rho(x_0)  \\
	&\;  +\frac{a^{\Phi(x_0)}\rho(y)^{\Phi(y)}\ell(\wt \rho(y)/r)(2\gamma_0^{-1} + |\log (\wt \rho(y)/r)|)}{\wt \rho(y) r^{\Phi(y)}}   \partial_i \wt \rho(y) \\
	&\;   - \frac{af(y )^{\Phi(x_0)}\ell(f(y)/r)(2\gamma_0^{-1} + |\log (f(y)/r)|)}{f(y)r^{\Phi(x_0)}}\partial_i \rho(x_0)     + \frac{ a^{\Phi(x_0)} \rho(y)^{\Phi(y)} L_\ell (\wt \rho(y)/r) \log (\rho(y)/r)}{r^{\Phi(y)}} \partial_i\Phi(y)
\end{align*}
 and 
\begin{align*}
& 	\partial_{ij} \left( a^{\Phi(x_0)} \psi_{r,\sigma}  -  H_2 \right)(y)\\
&= \frac{a^{\Phi(x_0)}\Phi(y)\rho(y)^{\Phi(y)-1} L_\ell(\wt \rho(y)/r)}{r^{\Phi(y)}}  \partial_{ij} \rho(y)  + \frac{a^{\Phi(x_0)}\rho(y)^{\Phi(y)-1} L_\ell(\wt \rho(y)/r)}{r^{\Phi(y)}}   \partial_{i} \rho(y) \partial_j \Phi(y) \\
&\quad +    \frac{a^{\Phi(x_0)}\Phi(y) (\Phi(y)-1)\rho(y)^{\Phi(y)-2}  L_\ell(\wt \rho(y)/r) }{r^{\Phi(y)}} \partial_i \rho(y)   \partial_j \rho(y)\\
&\quad - \frac{a^2\Phi(x_0) (\Phi(x_0)-1)f(y )^{\Phi(x_0)-2}  L_\ell(f(y)/r)}{r^{\Phi(x_0)}}  \partial_i \rho(x_0) \partial_j \rho(x_0) \\
&\quad +\frac{a^{\Phi(x_0)}\Phi(y)\rho(y)^{\Phi(y)-1} \ell(\wt \rho(y)/r)(2\gamma_0^{-1} + |\log (\wt \rho(y)/r)|)}{\wt \rho(y)r^{\Phi(y)}}  \partial_i \rho(y)  \partial_j \wt \rho(y) \\
&\quad - \frac{a\Phi(x_0)f(y )^{\Phi(x_0)-1}   \ell(f(y)/r)(2\gamma_0^{-1} + |\log (f(y)/r)|)}{r^{\Phi(x_0)}}\partial_i \rho(x_0)\partial_j \rho(x_0)  \\
&\quad + \frac{a^{\Phi(x_0)}\Phi(y)\rho(y)^{\Phi(y)-1}L_\ell(\wt \rho(y)/r) \log(\rho(y)/r)}{r^{\Phi(y)}}  \partial_i \rho(y)  \partial_j \Phi(y)\\
&\quad + \frac{a^{\Phi(x_0)}\rho(y)^{\Phi(y)}\ell(\wt \rho(y)/r)(2\gamma_0^{-1} + |\log (\wt \rho(y)/r)|)}{\wt \rho(y) r^{\Phi(y)}}   \partial_{ij} \wt \rho(y) \\
&\quad  +\frac{a^{\Phi(x_0)}\Phi(y)\rho(y)^{\Phi(y)-1}\ell(\wt \rho(y)/r)(2\gamma_0^{-1} + |\log (\wt \rho(y)/r)|)}{\wt \rho(y) r^{\Phi(y)}}   \partial_i \wt \rho(y) \partial_j \rho(y) \\
&\quad -  \frac{a^2 \Phi(x_0)f(y )^{\Phi(x_0)-1}\ell(f(y)/r)(2\gamma_0^{-1} + |\log (f(y)/r)|)}{f(y)r^{\Phi(x_0)}}\partial_i \rho(x_0) \partial_j \rho(x_0)     \\
&\quad  +\frac{a^{\Phi(x_0)}\rho(y)^{\Phi(y)}\ell'(\wt \rho(y)/r)(2\gamma_0^{-1} + |\log (\wt \rho(y)/r)|)}{\wt \rho(y) r^{\Phi(y) + 1 }}   \partial_i \wt \rho(y)  \partial_j \wt \rho(y)\\
&\quad   - \frac{af(y )^{\Phi(x_0)}\ell'(f(y)/r)(2\gamma_0^{-1} + |\log (f(y)/r)|)}{f(y)r^{\Phi(x_0) + 1}}\partial_i \rho(x_0) \partial_j \rho(x_0)     \\
	&\quad  -\frac{a^{\Phi(x_0)}\rho(y)^{\Phi(y)}\ell(\wt \rho(y)/r)}{\wt \rho(y)^2 r^{\Phi(y)}}   \partial_i \wt \rho(y) \partial_j \wt \rho(y) + \frac{a^2f(y )^{\Phi(x_0)}\ell(f(y)/r)}{f(y)^2 r^{\Phi(x_0)}}\partial_i \rho(x_0) \partial_j \rho(x_0)     \\
		&\quad  -\frac{a^{\Phi(x_0)}\rho(y)^{\Phi(y)}\ell(\wt \rho(y)/r)(2\gamma_0^{-1} + |\log (\wt \rho(y)/r)|)}{\wt \rho(y)^2 r^{\Phi(y)}}   \partial_i \wt \rho(y)   \partial_j \wt \rho(y) \\
	&\quad   + \frac{a^2f(y )^{\Phi(x_0)}\ell(f(y)/r)(2\gamma_0^{-1} + |\log (f(y)/r)|)}{f(y)^2r^{\Phi(x_0)}}\partial_i \rho(x_0) \partial_j \rho(x_0)     \\
&\quad + \frac{a^{\Phi(x_0)}\rho(y)^{\Phi(y)}\ell(\wt \rho(y)/r)(2\gamma_0^{-1} + |\log (\wt \rho(y)/r)|) \log(\rho(y)/r)}{\wt \rho(y) r^{\Phi(y)}}   \partial_i \wt \rho(y)\partial_j \Phi(y)\\
&\quad + \partial_j \left[ \frac{ a^{\Phi(x_0)} \rho(y)^{\Phi(y)} L_\ell (\wt \rho(y)/r) \log (\rho(y)/r)}{r^{\Phi(y)}} \partial_i\Phi(y)\right].
\end{align*}
On the other hand, by Lemma \ref{l:double-Dini-limit},  we have
\begin{align}\label{e:barriers-nonsymmetric-2}
\sup_{u\in (0,1]}	 u^{-\ell (u)} = \sup_{u\in (0,1]}	 e^{-\ell (u) \log u}  <\infty.
\end{align}
Hence, by \eqref{e:exponent-Phi-2}, we have for any $y\in B(x_0,bs)$,
$$
	1\le r^{-|\Phi(x_0)-\Phi(y)|} \le 	s^{-|\Phi(x_0)-\Phi(y)|} \le s^{-c_{7}\ell(s)}\le c_{8}.
$$
 Thus, using \eqref{e:regular-distance-1}--\eqref{e:regular-distance-3}, \eqref{e:barriers-critical-2}, \eqref{e:exponent-Phi-2},  \eqref{e:exponent-Phi-3}, the scaling properties of $L_\ell$ and $\ell$, and the facts that $u\ell'(u)\le c\ell(u)\le c\ell(u)(1+|\log u|) \le cL_\ell(u)$ and $|u^2\ell''(u)| \le c\ell(u)$ for all $u>0$,  we get that  for any $y \in B(x_0,bs)$,
\begin{align}\label{e:barriers-nonsymmetric-3}
	\left|D \left(a^{\Phi(x_0)} \psi_{r,\sigma}  -  H_2\right) (y)\right|
	&\le  \frac{c_{9}  s^{\Phi(x_0)-1} L_\ell(s/r)  \ell(s) |\log (s/r)| }{r^{\Phi(x_0)}}  
\end{align}
and
\begin{align*}
	\left|D^2 \left(a^{\Phi(x_0)} \psi_{r,\sigma}  -  H_2\right) (y)\right|& \le  \frac{c_{10} s^{\Phi(x_0)-2}  \left( L_\ell(s/r)|\log (s/r)| + \ell(s/r) |\log (s/r)|^2 \right)  \ell(s) }{r^{\Phi(x_0)}} \\
	&\le \frac{c_{11} s^{\Phi(x_0)-2}   L_\ell(s/r)\ell(s) |\log (s/r)|   }{r^{\Phi(x_0)}}.
\end{align*}
Since $a^{\Phi(x_0)}\psi_{r,\sigma}(x_0)=H_2(x_0)$, we deduce that   for any $y \in B(x_0,bs)$,
\begin{align}\label{e:barriers-nonsymmetric-4}
	\left|  a^{\Phi(x_0)} \psi_{r,\sigma}(y)  -  H_2(y) \right| \le  \frac{c_{9} |x_0-y| s^{\Phi(x_0)-1} L_\ell(s/r)  \ell(s) |\log (s/r)|  }{r^{\Phi(x_0)}}    
\end{align}
and
\begin{align}\label{e:barriers-nonsymmetric-5}
&	\left|  a^{\Phi(x_0)} \psi_{r,\sigma} (y)  - H_2(y)- (y-x_0) \cdot \nabla \left(a^{\Phi(x_0)} \psi_{r,\sigma}   - H_2\right)(x_0)   \right| \nn\\
&\le \frac{c_{11}|x_0-y|^2 s^{\Phi(x_0)-2}   L_\ell(s/r)\ell(s) |\log (s/r)|   }{2r^{\Phi(x_0)}}.
\end{align}
Using   \eqref{e:nonsymmetric-nu-density},  \eqref{e:barriers-nonsymmetric-4} if $\alpha<1$ and \eqref{e:barriers-nonsymmetric-5} if $\alpha\ge 1$, we obtain
\begin{align*}
	|I_1|&\le \frac{c_{12}s^{\Phi(x_0)}   L_\ell(s/r)\ell(s) |\log (s/r)|  }{r^{\Phi(x_0)}} \int_{B(x_0,bs)}  \frac{ \left(  ( |x_0-y| /s) \ind_{\alpha<1}  + ( |x_0-y| /s)^2 \ind_{\alpha \ge 1} \right)  }{|x_0-y|^{d+\alpha}} dy\\
	&=  \frac{c_{13}s^{\Phi(x_0)-\alpha}   L_\ell(s/r)\ell(s) |\log (s/r)|  }{r^{\Phi(x_0)}} .
\end{align*}
Moreover, by \eqref{e:barriers-nonsymmetric-3} and  \eqref{e:nonsymmetric-nu-density}, we  obtain
\begin{align*}
	\ind_{\alpha=1}|I_4| &\le \frac{c_{14}  s^{\Phi(x_0)-\alpha} L_\ell(s/r)  \ell(s) |\log (s/r)| }{r^{\Phi(x_0)}} 
\end{align*}
and
\begin{align*}
	\ind_{\alpha>1} |I_5| &\le \frac{c_{15}  s^{\Phi(x_0)-1} L_\ell(s/r)  \ell(s) |\log (s/r)| }{r^{\Phi(x_0)}}  \int_{B(x_0,bs)^c}  \frac{\ind_{\alpha>1} dy}{|x_0-y|^{d+\alpha-1}} \nn\\
	&= \frac{c_{16}  s^{\Phi(x_0)-\alpha} L_\ell(s/r)  \ell(s) |\log (s/r)| }{r^{\Phi(x_0)}} .
\end{align*}

For $I_3$, we note that for all $y \in B(x_0,\sigma r)^c$, by \eqref{e:barriers-critical-1} and the scaling of $L_\ell$,  $\psi_{r,\sigma}(y) \le c( |x_0-y|/r)^{\Phi(y)} L_\ell(\sigma)$ and $H_2(y) \le c (|x_0-y|/r)^{\Phi(x_0)} L_\ell(|x_0-y|/r))$. Hence, using
\eqref{e:nonsymmetric-nu-density}   and the scaling of $L_\ell$ (with the fact that $2\eps_0\le (\alpha-\gamma_1)/2 \le (\alpha-\Phi(x_0))/2$), since $s\le \sigma r$ and  $u^{\gamma_1-\alpha}L_\ell(u/r)\ell(u/r)$ is decreasing, we obtain
\begin{align*}
	|I_3|&\le c_{17} \int_{B(x_0,\sigma r)^c} \frac{( |x_0-y|/r)^{\Phi(y)} L_\ell(\sigma) +  (|x_0-y|/r)^{\Phi(x_0)} L_\ell(|x_0-y|/r))}{|x_0-y|^{d+\alpha}} dy  \\
	&\le c_{17} \int_{B(x_0,\sigma r)^c} \frac{\left(( |x_0-y|/r)^{\gamma_0}  \vee ( |x_0-y|/r)^{\gamma_1} \right) L_\ell(\sigma) +  (|x_0-y|/r)^{\Phi(x_0)} L_\ell(|x_0-y|/r))}{|x_0-y|^{d+\alpha}} dy  \\
	& \le \frac{c_{18}\sigma^{\gamma_0} L_\ell(\sigma)}{(\sigma r)^\alpha}  +  \frac{c_{18}  L_\ell(\sigma)}{\sigma^{(\alpha-\Phi(x_0))/2} r^{(\alpha+\Phi(x_0))/2}} \int_{B(x_0,\sigma r)^c} \frac{dy}{|x_0-y|^{d+(\alpha-\Phi(x_0))/2}} \\
	&\le \frac{c_{19}\sigma^{\gamma_0-\alpha} L_\ell(\sigma)}{r^\alpha} \le \frac{c_{19}L_\ell(s/r)\ell(s/r)}{\sigma^{\gamma_1-\gamma_0}\ell(\sigma)r^\alpha } \bigg(\frac{s}{r}\bigg)^{\gamma_1-\alpha}  \le \frac{c_{19}L_\ell(s/r)\ell(s/r)}{\sigma^{\gamma_1-\gamma_0}\ell(\sigma)r^\alpha } \bigg(\frac{s}{r}\bigg)^{\Phi(x_0)-\alpha}.
\end{align*}

Observe that for any $y \in D\cap (B(x_0,\sigma r)\setminus  B(x_0,bs))$, 
\begin{align*}
		\left|  a^{\Phi(x_0)} \psi_{1,r}(y) - H_2(y)  \right| 	&\le a^{\Phi(x_0)}	 L_\ell( \wt \rho(y)/r) \bigg| \frac{ \rho(y)^{\Phi(y)}}{r^{\Phi(y)}}  - \frac{ \rho(y)^{\Phi(x_0)} }{r^{\Phi(x_0)}}   \bigg| \nn\\
		&\quad + \frac{ a^{\Phi(x_0)} \rho(y)^{\Phi(x_0)}}{r^{\Phi(x_0)}}	\left|  L_\ell( \wt \rho(y)/r) -  L_\ell(a\rho(y)/r)   \right| \nn\\
	&\quad  +  \frac{| a^{\Phi(x_0)} \rho(y)^{\Phi(x_0)} L_\ell(a\rho(y)/r)  - f(y)^{\Phi(x_0)} L_\ell(f(y)/r)|}{r^{\Phi(x_0)}}\\
		&=:g_1(y)+g_2(y)+g_3(y).
\end{align*}
Following the argument for \eqref{e:critical-I3-diff}, since $2\eps_0 < \gamma_0 \le \Phi(x_0)$, we see that 
\begin{align*}
	g_3(y)&\le  \frac{| a^{\Phi(x_0)} \rho(y)^{\Phi(x_0)}  - f(y)^{\Phi(x_0)}|}{r^{\Phi(x_0)}}\bigg|  L_\ell\left(\frac{a\rho(y)}{r}\right) +  L_\ell\left(\frac{f(y)}{r}\right) \bigg| \nn\\
	&\le \frac{c_{20} L_\ell( |x_0-y|/r)| a^{\Phi(x_0)} \rho(y)^{\Phi(x_0)}  - f(y)^{\Phi(x_0)}|}{r^{\Phi(x_0)}}.
\end{align*}
 Hence, repeating the arguments for \eqref{e:barriers-I3-case1} if $\Phi(x_0)>1$,  and for \eqref{e:barriers-I3-case1+}--\eqref{e:barriers-I33} 
 if $\Phi(x_0)\le 1$, 
 and using the fact $u^{-(\alpha-\Phi(x_0))/2} L_\ell(u/r)\ell(u)$ is non-increasing, we obtain
\begin{align}\label{e:nonsymmetric-A3}
\int_{ D\cap (B(x_0,\sigma r)\setminus  B(x_0,bs))} 
g_3(y)\nu(y-x_0)dy \le  \frac{c_{21} s^{\Phi(x_0)-\alpha}L_\ell(s/r)\ell(s/r)}{r^{\Phi(x_0)} }.
\end{align}
For any $y\in D\cap (B(x_0,\sigma r)\setminus  B(x_0,bs))$, using $a\le 4$,  \eqref{e:barriers-critical-1}, the fact that $\sup_{v\in (0,u)} v^{\gamma_0}|\log v| \asymp u^{\gamma_0} |\log u|$ for all $u \in (0,1/4]$, \eqref{e:inward-normal-result-1}, the scaling of $L_\ell$ and  \eqref{e:exponent-Phi-2}, we get
\begin{align*}
g_1(y)&\le 4 L_\ell (\wt \rho(y)/r) |\Phi(y)-\Phi(x_0)| \bigg[ \bigg( \frac{\rho(y)}{r}\bigg)^{\Phi(y)} \vee  \bigg( \frac{\rho(y)}{r}\bigg)^{\Phi(x_0)} \bigg]\left| \log \frac{\rho(y)}{r} \right|   \\
	& \le c_{22}  L_\ell (|x_0-y|/r) |\Phi(y)-\Phi(x_0)|   \bigg[ \bigg( \frac{|x_0-y|}{r}\bigg)^{\Phi(y)} \vee  \bigg( \frac{|x_0-y|}{r}\bigg)^{\Phi(x_0)} \bigg] \left| \log \frac{|x_0-y|}{r} \right|\\
	&\le c_{23} L_\ell(|x_0-y|/r)   \ell(|x_0-y|) \bigg( \frac{|x_0-y|}{r}\bigg)^{\Phi(x_0)- c_{24}\ell(|x_0-y|)} \left| \log \frac{|x_0-y|}{r} \right|.
\end{align*}
By \eqref{e:barriers-nonsymmetric-2},   $(r/|x_0-y|)^{c_{24} \ell(|x_0-y|)} \le (1/|x_0-y|)^{c_{24} \ell(|x_0-y|)} \le c_{25}$. Therefore, using  \eqref{e:nonsymmetric-nu-density} and the scaling properties of $L_\ell$ and $\ell$, since $u^{-(\alpha-\gamma_1)/2}L_\ell(u/r) \ell(u)$ is decreasing,  we obtain
\begin{align}\label{e:nonsymmetric-A1}
		&\int_{D\cap (B(x_0,\sigma r)\setminus  B(x_0,bs))}
		g_1(y)\nu(y-x_0) dy  \nn\\
		&\le c_{26} \int_{B(x_0,bs)^c} \frac{L_\ell(|x_0-y|/r)  \ell(|x_0-y|)}{|x_0-y|^{d+\alpha}} \bigg( \frac{|x_0-y|}{r}\bigg)^{\Phi(x_0)} \left| \log \frac{|x_0-y|}{r} \right| dy\nn\\
			&\le \frac{c_{27}s^{(\Phi(x_0)-\alpha)/2}L_\ell(s/r) \ell(s)}{r^{\Phi(x_0)}} \int_{B(x_0,bs)^c} \frac{1}{|x_0-y|^{d+(\alpha-\Phi(x_0))/2}}  \left| \log \frac{|x_0-y|}{r} \right| dy \nn\\
			&\le \frac{c_{28}s^{\Phi(x_0)-\alpha}L_\ell(s/r) \ell(s) |\log (s/r)|}{r^{\Phi(x_0)}}.
\end{align}
For 
$g_2$,  using $a\le 4$, \eqref{e:inward-normal-result-1}, the scaling of $L_\ell$ and $\ell$, $\ell \le L_\ell$, \eqref{e:regular-distance-2} and
 \eqref{e:barriers-critical-1},  we get that for any $y\in D\cap (B(x_0,\sigma r)\setminus  B(x_0,bs))$,
\begin{align*}
g_2(y)&\le \frac{ c_{29} \rho(y)^{\Phi(x_0)}  \ell ( \rho(y)/r) (1+ |\log \rho(y)/r|)}{r^{\Phi(x_0)}}	\left|  \frac{1}{|\nabla \rho(y)|} - \frac{1}{|\nabla \rho(x_0)|}   \right|\\
	&\le  \frac{ c_{30}|x_0-y|^{\Phi(x_0)}  L_\ell (|x_0-y|/r) \ell(|x_0-y|) (1+ |\log |x_0-y|/r|)}{r^{\Phi(x_0)}}.
\end{align*}
Thus, by repeating the argument for \eqref{e:nonsymmetric-A1},  we obtain
\begin{align}\label{e:nonsymmetric-A2}
	&\int_{D\cap (B(x_0,\sigma r)\setminus  B(x_0,bs))} 
	g_2(y) \nu(y-x_0) dy  \le \frac{c_{31}s^{\Phi(x_0)-\alpha}L_\ell(s/r) \ell(s) |\log (s/r)|}{r^{\Phi(x_0)}}.
\end{align}
By \eqref{e:nonsymmetric-A1}, \eqref{e:nonsymmetric-A2} and \eqref{e:nonsymmetric-A3}, we arrive at $|I_2|\le c_{32}r^{-\Phi(x_0)}s^{\Phi(x_0)-\alpha}L_\ell(s/r) \ell(s) |\log (s/r)|.$

Combining the estimates for $|I_1|, |I_2|, |I_3|, \ind_{\alpha=1}|I_4|, \ind_{\alpha>1}|I_5|$ with \eqref{e:barriers-nonsymmetric-H1}, and using $a^{\Phi(x_0)}\le 4$, $\ell(s)\le \ell(s/r)$ and $\sigma^{\gamma_1-\gamma_0} \ell(\sigma)\le c$, we conclude  that \eqref{e:barriers-nonsymmetric-result-2} holds.

\smallskip

(ii) To obtain \eqref{e:barriers-nonsymmetric-result-1}, we follow the arguments in (i). Using $a^{\Phi(x_0)}h_{r,\sigma}(x_0)=H_1(x_0)$, we obtain
\begin{align*}
	&	\sL \left( a^{\Phi(x_0)} h_{r,\sigma}  -  H_1 \right)(x_0)\\
	&= \int_{B(x_0,bs)} \left(a^{\Phi(x_0)} h_{r,\sigma} (y) - H_1(y) -\ind_{\alpha \ge 1} (y-x_0) \cdot \nabla \left(a^{\Phi(x_0)} h_{r,\sigma}   - H_1\right)(x_0)  \right) \nu ( y-x_0)dy\\
	&\quad + \int_{D \cap (B(x_0,\sigma r) \setminus B(x_0,bs))}  \left(a^{\Phi(x_0)} h_{r,\sigma} (y) - H_1(y) \right) \nu ( y-x_0)dy \\
	&\quad + \int_{D\setminus B(x_0,\sigma r)} \left(a^{\Phi(x_0)} h_{r,\sigma} (y) - H_1(y) \right) \nu ( y-x_0)dy  + \ind_{\alpha =1} \mu \cdot \nabla \left(a^{\Phi(x_0)} h_{r,\sigma}   - H_1\right)(x_0)\\
	&\quad   -\ind_{\alpha>1} \int_{B(x_0,bs)^c} (y-x_0) \cdot  \nabla \left(a^{\Phi(x_0)} h_{r,\sigma}   - H_2\right)(x_0)  \, \nu(y-x_0)dy\\
	& =:I_1'+I_2'+I_3' + \ind_{\alpha=1} I_4' - \ind_{\alpha>1}I_5'.
\end{align*}
Following the arguments in (i), we obtain
\begin{align*}
	|I_1'|+|I_2'|+  \ind_{\alpha=1}| I_4'| + \ind_{\alpha>1}|I_5'| \le c_{33} r^{-\Phi(x_0)} s^{\Phi(x_0)-\alpha} \ell(s) |\log (s/r)|.
\end{align*}
 For $I_3'$,   by \eqref{e:nonsymmetric-nu-density} and  \eqref{e:barriers-critical-1}, we have
\begin{align*}
	|I_3'|&\le c_{34} \int_{B(x_0,\sigma r)^c} \frac{\left(( |x_0-y|/r)^{\gamma_0}  \vee ( |x_0-y|/r)^{\gamma_1} \right)+  (|x_0-y|/r)^{\Phi(x_0)}}{|x_0-y|^{d+\alpha}} dy  \le \frac{c_{35}\sigma^{\gamma_0}}{(\sigma r)^\alpha}.
\end{align*}
Thus, combining \eqref{e:barriers-nonsymmetric-H2} and the fact that $a^{-\Phi(x_0)}\le 2$, we arrive at the desired result.  \qed

	\begin{cor}\label{c:nonsymmetric-barrier}
There exist $r_3 \in (0,r_1]$ and  comparison constants depending only on $d,r_0,\Psi$ and $\ell$ such that for  $r\in (0,r_3]$ and $x\in D(r)$,
		\begin{align*}
			\P_{x} (X_{\tau_{D(r)}} \in D)  \asymp ( \delta_D(x) /r)^{\gamma(n_{Q_x})} \quad \text{ and } \quad 	\wh\P_{x} (\wh X_{\wh \tau_{D(r)}} \in D)  \asymp (\delta_D(x)/r)^{\wh \gamma(n_{Q_x})}.
		\end{align*}
	\end{cor}
	\pf Since the L\'evy exponent of the dual process $\wh Y$ is the complex conjugate of $\psi(\xi)$, 
it suffices to prove the statement for $X$ only, 		by interchanging the roles of $Y$ and $\wh Y$.

Let  $R\in (0,r_2]$, where 
$r_2$ is the constant in Proposition \ref{p:barriers-nonsymmetric}, and let $\sigma \in (0,1/4]$ be a  constant to be determined later. Define $h_{R, \sigma}$  and $\psi_{R,\sigma}$ as \eqref{e:def-h-r-sigma-nonsymmetric} and \eqref{e:def-psi-r-sigma-nonsymmetric}, respectively.   
Applying Proposition \ref{p:barriers-nonsymmetric}, we obtain for all $\eta \in (0,\sigma)$ and  $x\in D( \eta R)$,
\begin{align}\label{e:nonsymmetric-barrier-generators}
	|	\sL h_{R,  \sigma}(x)| &\le c_1 R^{-\Phi(x)} \delta_D(x)^{\Phi(x)-\alpha} \ell(\delta_D(x)/R) |\log (\delta_D(x)/R)| + c_2\sigma^{\gamma_0-\alpha} R^{-\alpha},\nn\\
		\sL \psi_{R, \sigma }(x) &\ge  \left( c_3  - \frac{c_4 L_\ell(\eta)}{\sigma^{\gamma_1-\gamma_0}\ell(\sigma)}  \right)   R^{-\Phi(x)} \delta_D(x)^{\Phi(x)-\alpha} \ell(\delta_D(x)/R)  |\log (\delta_D(x)/R)|.
\end{align}
Moreover, by  \eqref{e:regular-distance-1}, \eqref{e:inward-normal-result-1} and the scaling of $L_\ell$, we see that
$\psi_{R, \sigma}(y)/h_{R, \sigma}(y) \le c_5 L_\ell( \sigma)$ for all $y \in D$. By choosing $ \sigma$ small enough,  we obtain  $8c_1c_3^{-1}\psi_{R, \sigma}(y) \le h_{R,\sigma}(y)$ for all $y \in D$.

Note that $s^{\gamma_1-\alpha}\ell(s)$ is decreasing and $\lim_{s\to 0} s^{\gamma_1-\alpha}\ell(s)\ge \lim_{s\to 0} s^{(\gamma_1-\alpha)/2}\ell(1)=\infty$. Hence, there exists $\sigma_0 \in (0,\sigma)$ such that   $c_2\sigma^{\gamma_0-\alpha} \le c_1 \inf_{s\in (0,\sigma_0]} s^{\gamma_1-\alpha}\ell(s)$ and  $c_4L_\ell(\sigma_0) \le c_3\sigma^{\gamma_1-\gamma_0}\ell(\sigma)/2$.  By \eqref{e:nonsymmetric-barrier-generators}, it follows that  for all $x \in D(\sigma_0R)$,
\begin{align}\label{e:nonsymmetric-barrier-1}
	|	\sL h_{R,  \sigma}(x)| &\le c_1 R^{-\Phi(x)} \delta_D(x)^{\Phi(x)-\alpha} \ell(\delta_D(x)/R) |\log (\delta_D(x)/R)| + c_1R^{-\alpha}   \inf_{s\in (0,\sigma_0R]} (s/R)^{\gamma_1-\alpha}\ell(s/R) \nn\\
	&\le c_1 R^{-\Phi(x)} \delta_D(x)^{\Phi(x)-\alpha} \ell(\delta_D(x)/R) |\log (\delta_D(x)/R)| + c_1R^{-\gamma_1} \delta_D(x)^{\gamma_1-\alpha} \ell(\delta_D(x)/R)\nn\\
		&\le 2c_1 R^{-\Phi(x)} \delta_D(x)^{\Phi(x)-\alpha} \ell(\delta_D(x)/R) |\log (\delta_D(x)/R)|,\nn\\
			\sL \psi_{R, \sigma }(x) &\ge 2^{-1}c_3 R^{-\Phi(x)} \delta_D(x)^{\Phi(x)-\alpha} \ell(\delta_D(x)/R) |\log(\delta_D(x)/R)|.
\end{align}

Define $u_1(y):=h_{R,\sigma}(y) + 4c_1c_3^{-1}\psi_{R, \sigma}(y)$ and $u_2(y):=h_{R,\sigma}(y) -4c_1c_3^{-1}\psi_{R, \sigma}(y)$. By \eqref{e:regular-distance-1}, \eqref{e:exponent-Phi-1} and  \eqref{e:exponent-Phi-2}, we have for all   $x \in D(\sigma_0R)$, 
\begin{align*}
&\frac{	h_{R,\sigma}(x)}{(\delta_D(x)/R)^{\gamma(n_{Q_x})}} \vee \frac{(\delta_D(x)/R)^{\gamma(n_{Q_x})}}{	h_{R,\sigma}(x)}\\
& \le c_6 \bigg( \frac{\delta_D(x)}{R} \bigg)^{-|\Phi(x) - \Phi(Q_x)|}   \le c_6 \bigg( \frac{\delta_D(x)}{R} \bigg)^{-c_7\ell(\delta_D(x))}  \le c_6 \bigg( \frac{\delta_D(x)}{R} \bigg)^{-c_7\ell(\delta_D(x)/R)}.
\end{align*}
Since $\sup_{s \in (0,1/4]} s^{-\ell(s)}<\infty$ by Lemma \ref{l:double-Dini-limit}, it follows that
\begin{align}\label{e:barriers-nonsymmetric-first}
	u_1(x) \asymp u_2(x) \asymp h_{R,\sigma}(x) \asymp (\delta_D(x)/R)^{\gamma(n_{Q_x})}\quad \text{for $x \in D(\sigma_0R)$}.
\end{align}
On the other hand, by  \eqref{e:regular-distance-1} and  \eqref{e:nonsymmetric-barrier-1},  we have $c_8^{-1} \le u_2(y) \le u_1(y) \le c_8$ for all $y \in D\setminus D(\sigma_0 R)$, $\sL u_1 \ge 0$ in $D(\sigma_0R)$ and $\sL u_2 \le 0$ in $D(\sigma_0R)$.
Hence, by Corollary \ref{c:maximum-nonsymmetric}, we obtain for all $x\in D(\sigma_0R)$, 
\begin{align*}
	u_1(x) \le \E_x\big[ u_1(X_{\tau_{D(\sigma_0R)}}) \big] = \E_x\big[ u_1(X_{\tau_{D(\sigma_0R)}}): X_{\tau_{D(\sigma_0R)}} \in D \setminus D(\sigma_0R) \big] \le c_8 \P_x(X_{\tau_{D(\sigma_0r)}} \in D ) 
\end{align*}
and $	u_2(x) \ge \E_x[ u_2(X_{\tau_{D(\sigma_0R)}}) ] = \E_x[ u_2(X_{\tau_{D(\sigma_0R)}}): X_{\tau_{D(\sigma_0R)}} \in D \setminus D(a_0R) ] \ge c_8^{-1} \P_x(X_{\tau_{D(\sigma_0R)}} \in D ) .$ Combining these with \eqref{e:barriers-nonsymmetric-first},  we obtain the desired result with $r_3= \sigma_0r_2$. \qed

	 \medskip
	
	\noindent \textbf{Proof of Theorem \ref{t:nonsymmetricLevy}.}  
 Following the proof of Theorem \ref{t:critical} with Corollary \ref{c:critical-barrier} in place of Corollary \ref{c:nonsymmetric-barrier},  we obtain the desired results  from Theorems \ref{t:factorization-heatkernel}, \ref{t:largetime-bounded}, \ref{t:factorization-Green}, and Proposition \ref{p:boundary-term-scale}.   \qed

	\subsection{Killed stable-like processes
		 in $C^{1,\eps}$ open sets  with low regularity coefficients}\label{ss:5.3}
	
	Let $\alpha\in (0,2)$, 	$d\ge 2$ and $D\subset \R^d$ be a bounded  $C^{1,\eps}$ open set for some $\eps \in (0,1]$. Let $K$ be a symmetric Borel function on $\R^d\times \R^d$ satisfying 
	\begin{align}
	A_7 \le 	K(x,y) \le  A_8 &\quad \text{for all $x,y \in \R^d$},\label{e:low-regular-1}\\
	|K(x+h,y+h)-K(x,y)| \le 
	A_9 |h|^\theta &\quad \text{for all $x,y \in D$ and $|h|<1$},\label{e:low-regular-2}
	\end{align}
	for some constants 
	$A_8 \ge A_7>0$,  $A_9>0$  and $\theta>0$.
	Consider the bilinear form $(\overline \sE,  \sF)$ on  $L^2(\R^d,dx)$ defined by
	\begin{align*}
		\overline \sE(u,v) &= \int_{\R^d\times \R^d\setminus \diag} (u(x)-u(y))(v(x)-v(y)) J(x,y) dxdy, \\
			\sF&=\bigg\{ u \in L^2(\R^d,dx): \int_{\R^d\times \R^d\setminus \diag}  \frac{(u(x)-u(y))^2}{|x-y|^{d+\alpha}} dxdy<\infty\bigg\}.
	\end{align*}

By the arguments at the beginning of Subsection \ref{ss:critical-killing},   $(\overline \sE,\sF)$ is a regular symmetric Dirichlet form on $L^2(\R^d,dx)$. Further, by \cite[Theorem 1.1]{CK03}, the symmetric Hunt process $Z$ associated with $(\overline \sE,\sF)$ satisfies  hypothesis {\bf (A)} with $\phi(r)=r^\alpha$ and $T_0=1$.	Let  $Z^D$ denote the  killed  process of $Z$ on $D$. Compared to the setting of Subsection \ref{ss:critical-killing}, we now assume that  $D$ is  a bounded $C^{1,\eps}$ 
open set for some $\eps\in (0, 1]$
and consider only the killed process, while allowing the function  $K(x,y)$ to be less regular. Notably, if only \eqref{e:low-regular-1} and \eqref{e:low-regular-2} are imposed, the principal value integral
	$$
	p.v. \int_D (u(y)-u(x)) \frac{K(x,y)}{|x-y|^{d+\alpha}}dy
	$$ 
may fail to exist even for
 $u\in C_c^\infty(D)$. Under this setting, the following Green function estimates for $Z^D$ were recently established in \cite{KW24}.
 
 \begin{thm}\label{t:low-regular-green}
 	{\bf \!(\!\!\cite[Theorem 1.1]{KW24} \!\!).}	Let $\alpha \in (0,2)$, $d\ge 2$ and $D\subset \R^d$ be a  	bounded $C^{1,\eps}$ open set for some $\eps \in (0,1]$.	
	Assume \eqref{e:low-regular-1}  and \eqref{e:low-regular-2}.  Then the killed process  $Z^D$ admits a Green function $G(x,y)$ on $D\times D$  satisfying the following estimates: For  all $(x,y) \in D\times D$,
 		\begin{align*}
 			G(x,y) \asymp\bigg(1 \wedge\frac{\delta_D(x)}{|x-y|} \bigg)^{\alpha/2}\bigg(1 \wedge\frac{\delta_D(y)}{|x-y|} \bigg)^{\alpha/2}  |x-y|^{\alpha-d} .
 		\end{align*}
 \end{thm}

 Since every $C^{1,\eps}$ open set  is  $\eta$-fat, an application of Corollary \ref{c:factorization-equivalence} yields the following  heat kernel estimates for $Z^D$.
 
 \begin{thm}\label{t:low-regular-heatkernel}
 Under the setting of Theorem \ref{t:low-regular-green}, the process $Z^D$ admits a  transition density $p(t,x,y)$ on $(0,\infty) \times D \times D$ satisfying the following estimates: For  all $(t,x,y) \in (0,3]\times D\times D$,
 	\begin{align*}
 		p(t,x,y) \asymp \bigg(1 \wedge\frac{\delta_D(x)}{t^{1/\alpha}} \bigg)^{\alpha/2}\bigg(1 \wedge\frac{\delta_D(y)}{t^{1/\alpha}} \bigg)^{\alpha/2} \bigg( t^{-d/\alpha} \wedge \frac{t}{|x-y|^{d+\alpha}}\bigg).
 	\end{align*}
 Moreover, there exists $C\ge 1$ such that for all $(t,x,y) \in [3,\infty)\times D\times D$,
 	\begin{align*}
 		C^{-1}\delta_D(x)^{\alpha/2}\delta_D(y)^{\alpha/2}e^{-\lambda_1t}\le 		p(t,x,y) \le C\delta_D(x)^{\alpha/2}\delta_D(y)^{\alpha/2}e^{-\lambda_1t},
 	\end{align*}
	where $-\lambda_1:=\sup \text{\rm Re\,} \sigma(\sA)$
	and $\sigma(\sA)$ is the spectrum of the generator $\sA$ of $Z^D$. 
 \end{thm}

	\appendix
	
	\section{Appendix: Regularization of Dini and double Dini functions}\label{appendix}
	
	Recall that $\Dini$ denotes the family of all Dini functions, and $\DDini$ the family of all double Dini functions. Clearly, $\DDini\subset \Dini$.
		We extend every $\ell \in \Dini$ to $(0,\infty)$ by setting $\ell(r)=\ell(1)$ for $r\ge1$.  		For $f:(0, 1]\to (0, \infty)$ and $\eps\in (0,1)$, define 
		\begin{align}\label{e:def-f-eps} f_\eps (r):=\inf_{s\in (0,1]}\left\{ f(s)  + f(1) \bigg( \frac{r}{s}\bigg)^\eps  \right\}, \quad r\in (0,1]. 
		\end{align}

	\begin{lemma}\label{l:Dini-regularizing-pre}
		If $f\in \Dini$, then for any $\eps\in (0,1)$,
		$f_\eps(r)/r^\eps$ is non-increasing, $f(r)\le f_\eps(r)$ for all $r\in (0,1]$, and $f_\eps \in \Dini$. Moreover, if $f\in \DDini$, then $f_\eps \in \DDini$.
	\end{lemma}
	\pf  Clearly, $f_\eps(r)$ is increasing and  $f_\eps(r)/r^\eps$ is non-increasing.	For all $s,r \in (0,1]$, we have $f(r)\le f(s)$ if $r\le s$ and $f(r)\le f(1) \le f(1)(r/s)^\eps $ if $r\ge s$. Thus, $f(r) \le f_\eps(r)$ for all $r\in (0,1]$. Moreover, using $ f_\eps (r) \le f(r |\log r|^{3/\eps})  + f(1) / |\log r|^3$ for all  sufficiently small $r$ and  the change of the variables $s= r|\log r|^{3/\eps}$, we obtain
	\begin{align*}
		\int_0^1 \frac{ f_\eps(r)}{r} dr \le  \int_0^{c_1} \frac{f(r|\log r|^{3/\eps})}{r} dr + \int_0^{c_1} \frac{f(1)}{r|\log r|^3} dr + c_2 \le c_3 \int_0^{1} \frac{f(s)}{s}ds + c_2<\infty.
	\end{align*}
	Hence, $f_\eps \in \Dini$. Further, we also have
	\begin{align*}
			\int_0^1 \frac{ f_\eps(r) |\log r|}{r} dr \le  \int_0^{c_1} \frac{f(r|\log r|^{3/\eps}) |\log r|}{r} dr + \int_0^{c_1} \frac{f(1)}{r|\log r|^2} dr + c_2 \le c_4 \int_0^{1} \frac{f(s)|\log s|}{s}ds + c_5.
	\end{align*}
	Therefore, $f\in \DDini$ implies $f_\eps\in \DDini$.	  The proof is complete. \qed

	\begin{lemma}\label{l:Dini-regularizing}
	Let $f\in \Dini$ (resp. $f\in \DDini$). For every $\eps \in (0,1/4]$, there exists  $\ell \in \Dini \cap  C^2((0,\infty))$ (resp.  $\ell \in \DDini \cap  C^2((0,\infty))$) satisfying the following properties:

	 (a) $\ell(1)\le 4f(1)$ and  $f(r)\le \ell(r)$   for all $r\in (0,1]$. 
	 
	 (b) $ r\ell'(r) \le 2\eps\ell(r)$ and $|r^2\ell''(r)| \le 6\eps \ell(r)$ for all $r\in (0,1]$.
	 
	 (c) $\ell(r)/r^\eps$ is non-increasing on $(0,1]$.
	 
	 (d) $\ell(r) \le \eps \int_0^r s^{-1}\ell(s)ds$ for all $r\in (0,1]$.
	\end{lemma}
	\pf Define $f_\eps$ by \eqref{e:def-f-eps} and
	\begin{align*}
		\ell(r):=\frac{2}{r}\int_{0}^{r}\frac{1}{u} \int_0^u \frac{1}{t}\int_0^t  f_\eps(s)\,dsdtdu, \quad r\in (0,1].
	\end{align*}
	   By taking $s=r=1$ in \eqref{e:def-f-eps}, we get $f_\eps(1) \le 2f(1)$. Hence,  $\ell(1) \le 4f(1)$.
	Since  $f_\eps(s)$ is increasing and  $f_\eps(s)/s^\eps$ is non-increasing,  	we have
	\begin{align}\label{e:Dini-regularizing-1}
	f_\eps(t) \ge 	\frac{1}{t}	\int_0^t f_\eps (s) ds \ge \frac{f_\eps(t)}{t^{1+\eps}} \int_0^t s^\eps ds =  \frac{f_\eps(t)}{1+\eps}, \quad t\in (0,1].
	\end{align}
	Using \eqref{e:Dini-regularizing-1} three times, we obtain
	\begin{align*}
		\ell(r) \ge \frac{2}{(1+\eps)r}\int_{0}^{r}\frac{2}{u} \int_0^u   f_\eps(t)\,dtdu \ge \frac{1}{(1+\eps)^2 r}\int_0^r f_\eps(u) du \ge \frac{2 f_\eps(r)}{(1+\eps)^3} \ge f_\eps(r) \ge f(r). 
	\end{align*}
	Moreover, using the monotonicity of $f_\eps$ and  \eqref{e:Dini-regularizing-1}, we also get that
	\begin{align*}
		r\ell'(r) &= \frac{2}{r} \int_0^r \frac{1}{t}\int_0^t  f_\eps(s)\,dsdt - \frac{2}{r}\int_{0}^{r}\frac{1}{u} \int_0^u \frac{1}{t}\int_0^t  f_\eps(s)\,dsdtdu\\
		&\ge \frac{2}{r} \int_0^r \frac{1}{t}\int_0^t  f_\eps(s)\,dsdt - \frac{2}{r}\int_{0}^{r}\frac{1}{u} \int_0^u  f_\eps(t)\,dtdu=0,\\
			r\ell'(r) &\le  \frac{2}{r} \int_0^r \frac{1}{t}\int_0^t  f_\eps(s)\,dsdt - \frac{2}{(1+\eps )r}\int_{0}^{r}\frac{1}{u} \int_0^u f_\eps(t) \, dtdu \\
			&= \frac{2\eps }{(1+\eps)r}\int_0^r \frac{1}{t}\int_0^t  f_\eps(s)\,dsdt \le \frac{2\eps f_\eps(r)}{1+\eps}  \le 2\eps \ell(r)
	\end{align*}
	and
	\begin{align*}
		|r^2 \ell''(r)| &\le \frac{2}{r} \left|  \int_0^r f_\eps (s) ds  - \int_0^r \frac{1}{t}\int_0^t f_\eps (s)ds dt  \right|\\
		&\quad +  \frac{4}{r} \left| \int_0^r \frac{1}{u}\int_0^u \frac{1}{t} \int_0^t f_\eps (s)dsdtdu - \int_0^r \frac{1}{t}\int_0^t f_\eps(s) ds dt \right| \\
	&\le \frac{2\eps }{(1+\eps)r} \int_0^r f_\eps(s)ds + \frac{4\eps }{(1+\eps)r}\int_0^r \frac{1}{t}\int_0^t f_\eps(s) ds dt  \le \frac{6\eps f_\eps(r)}{1+\eps} \le 6\eps \ell(r).
	\end{align*}
	 Thus, $\ell$ is non-increasing, and (a) and (b) hold. 
	 Since  $f_\eps(s)$ is increasing, $\ell(r) \le 2f_\eps(r)$. 	 It follows from Lemma \ref{l:Dini-regularizing-pre} that $f\in \Dini$ (resp. $f\in \DDini$) implies $f_\eps \in \Dini$  (resp. $f_\eps\in \DDini$). Thus $\ell \in \Dini$ if $f\in \Dini$ and $\ell \in \DDini$ if $f\in \DDini$.

For all $0<r\le R\le 1$, using  a change of the variables and the monotonicity of $f_\eps(u)/u^\eps$, we get
\begin{align*}
	&\frac{\ell(R)}{R^\eps } =\frac{2}{R^{1+\eps}}\int_{0}^{R}\frac{1}{u} \int_0^u \frac{1}{t}\int_0^t  f_\eps(s)\,dsdtdu =\frac{2}{R^{1+\eps}}\int_{0}^{r}\frac{1}{\wt u} \int_0^{R\wt u /r} \frac{1}{t}\int_0^t  f_\eps(s)\,dsdtd\wt u \\
	&=\frac{2}{R^{1+\eps}}\int_{0}^{r}\frac{1}{\wt u} \int_0^{\wt u } \frac{1}{\wt t}\int_0^{R \wt t/r}  f_\eps(s)\,dsd\wt td\wt u =\frac{2}{R^{\eps}r}\int_{0}^{r}\frac{1}{\wt u} \int_0^{\wt u } \frac{1}{\wt t}\int_0^{\wt t}  f_\eps(R \wt s/r)\,d\wt sd\wt td\wt u \le \frac{\ell(r)}{r^\eps}.
\end{align*}
This proves (c). For  (d), by  (c),  we have for all $r\in (0,1]$,
\begin{align*}
	\int_0^r \frac{\ell(s)}{s} ds \ge \frac{\ell(r)}{r^\eps}\int_0^r \frac{ds}{s^{1-\eps}}  =  \frac{\ell(r)}{\eps}.
\end{align*}	
	The proof is complete. \qed 
	
	\begin{lemma}\label{l:double-Dini-limit}
	If  $f\in \DDini$, then $\lim_{r\to 0}  f(r)|\log r| =0$. 
	\end{lemma}
	\pf Suppose $\limsup_{r\to 0} f(r)|\log r|\ge c_1>0$. Then there exists a sequence $(a_n)_{n\ge 1} \subset (4,\infty)$ such that $a_{n+1}\ge 2a_{n}$ and $f(1/a_{n}) \log a_{n} \ge c_1/2$ for all $n\ge 1$. It follows that
	\begin{align*} 
		\int_0^1 \frac{f(r)|\log r| }{r}dr \ge \sum_{n=1}^\infty \frac{f(1/a_n)\log (a_n/2)}{2/a_n} \int_{1/a_n}^{2/a_n}dr \ge \sum_{n=1}^\infty \frac{c_1}{8}  = \infty,
	\end{align*}
	which is a contradiction.  \qed

		\vskip 0.4truein

\noindent {\bf Soobin Cho:} Department of Mathematics,
University of Illinois Urbana-Champaign,
Urbana, IL 61801, U.S.A.
Email: \texttt{soobinc@illinois.edu}
	
	\medskip

\noindent {\bf Renming Song:} Department of Mathematics,
University of Illinois Urbana-Champaign,
Urbana, IL 61801, U.S.A.
Email: \texttt{rsong@illinois.edu}

\end{document}